\documentclass[11pt]{article}

\usepackage[margin=1.5in]{geometry}

\usepackage[english]{babel}
\usepackage{amsmath,amssymb,amsthm}
\usepackage{enumerate}
\usepackage[affil-it]{authblk}
\usepackage{color}
\usepackage{hyperref}
\usepackage{comment}
\usepackage{bbm}

\usepackage{pgf}
\usepackage{tikz}
\usetikzlibrary{arrows,automata}
\usepackage[latin1]{inputenc}
\usepackage{verbatim}
\usepackage{authblk}
\usepackage{blkarray, bigstrut}
\usepackage{booktabs}

\newtheorem{theorem}{Theorem}[section]
\newtheorem{proposition}[theorem]{Proposition}
\newtheorem{lemma}[theorem]{Lemma}
\newtheorem{corollary}[theorem]{Corollary}

\theoremstyle{definition}
\newtheorem{definition}[theorem]{Definition}
\newtheorem{example}[theorem]{Example}
\newtheorem*{remark}{Remark}

\newcommand{\ZZ}{\mathbb{Z}}

\newcommand{\R}{\mathbb{R}}
\newcommand{\C}{\mathbb{C}}

\newcommand{\mx}{\mathbf{x}}

\newcommand{\me}{\mathbf{e}}

\newcommand{\ma}{\mathbf{a}}
\newcommand{\mv}{\mathbf{v}}

\newcommand{\FF}{\mathbb{F}}
\def\cL{\mathcal{L}}
\def\cF{\mathcal{F}}
\def\cB{\mathcal{B}}
\def\id{\mathrm{id}}
\def\CC{\mathbb{C}}
\renewcommand\O{\mathcal{O}}
\def\cC{\mathcal{C}}

\renewcommand{\P}{\operatorname{\mathbb{P}}}

\newcommand{\T}{\mathcal{T}}

\newcommand{\lB}{\mathbbm{1}_\cB}
\newcommand{\cX}{\mathcal{X}}

\def\Card{\mathrm{Card}}
\def\End{\mathrm{End}}
\def\cC{\mathcal{C}}
\def\Ind{\mathrm{Ind}}

\begin{document}

\title{Double Coset Markov Chains  \\ \small{Expository Version} }

\author[1,2]{Persi Diaconis}
\author[3]{Arun Ram}
\author[2]{Mackenzie Simper \thanks{Correspondence: \href{mailto:msimper@stanford.edu}{\nolinkurl{msimper@stanford.edu}}}}
\affil[1]{Department of Statistics, Stanford University}
\affil[2]{Department of Mathematics, Stanford University}
\affil[3]{School of Mathematics and Statistics, University of Melbourne}

\date{}

\maketitle

\begin{abstract}
Let $G$ be a finite group. Let $H, K$ be subgroups of $G$ and $H \backslash G / K$ the double coset space. Let $Q$ be a probability on $G$ which is constant on conjugacy classes ($Q(s^{-1} t s) = Q(t)$). The random walk driven by $Q$ on $G$ projects to a Markov chain on $H \backslash G /K$. This allows analysis of the lumped chain using the representation theory of $G$. Examples include coagulation-fragmentation processes and natural Markov chains on contingency tables. Our main example projects the random transvections walk on $GL_n(q)$ onto a Markov chain on $S_n$ via the Bruhat decomposition. The chain on $S_n$ has a Mallows stationary distribution and interesting mixing time behavior. The projection illuminates the combinatorics of Gaussian elimination. Along the way, we give a representation of the sum of transvections in the Hecke algebra of double cosets. Some extensions and examples of double coset Markov chains with $G$ a compact group are discussed.
\end{abstract}

\tableofcontents

\section{Introduction} \label{sec: intro}

The main results of this paper develop tools which allow projecting a random walk on a group to a Markov chain on special equivalence classes of the group. Then Fourier analysis on the group can be harnessed to give sharp analysis of rates of convergence to stationary for the Markov chain on equivalence classes.

\begin{example}[Coagulation-Fragmentation Processes]
In chemistry and physics, coagulation-fragmentation processes are models used to capture the behavior of `blobs' that combine and break up over time. These processes are used in population genetics to model the merging and splitting of family groups. A simple mean field model considers $n$ unlabeled particles in a partition $\lambda = (\lambda_1, \dots, \lambda_k)$, $\lambda_1 \ge \lambda_2 \ge \hdots \ge \lambda_k > 0$, $\sum_{i = 1}^k \lambda_i = n$. At each step of the process, a pair of particles is chosen uniformly at random and the partition evolves according to the rules:
\begin{enumerate}
    \item If the particles are in distinct blocks, combine the blocks.
    \item If the particles are in the same block, break the block uniformly into two blocks.
    \item If the same particle is chosen twice, do nothing.
\end{enumerate}
This defines a Markov chain on partitions of $n$. Natural questions are:
\begin{itemize}
    \item What is the stationary distribution $\pi(\lambda)$?
    \item How does the process evolve?
    \item How long to reach stationarity?
\end{itemize}
  All of this is available by considering the random transpositions process on the symmetric group $S_n$. The transition probabilities for this process are constant on conjugacy classes; the conjugacy classes are indexed by partitions; and the conjugacy class containing the current permutation of the walk evolves as the coagulation-fragmentation process on partitions. The answers are:
\begin{itemize}
    \item The stationary distribution is $\pi(\lambda) = \prod_{i=1}^n 1/(i^{a_i} a_i!)$ for a partition $\lambda$ with $a_i$ parts of size $i$.
    
    \item Starting at $\lambda = 1^n$, the pieces evolve as the connected components of a growing Erd\H{o}s-R\'{e}nyi random graph. 
    
    \item It takes order $n \log(n)$ steps to reach stationarity.
\end{itemize}

\end{example}

The coagulation-fragmentation process is a special case of a double coset walk: Let $G$ be a finite group, $H, K$ subgroups of $G$. The equivalence relation 
\[
s \sim t \iff hsk = t, \quad h \in H, \quad k \in K
\]
partitions the group into double cosets $H \backslash G / K$. Let $Q(s)$ be a probability on $G$. That is, $0 \le Q(s) \le 1$ for all $s \in G$ and $\sum_{s \in G} Q(s) = 1$. Further assume $Q(s)$ a \emph{class function}: it is constant on conjugacy classes, i.e.\ $Q(s) = Q(t^{-1} s t)$ for all $s, t \in G$. The probability $Q$ defines a random walk on $G$ by multiplication by random elements chosen according to $Q$. In other words, the random walk is induced by convolution, $Q^{*k}(s) = \sum_{t \in G} Q(t) Q^{*(k-1)}(s t^{-1})$ and a single transition step has probability $P(x, y) = Q(yx^{-1})$.

This random walk induces a random process on the space of double cosets. While usually a function of a Markov chain is no longer a Markov chain, in this situation the image of the random walk on $H \backslash G / K$ is Markov. Section \ref{sec: MC} proves the following general result. Throughout, we pick double coset representatives $x \in G$ and write $x$ for $H x K$.

\begin{theorem} \label{thm: general}
For $Q(s) = Q(t^{-1} s t)$ a probability on $G$, the induced process on $H \backslash G / K$ is Markov with the following properties.
\begin{enumerate}[(1)]
    \item The transitions are
    \[
    P(x, y) = Q(HyKx^{-1}), \quad x, y \in H \backslash G /K.
    \]
    \item The stationary distribution is 
    \[
    \pi(x) = \frac{|HxK|}{|G|}.
    \]
    \item If $Q(s^{-1}) = Q(s)$, then $P$ is reversible with respect to $\pi$:
    \[
    \pi(x)P(x, y) = \pi(y)P(y, x).
    \]
    
    \item Suppose $Q$ is concentrated on a single conjugacy class $\mathcal{C}$ (that is, $Q(s) = \delta_{\mathcal{C}}(s)/|\mathcal{C}|$). Then the eigenvalues of $P(x, y)$ are among the set 
    \[
    \left\lbrace \frac{\chi_\lambda(\mathcal{C})}{\chi_\lambda(1)} \right\rbrace_{\lambda \in \widehat{G}},
    \]
    where $\widehat{G}$ is the set of all irreducible representations of $G$ and $\chi_\lambda$ is the character of the irreducible representation indexed by $\lambda$.
    
    \item The multiplicity of $\chi_\lambda(\mathcal{C})/\chi_\lambda(1)$ is
    \[
    m_\lambda = \left\langle \chi_\lambda|_H, 1 \right\rangle \cdot \left\langle \chi_\lambda|_K, 1 \right\rangle.
    \]
    Here $\left\langle \chi_\lambda|_H, 1 \right\rangle$ is the number of times the trivial representation appears when $\chi_\lambda$ is restricted to $H$.
    
    \item For any time $\ell > 0$,
    \[
    \sum_{x \in H \backslash G /K} \pi(x) \|P_x^\ell - \pi \|_{TV}^2 \le \frac{1}{4} \sum_{\lambda \neq 1} m_\lambda \left| \frac{\chi_\lambda(\mathcal{C})}{\chi_\lambda(1)} \right|^{2 \ell}.
    \]
    \end{enumerate}
\end{theorem}

This theorem shows that the properties of the induced chain are available via the character theory of $G$. It is proved with variations and extensions in Section \ref{sec: MC}. The main example is introduced next.

\begin{example}[$GL_n(q)$ and Gaussian Elimination] \label{ex: gln}
Fix a prime power $q$ and let $GL_n(q)$ be the invertible $n \times n$ matrices over $\mathbb{F}_q$. Let $H = K = \cB$ be the Borel subgroup: upper-triangular matrices in $GL_n(q)$. A classical result is the Bruhat decomposition,
\[
GL_n(q) = \bigsqcup_{\omega \in S_n} \cB \omega \cB,
\]
where $\omega$ is the permutation matrix for the permutation $\omega \in S_n$. This means the double cosets $\cB \backslash GL_n(q) / \cB$ are indexed by permutations. As explained in Section \ref{sec: gaussianElim} below, the permutation $\omega$ associated to $M \in GL_n(q)$ is the `pivotal' permutation when $M$ is reduced to upper-triangular form by row reduction (Gaussian elimination).

\end{example}

The set of \emph{transvections} $\mathcal{T}_{n,q}$ is the conjugacy class containing the basic row operations $I + \theta E_{ij}$; here $\theta \in \FF_q$, and $E_{ij}$ is the matrix with a $1$ in position $(i, j)$ and zeroes everywhere else (so $I + \theta E_{ij}$ acts by adding $\theta$ times row $i$ to row $j$). Hildebrand \cite{hildebrand} gave sharp convergence results for the Markov chain on $GL_n(q)$ generated by the class function $Q(M) = \delta_{\mathcal{T}}(M)/|\mathcal{T}_{n,q}|$. He shows that $n$ steps are necessary and sufficient for convergence to the uniform distribution, for any $q$. Of course, convergence of the lumped chain on $\cB \backslash GL_n(q) / \cB$ might be faster. The results we found surprised us. Careful statements are included below. At a high level, we found:
\begin{itemize}
    \item Starting from a `typical' state $x \in S_n$, order $\log(n)/\log(q)$ steps are necessary and sufficient for convergence. This is an exponential speed-up from the original chain.
    \item Starting from $id$, order $n$ steps are necessary and sufficient for convergence.
    \item Starting from $\omega_0$, the reversal permutation, order $\log(n)/(2\log(q))$ steps are necessary and sufficient for convergence.
\end{itemize}
To simplify the statement of an honest theorem, let us measure convergence in the usual chi-square or $L^2$ distance:
\[
\chi^2_x(\ell) = \sum_y \frac{\left( P^\ell(x, y) - \pi_q(y) \right)^2}{\pi_q(y)}.
\]
In Section \ref{sec: MT_results}, the usual total variation distance is treated as well.

\begin{theorem} \label{thm: transvections}
The random transvections walk on $GL_n(q)$ induces a Markov chain $P(x, y)$ on $S_n \cong \cB \backslash GL_n(q) / \cB$ with stationary distribution.
    \[
    \pi_q(\omega) = q^{I(\omega)}/[n]_q!, \quad [n]_q! := (1 + q)(1 + q + q^2) \hdots (1 + q + \hdots + q^{n-1}),
    \]
    for $\omega \in S_n$, where $I(\omega)$ is the number of inversions in $\omega$.
    
Furthermore, if $\log(q) > 6/n$ then the following statements are true. 
\begin{enumerate}[(a)]
    
    \item (Typical Start) If $\ell \ge (\log(n) + c)/(\log(q) - 6/n)$ for any $c > 0$, then
    \[
    \sum_{x \in S_n} \pi(x) \chi_x^2(\ell) \le (e^{e^{-c}} - 1) + e^{-2c}.
    \] 
    Conversely, for any $\ell$
    \[
    \sum_{x \in S_n} \pi(x) \chi_x^2(\ell)  \ge (n-1)^2 q^{-4 \ell}.
    \]
    These results show order $\log_q(n)$ steps are necessary and sufficient for convergence.
    
    \item (Starting from $id$) If $\ell \ge (n \log(q)/2 + c)/(\log(q) - 6/n), c > 0$, then
    \begin{align*}
        \chi_{id}^2(\ell) \le (e^{e^{-2c}} - 1) + e^{-2c}.
    \end{align*}
    Conversely, for any $\ell$,
    \[
        \chi_{id}^2(\ell) \ge (q^{n-1} - 1)(n - 1)q^{-4 \ell}.
    \]
    These results show that order $n$ steps are necessary and sufficient for convergence starting from the identity. 
    
    \item (Starting from $\omega_0$) If $\ell \ge (\log(n)/2 + c)/(\log(q) - 6/n)$ for $c \ge 2 \sqrt{2}$, then there is a universal constant $K > 0$ (independent of $q, n$) such that 
    \begin{align*}
        \chi_{\omega_0}^2(\ell) \le - 2K \log(1 - e^{-c}) + K e^{-nc}.
    \end{align*}
    Conversely, for any $\ell$,
    \[
        \chi_{\omega_0}^2(\ell) \ge q^{-(n-2)}(n-1) (q^{n-1} -1) q^{-4 \ell}.
    \]
    These results show that order $\log_q(n)/2$ steps are necessary and sufficient for convergence starting from $\omega_0$.

\end{enumerate}
\end{theorem}

\begin{remark}
Note that while Hildebrand's result of order $n$ convergence rate was independent of $n$, the rates in Theorem \ref{thm: transvections} depend on $q$.
\end{remark}

The stationary distribution $\pi_q$ is the Mallows measure on $S_n$. This measure has a large enumerative literature; see \cite{diaconisSimper} Section 3 for a review or \cite{chenyangThesis}. It is natural to ask what the induced chain `looks like' on $S_n$. After all, the chain induced by random transpositions on partitions has a simple description and is of general interest. Is there a similarly simple description of the chain in $S_n$? This question is treated in Section \ref{sec: heckeComputations} using the language of Hecke algebras.

\begin{theorem} \label{thm: D}
Let $H_n(q)$ be the Hecke algebra corresponding to the $\cB \backslash GL_n(q) / \cB$ double cosets and $D = \sum_{T \in \mathcal{T}_{n,q}} T \in H_n(q)$ be the sum of all transvections. Then,
\[
D = (n - 1)q^{n-1} - [n-1]_q + (q-1) \sum_{1 \le i < j \le n} q^{n - 1 - (j-i)} T_{ij},
\]
with $T_{ij}$ in the Hecke algebra. 
\end{theorem}

This gives a probabilistic description of the induced chain on $S_n$. Roughly stated, from $\omega \in S_n$ pick $(i,j), j < i,$ with probability proportional to $q^{-(j - i)}$ and transpose $i$ and $j$ in $\omega$ using the Metropolis algorithm. This description is explained in Section \ref{sec: heckeComputations}; see also \cite{diaconisRam}. The probabilistic description is crucial in obtaining good total variation lower bounds for Theorem \ref{thm: transvections}.

\paragraph{Outline} 
 Section \ref{sec: background} develops and surveys background material on double cosets, Markov chains (proving Theorem \ref{thm: general}), transpositions and coagulation-fragmentation processes, transvections, and Gaussian elimination. Theorem \ref{thm: transvections} is proven in Section  \ref{sec: MT_results}. Theorem \ref{thm: D} is proved in Section \ref{sec: heckeComputations} using a row reduction; an alternative proof is contained Section \ref{sec: symmetricFunction} using symmetric function theory. Section \ref{sec: topCard} contains another Markov chain from a lumping of the transvections chain, and Section \ref{sec: extensions} surveys further examples -- contingency tables and extensions of the $GL_n$ results to finite groups of Lie type, for which the Bruhat decomposition holds and there are natural analogs of transvections. Of course, there are an infinite variety of groups $G, H, K$, and we also indicate extensions to compact groups.

\paragraph{Notation}
Throughout, $q$ will be a prime power. For a positive integer $n$, define the quantities
\[
[n]_q := \frac{q^n - 1}{q -1} = q^{n-1} + q^{n-2} + \hdots + q + 1, \quad [n]_q! := [n]_q [n-1]_q \hdots [1]_q.
\]
Sometimes the subscript $q$ will be omitted when it is obvious from context.

\paragraph{Acknowledgements}

The authors thank Sourav Chatterjee, David Craven, Jason Fulman, Bob Guralnick, Jimmy He, Marty Isaacs, Zhihan Li, Martin Liebeck, James Parkinson, and Nat Theme and for helpful discussions. M.S.\ was supported by a Lieberman Fellowship at Stanford University. P.D.\ was partially supported by NSF grant DMS-1954042.

\section{Background} \label{sec: background}

This section gives the basic definitions and tools needed to prove our main results. Section \ref{sec: doubleCosets} gives background on double cosets. In Section \ref{sec: MC}, Markov chains are reviewed and Theorem \ref{thm: general} is proved, along with extensions. Section \ref{sec: RT} reviews the coagulation-fragmentation literature along with the random transpositions literature. Section \ref{sec: transvections} develops what we need about transvections and Section \ref{sec: gaussianElim} connects the Bruhat decomposition to Gaussian elimination.

\subsection{Double Cosets} \label{sec: doubleCosets}

Let $H, K$ be subgroups of a finite group $G$. The double coset decomposition is a standard tool of elementary group theory with applications in representation theroy and number theory (Hecke algebras). For a detailed survey, see \cite{diaconisSimper} or \cite{curtisReiner}. Double cosets can have very different sizes and \cite{diaconisSimper}, \cite{paguyo} develop a probabilistic and enumerative theory. For present applications, an explicit description of the double cosets are needed.

\begin{example}
Let $S_\lambda, S_\mu$ be parabolic subgroups of the symmetric group $S_n$. Here $\lambda = (\lambda_1, \dots, \lambda_I)$ and $\mu = (\mu_1, \dots, \mu_J)$ are partitions of $n$. The subgroup $S_\lambda$ consists of all permutations in which the first $\lambda_1$ elements may only be permuted amongst each other, the next $\lambda_1 + 1, \dots, \lambda_1 + \lambda_2$ elements may only be permuted amongst each other, and so on. It is a classical fact that the double cosets $S_\lambda \backslash S_n / S_\mu$ are in bijection with `contingency tables' -- arrays of non-negative integers with row sums $\lambda$ and column sums $\mu$. See \cite{jamesKerber}, Section 1.3. For proofs and much discussion of the connections between the group theory and applications and statistics see \cite{diaconisSimper}, Section 5. Random transpositions on $S_n$ induces a natural Markov chain on these tables, see Chapter 3 of \cite{simperThesis}.
\end{example}

\begin{example}
Let $M$ be a finite group, $G = M \times M$, and $H = K = M$ embedded diagonally as subgroups of $G$ (that is, $\{(m, m): m \in M \}$). The conjugacy classes in $G$ are products of classes in $M$. In the double coset equivalency classes, note that
\[
(s, t) \sim (id, s^{-1}t) \sim (id, k^{-1} s^{-1} t k),
\]
and so double cosets can be indexed by conjugacy classes of $M$. If $Q_1$ is a conjugacy invariant probability on $M$, then $Q = Q_1 \times \delta_{id}$ is conjugacy invariant on $G$. The random walk on $G$ induced by $Q$ maps to the random walk on $M$ induced by $Q_1$. In this way, the double coset walks extend conjugacy invariant walks on $M$. Example 1 in the introduction is a special case.

Of course, the conjugacy classes in $M$ (and so the double cosets) can be difficult to describe. Describing the conjugacy classes of $U_n(q)$ -- the unit upper-triangular matrices in $GL_n(q)$ -- is a well known `wild' problem. See \cite{20author} or \cite{diaconisMalliaris} for background and details.
\end{example}

\begin{example}
Let $G$ be a finite group of Lie type, defined over $\FF_q$, with Weyl group $W$. Let $B$ be the Borel subgroup (maximal solvable subgroup). Take $H = K = B$. The Bruhat decomposition gives
\[
G = \bigsqcup_{\omega \in W} B \omega B,
\]
so the double cosets are indexed by $W$. See \cite{carterBook}, Chapter 8, for a clear development in the language of groups with a $(B, N)$ pair. 

Conjugacy invariant walks on $G$ have been carefully studied in a series of papers by David Gluck, Bob Guralnick, Michael Larsen, Martin Liebeck, Aner Shalev, Pham Tiep and others. These authors develop good bounds on the character ratios needed. See \cite{guralnickLarsen} for a recent paper with careful reference to earlier work. Of course, Example 2 with $G = GL_n(q), W = S_n$ is a special case. The present paper shows what additional work is needed to transfer results from $G$ to $W$.
\end{example}

The double cosets form a basis for the algebra of $H-K$ bi-invariant functions $L(H \backslash G /K)$ with product
\[
f \star g(s) = \sum_t f(t) g(st^{-1}).
\]
This is usually developed for $H = K$ \cite{curtisReiner}, \cite{ceccherini}, \cite{diaconisBook} but the extra flexibility is useful. We add a caveat: When $H = K$, the algebra of bi-invariant functions (into $\CC$) is semi-simple and with a unit. This need not be the case for general $H$ and $K$. David Craven tried all pairs of subgroups of $S_4$ and found examples which were not semi-simple. For instance, $H = D_8$ (generated by $\{(34), (14)(23), (13)(24) \}$) and $K = C_2 \times C_2 \subset S_4$ (generated by $\{(14)(23), (13)(24) \}$) give an algebra with no unit which is not semi-simple. This occurs even for some pairs of distinct parabolic subgroups of $S_n$. There are also distinct pairs of parabolics where the algebra is semi-simple. Determining when this occurs is an open question. 

Further examples are in Section \ref{sec: extensions}. Since the theory is developed for general $H, K, G$ there is a large set of possibilities. What is needed are examples where the double cosets are indexed by familiar combinatorial objects and the walks induced on $H \backslash G / K$ are of independent interest. 


\subsection{Markov Chain Theory} \label{sec: MC}

Let $H, K$ be subgroups of a finite group $G$, and $Q$ a probability on $G$. See \cite{levinPeres} for an introduction to Markov chains; see \cite{diaconisBook} or \cite{saloffcoste2004} for random walks on groups.

\begin{proposition} \label{prop: inducedChain}
Let $Q$ be a probability on $G$ which is $H$-conjugacy invariant ($Q(s) = Q(h^{-1}s h) $ for $h \in H, s \in G$). The image of the random walk driven by $Q$ on $G$ maps to a Markov chain on $H \backslash G /K$ with transition kernel
\[
P(x, y) = Q(HyKx^{-1}).
\]
The stationary distribution of $P$ is $\pi(x) = |HxK|/|G|$. If $Q(s) = Q(s^{-1})$ then $(P, \pi)$ is reversible.
\end{proposition}

\begin{proof}
The kernel $P$ is well-defined; that is, it is independent of the choice of double coset representatives for $x, y$. Dynkin's criteria (\cite{kemenySnell} Chapter 6, \cite{pang}) says that the image of a Markov chain in a partitioning of the state space is Markov if and only if for any set in the partition and any point in a second set, the chance of the original chain moving from the point to the first set is constant for points in the second set. 

Fixing $x, y$, observe
\[
Q(HyK(hxk)^{-1}) = Q(HyK x^{-1}h^{-1}) = Q(HyKx^{-1}).
\]
If the uniform distribution on $G$ is stationary for the walk generated by $Q$, then the stationary distribution of the lumped chain is $\pi(x) = |HxK|/|G|$. Finally, any function of a reversible chain is reversible and $Q(s) = Q(s^{-1})$ gives reversibility of the walk on $G$.
\end{proof}

\begin{remark}
A different sufficient condition for Proposition \ref{prop: inducedChain} is $Q(sh) = Q(s)$ for all $s \in G, h \ in H$. 
\end{remark}

\begin{remark}
Usually, a function of a Markov chain is \emph{not} Markov. For relevant discussion of similar `orbit chains', see \cite{diaconisBoyd}.
\end{remark}

In all of our examples, the measure $Q$ is a class function ($Q(s) = Q(t^{-1}st)$ for all $s, t, \in G$), which is a stronger requirement than that in Proposition \ref{prop: inducedChain}. The eigenvalues of the walk on $G$ can be given in terms of the irreducible complex characters of $G$. Let $\widehat{G}$ be an index set for these characters. We write $\lambda \in \widehat{G}$ and $\chi_\lambda(\cC)$ for the character value at the conjugacy class $\cC$. Let
\[
\beta_\lambda = \frac{1}{\chi_\lambda(1)} \sum_{s \in G} Q(s) \chi_\rho(s).
\]
If $Q$ is simply concentrated on a single conjugacy class $\cC$, then $\beta_\lambda$ is the character ratio
\[
\beta_\lambda = \frac{\chi_\lambda(\cC)}{\chi_\lambda(1)}.
\]
For a review of a large relevant literature on character ratios and their applications, see \cite{guralnickLarsen}.

The restriction of $\chi_\lambda$ to $H$ is written $\chi_\lambda|_H$ and $\left\langle \chi_\lambda |_H, 1 \right\rangle$ is the number of times the trivial representation of $H$ appears in $\chi_\lambda|_H$. By reciprocity, this is $\left\langle \chi_\lambda, \mathrm{Ind}_H^G(1) \right\rangle$, where $\mathrm{Ind}_H^G$ is the induced representation from $H$ to $G$.

\begin{proposition}
Let $Q$ be a class function on $G$. The induced chain $P(x, y)$ of Proposition \ref{prop: inducedChain} has eigenvalues
\begin{equation}
    \beta_\lambda = \frac{1}{\chi_\lambda(1)} \sum_{s \in G} Q(s) \chi_\rho(s),
\end{equation}
with multiplicity
\begin{equation}
m_\lambda = \left\langle \chi_\lambda|_H, 1 \right\rangle \cdot \left\langle \chi_\lambda|_K, 1 \right\rangle.
\end{equation}
The average square total variation distance to stationarity satisfies
\begin{equation} \label{eqn: averageBound}
    \sum_x \pi(x) \|P_x^\ell - \pi \|_{TV}^2 \le \frac{1}{4} \sum_{\lambda \in \widehat{G}, \lambda \neq 1} m_\lambda \beta_\lambda^{2 \ell}.
\end{equation}
\end{proposition}

\begin{proof}
The eigenvalues of a lumped chain are always some subset of the eigenvalues of the original chain. To determine the multiplicity of the eigenvalue $\beta_\lambda$ in the lumped chain, fix $\lambda: G \to GL_{d_\lambda}(k)$ an irreducible representation of $G$ over the field $k$. Let $M = M^\lambda$ be the $d_\lambda \times d_\lambda$ matrix representation of $\lambda$. That is, each entry $M_{ij}: G \to k$ is a function of $G$. These functions are linearly independent and can be chosen to be orthogonal with respect to 
\[
\left\langle f_1, f_2 \right\rangle  := \frac{1}{|G|} \sum_{g \in G} f_1(g) \overline{f_2}(g)
\]
(see Chapter 3 of \cite{serre}). Let $V_\lambda$ be the space of all linear combinations of the functions $M_{ij}$. If $f \in V_\lambda$, then
\begin{align*}
Pf(x) = \sum_{y \in G} P(x, y)f(y) = \beta_\lambda f(x)
\end{align*}
That is, $V_\lambda$ is the eigenspace for the eigenvalue $\beta_\lambda$ and it has dimension $d_\lambda^2 = \chi_\lambda(1)^2$.

In the lumped chain on $H \backslash G / K$, a basis for the eigenspace for eigenvalue $\beta_\lambda$ are the $H \times K$ invariant functions in $V_\lambda$ \cite{diaconisBoyd}. To determine the dimension of this subspace, note that $G \times G$ can act on $V_\lambda$ by $f^{g_1, g_2}(x) = f(g_1^{-1} x g_2)$. This gives a representation of $G \times G$ on $V_\lambda$. The matrix of this representation is isomorphic to $M \otimes M$, since $M_{ij}(s^{-1} t u) = M_{ij}(s^{-1})M_{ij}(t)M_{ij}(u)$.


This representation restricts to a representation $M_H \otimes M_K$ of $H \times K$, and the dimension of the $H \times K$ invariant functions in $V_\lambda$ is the multiplicity of the trivial representation on $M_H \otimes M_K$. This is
\[
\left\langle \chi_\lambda|_H \otimes \chi_\lambda|_K, 1 \right\rangle = \left\langle \chi_\lambda|_H, 1 \right\rangle \cdot \left\langle \chi_\lambda|_K, 1 \right\rangle. 
\]

To note the total variation inequality, let $1 = \beta_1 \ge \beta_2 \ge  \dots \ge \beta_n \ge \dots \ge\beta_{|S_n|} \ge -1$ be the eigenvalues with eigenfunctions $f_j$ (chosen to be orthonormal with respect to $\pi$), we have
\[
\chi^2_x(\ell) = \sum_{j \neq 1} f_j(x)^2 \beta_j^{2\ell},
\]
where $\| \cdot \|_{2, \pi}$ denotes the $\ell^2$ norm with respect to the distribution $\pi$. Multiplying by $\pi(x)$ and summing over all $x$ in the state space gives
\begin{equation} \label{eqn: generalBound}
\sum_x \pi(x) \left\| \frac{P^\ell_x}{\pi} - 1 \right\|_{2, \pi}^2 = \sum_x \pi(x) \sum_{j \neq 1} f_j(x)^2 \beta_j^{2\ell} = \sum_{\lambda \neq 1} m_\lambda \beta_\lambda^{2 \ell},
\end{equation}
using orthonormality of $f_j$. The total variation bound arises since $4 \|P^\ell_x - \pi \|_{TV}^2 \le \|P^\ell_x/\pi - 1 \|_{2, \pi}^2$.
\end{proof}

\begin{remark}
By a classical bound, for any starting state $x$
\begin{equation*} 
    4 \|P_x^\ell - \pi \|_{TV}^2 \le \frac{1}{\pi(x)} \beta_{*}^{2 \ell},
\end{equation*}
with $\beta_* = \max_{\lambda \neq 1} |\beta_\lambda|$. Since $\pi(x)$ varies widely, this bound can be useful or very off! See the examples in Section \ref{sec: MT_results}.
\end{remark}

In applications, the rate of convergence to stationarity can vary greatly as a function of the starting state. The classical bound can be useful when $\pi(x)^{-1}$ is not too large. The average case bound \eqref{eqn: averageBound}, while sharp, is dominated by starting states with large $\pi(x)$. We do not have (at present) an accurate upper bound for a general starting state. However, for the special case where $H = K = B$ the Borel subgroup of a finite group of Lie Type (e.g.\ Example \ref{ex: gln} of the introduction) a useful bound is derived in \cite{diaconisRam}. This is detailed in Section \ref{sec: MT_results} where Example \ref{ex: gln} is treated. 

\subsection{Random Transpositions and Coagulation-Fragmentation Processes} \label{sec: RT}

Let $G = S_n$ be the symmetric group. The random transpositions Markov chain, studied in \cite{diaconisShah1981} is generated by the measure
\[
Q(\omega) = \begin{cases}
1/n & \omega = id \\
2/n^2 & \omega = (ij), \quad i < j \\
0 & \text{otherwise}
\end{cases}.
\]
This was the first Markov chain where a sharp cut-off for convergence to stationarity was observed. A sharp, explicit rate is obtained in \cite{saloff2008}. They show
\[
\|Q^{*\ell} - u \|_{TV} \le 2 e^{-c}, \,\,\,\,\,\, \ell = \frac{1}{2}n(\log(n) + c).
\]
The asymptotic `profile' (the limit of $\|Q^{*\ell} - u\|_{TV}$ as a function of $c$ for $n$ large) is determined in \cite{teyssier}. Schramm \cite{schramm} found a sharp parallel between random transpositions and the growth of an Erd\H{o}s-R\'{e}nyi random graph: Given vertices $1, 2, \dots, n$, for each transposition $(i, j)$ chosen, add an edge from vertices $i$ to $j$ to generate a random graph. See \cite{berestycki2011} for extensions and a comprehensive review. The results, translated by the coagulation-fragmentation description of the cycles, give a full and useful picture for the simple mean-field model described in the introduction.

It must be emphasized that this mean field model is a very special case of coagulation-fragmentation models studied in the chemistry-physics-probability literature. These models study the dynamics of particles diffusing in an ambient space, and allow general collision kernels (e.g.\ particles close in space may be more likely to join). The books by Bertoin \cite{bertoin} and Pitman \cite{pitman} along with the survey paper of Aldous \cite{aldous1999} are recommended for a view of the richness of this subject. On the other hand, the sharp rates of convergence results available for the mean field model are not available in any generality.

Other lumpings of random transpositions include classical urn models -- the Bernoulli-Laplace model \cite{diaconisShah1987}, \cite{nestoridi}, and random walks on phylogenetic trees \cite{diaconisHolmes}. The sharp analysis of random transpositions transfers, via comparison theory, to give good rates of convergence for quite general random walks on the symmetric group \cite{diaconisSC}, \cite{helfgott}. For an expository survey, see \cite{diaconis2003}.

There is a healthy applied mathematics literature on coagulation-fragmentation. A useful overview which treats discrete problems such as the ones treated here is \cite{ball}. A much more probabilistic development of the celebrated Becker-Doring version of the problem is in \cite{hingant}. This develops rates of convergence using coupling. See also \cite{durrett} for more models with various stationary distributions on partitions.


\subsection{Transvections} \label{sec: transvections}

Fix $n$, a prime $p$, and $q = p^a$ for some positive integer $a$. A \emph{transvection} is an invertible linear transformation of $\FF_q^n$ which fixes a hyperplane but is not the identity. Transvections are convenient generators for the group $SL_n(q)$ because they generalize the basic row operations of linear algebra. These properties are carefully developed in \cite{suzuki} Chapter 1, 9; \cite{artin} Chapter 4.

Using coordinates, let $\ma, \mv \in \FF_q^n$ be two non-zero vectors with $\ma^\top \mv = 0$. A transvection, denoted $T_{\ma, \mv} \in GL_n(q)$ is the linear map such that
\[
T_{\ma, \mv}(\mx) = \mx + \mv(\ma^\top \mx), \quad \mx \in \FF_q^n.
\]
It adds a multiple of $\mv$ to $\mx$, the amount depending on the `angle' between $\ma$ and $\mx$. As a matrix, $T_{\ma, \mv} = [a_1 \mv \, a_2 \mv \hdots a_n \mv] = I + \mv \ma^\top$. Multiplying $\ma$ by a non-zero constant and dividing $\mv$ by the same constant doesn't change $T_{\ma, \mv}$. Let us agree to normalize $\mv$ by making its last non-zero coordinate equal to $1$. Let $\mathcal{T}_{n, q} \subset SL_n(q)$ be the set of all transvections. 

An elementary count shows
\begin{equation}
    |\mathcal{T}_{n, q}| = \frac{(q^n - 1)(q^{n-1} - 1)}{q - 1}.
\end{equation}

It is easy to generate $T \in \mathcal{T}_{n, q}$ uniformly: Pick $\mv \in \FF_q^n$ uniformly, discarding the zero vector. Normalize $\mv$ so the last non-zero coordinate, say index $j$, is equal to $1$. Pick $a_1, a_2, \dots, a_{j-1}, a_{j+1}, \dots a_n$ uniformly in $\FF_q^{n-1} - \{0 \}$ and set $a_j$ so that $\ma^\top \mv = 0$. Clearly $T_{\ma, \mv}$ fixes the hyperplane $\{\mx: \ma^\top \mx = 0 \}$.

\begin{example}
Taking $\mv = \me_1, \ma = \me_2$ gives the transvection with matrix
\[
T_{\me_1, \me_2} = \begin{pmatrix} 1 & 1 & 0 & \hdots & 0 \cr 
0 & 1 & 0 & \hdots & 0 \cr
0 & 0 & 1 & \hdots & 0 \cr
\vdots & \vdots & \ddots & \hdots & 0 \cr
0 & 0 & 0 & \hdots & 1 \end{pmatrix}.
\]
This acts on $\mx$ by adding the second coordinate to the first. Similarly, the basic row operation of adding $\theta$ times the $i$th coordinate to the $j$th is given by $T_{\me_j, \theta \me_i}$.
\end{example}

\begin{lemma}
The set of transvections $\mathcal{T}_{n, q}$ is a conjugacy class in $GL_n(q)$.
\end{lemma}
\begin{proof}
Let $M \in GL_n(q)$, so $M T_{\me_2, \me_1} M^{-1}$ is conjugate to $T_{\me_2, \me_1}$. Then,
\[
M T_{\me_2, \me_1} M^{-1}(\mx) = M ( M^{-1} \mx + (M^{-1} \mx)_2 \me_1 ) = \mx + (M^{-1} \mx)_2 M \me_1
\]
Let $\ma$ be the second column of $(M^{-1})T_{\me_2, \me_1}$ and $\mv$ the first column of $M$ and check this last is $T_{\ma, \mv}$ (and $\ma^\top \mv = 0$). Thus, transvections form a conjugacy class.
\end{proof}

\subsection{Gaussian Elimination and the Bruhat Decomposition} \label{sec: gaussianElim}

The reduction of a matrix $M\in GL_n(q)$ to standard form by row operations is a 
classical topic in introductory linear algebra courses.  It gives efficient, numerically
stable ways to solve linear equations, compute inverses, and calculate determinants.  There are many
variations.  

\begin{example}
Consider the sequence of row operations
\[
M = \begin{pmatrix} 0 & 3 & 2 \cr 1 & 2 & 0 \cr 3 & 0 & 5 \end{pmatrix} \to \begin{pmatrix} 0 & 3 & 2 \cr 1 & 2 & 0 \cr 0 & -6 & 5 \end{pmatrix} \to \begin{pmatrix} 0 & 3 & 2 \cr 1 & 2 & 0 \cr 0 & 0 & 9 \end{pmatrix} \to \begin{pmatrix} 1 &2 & 0 \cr 0 & 3 & 2 \cr 0 & 0 & 9 \end{pmatrix}  = U.
\]
The first step subtracts $3$ times row $2$ from row $3$, multiplication by 
\[
L_1 = \begin{pmatrix} 1 & 0 & 0 \cr 0 & 1 & 0 \cr 0 & -3 & 1 \end{pmatrix}.
\]
The second step adds $2$ times row $1$ to row $3$, multiplication by 
\[
L_2 = \begin{pmatrix} 1 & 0 & 0 \cr 0 & 1 & 0 \cr 2 & 0 & 1 \end{pmatrix}.
\]
The third (pivot) step brings the matrix to upper triangular form by switching rows $1$ and $2$, multiplication by 
\[
\omega_1 = \begin{pmatrix} 0 & 1 & 0 \cr 1 & 0 & 0 \cr 0 & 0 & 1 \end{pmatrix}.
\]
This gives $\omega_1 L_2 L_1 M = U \implies M = L_2^{-1} L_1^{-1} \omega_1^{-1} U =  L \omega U$ with $L = L_2^{-1} L_1^{-1}, \omega = \omega_1^{-1}$.
\end{example}

If $\cL, \cB$ are the subgroups of lower and upper triangular matrices in $GL_n(q)$, this gives
\begin{equation} \label{eqn: LU}
    GL_n(q) = \bigsqcup_{\omega \in S_n} \cL \omega \cB
\end{equation}

Any linear algebra book treats these topics.  A particularly clear version which uses Gaussian elimination as a gateway to Lie theory is in Howe \cite{howe}. Articles by Lusztig \cite{lusztig} and Strang \cite{strang} have further historical, mathematical, and practical discussion.

Observe that carrying out the final pivoting step costs $d_c(\omega,\id)$ operations where
$d_c(\omega,\id)$, the Cayley distance of $\omega$ to the identity, is the minimum number of transpositions
required to sort $\omega$ (with arbitrary transpositions $(i,j)$ allowed).  Cayley proved
$d_c(\omega,\id) = n - \#\hbox{cycles in $\omega$}$ (see \cite{diaconis2003}). In the example above
$n=3$, $\omega = 213$ has two cycles and 
$3-2 =1$ -- one transposition sorts $\omega$.

How many pivot steps are needed `on average'? This becomes the question of the number of cycles in a pick from Mallows measure $\pi_q$. Surprisingly, this is a difficult question. Following partial answers by Gladkich and Peled \cite{gladkichPeledCycle}, this problem is solved in \cite{jimmyPreprint}. This shows that when $q > 1$ the limiting behavior of the number of even cycles under $\pi_q$ depends has an approximate normal distribution with mean and variance proportional to $n$. However, the number of odd cycles has bounded mean and variance.

\paragraph{The Bruhat Decomposition}

In algebraic group theory one uses
\begin{equation} \label{eqn: bruhatDecomp}
    GL_n(q) = \bigsqcup_{\omega \in S_n} \cB \omega \cB.
\end{equation}

This holds for any semi-simple group over any field with $\cB$ replaced by the Borel group (the largest solvable subgroup) and $S_n$ replaced by the Weyl group.

Let $\omega_0 = \left(\begin{smallmatrix} 
1 & 2 & \cdots & n\\
n & n-1 & \cdots & 1
\end{smallmatrix}\right)$ be the reversal permutation in $S_n$. Since $\cL = \omega_0 \cB \omega_0$, \eqref{eqn: bruhatDecomp} is equivalent to the LU decomposition \eqref{eqn: LU}. Given $M \in GL_n(q)$, Gaussian elimination on $\omega_0 M$ can be used to find $\omega_0 M \in L \omega' U$ and thus $M \in B \omega B$ with $\omega = \omega_0 \omega'$.

The subgroup $\cB$ gives rise to the quotient $GL_n(q)/\cB$. This may be pictured as the space of `flags'. Here a flag $F$ consists of an increasing sequence of subspaces $F = F_1 \subset F_2 \subset \hdots \subset F_n$ with $\text{dim}(F_i) = i$. Indeed, $GL_n(q)$ operates transitively on flags and the subgroup fixing the \emph{standard flag} $E = \left\langle \me_1 \right\rangle \subset \left\langle \me_1, \me_2 \right\rangle \hdots \subset \left\langle \me_1, \dots, \me_n \right\rangle$ is exactly $\cB$. This perspective will be further explained and used in Section \ref{sec: topCard} to study a function of the double coset Markov chain on $\cB \backslash GL_n(q) / \cB$.

\begin{remark}
The double-cosets of $GL_n(q)$ define equivalence classes for any subgroup of $GL_n(q)$. For the matrices $SL_n(q)$ with determinant $1$, these double-cosets again induce the Mallows distribution on permutations. More precisely, for $x \in SL_n(q)$, let $[x]_{SL_n(q)} = \{x' \in SL_n(q) : x' \in \cB x \cB \}$ be the equivalence class created by the double-coset relation $\cB \backslash GL_n(q) / \cB$, within $SL_n(q)$. Note that two matrices $x, x' \in SL_n(q)$ could be in the same double coset with $x' = b_1 x b_2$, but $b_1, b_2 \notin SL_n(q)$ (necessarily, $\mathrm{det}(b_1) = \mathrm{det}(b_2)^{-1}$). 

Then, $|[\omega]_{SL_n(q)}|/|SL_n(q)| = p_q(\omega)$. This follows since $|GL_n(q)| = (q-1) \cdot |SL_n(q)|$, and $|\cB \omega \cB| = (q - 1) \cdot |[\omega]_{SL_n(q)}|$. If $M \in GL_n(q)$ and $M \in B \omega B$, then $M/\det(M) \in [\omega]_{SL_n(q)}$. Conversely, for each $M \in SL_n(q)$ there are $(q - 1)$ unique matrices in $GL_n(q)$ created by multiplying $M$ by $1, 2, \dots, q - 1$.
\end{remark}

\section{Double Coset Walks on $\cB \backslash GL_n(q) / \cB$} \label{sec: MT_results}

Throughout this section, $\cB$ is the group of upper triangular matrices in $GL_n(q)$, $\mathcal{T}_{n,q}$ is the conjugacy class of transvections in $GL_n(q)$. This gives the probability measure on $GL_n(q)$ defined by
\[
Q(M) = \begin{cases}
\frac{1}{|\mathcal{T}_{n,q}|} & \text{if} \,\,\,\, M \in \mathcal{T}_{n,q} \\
0 & \text{else}
\end{cases}.
\]

Note the random transvections measure $Q$ is supported on $SL_n(q)$, a subgroup of $GL_n(q)$. This means that the random walk on $GL_n(q)$ driven by $Q$ is not ergodic (there is zero probability of moving $x$ to $y$ if $x, y$ are matrices with different determinants). However, $Q$ is a class function on $GL_n(q)$ since transvections form a conjugacy class. The image of the uniform distribution on $SL_n(q)$ mapped to $\cB \backslash GL_n(q) /\cB$ is the Mallows measure $\pi_q(\omega) = q^{I(\omega)}/[n]_q!$.

Section \ref{sec: heckAlgebraIntro} introduces the definition of the Markov chain as multiplication in the Hecke algebra, which is further explained in Section \ref{sec: heckeComputations}. Section \ref{sec: evs} gives combinatorial expressions for the eigenvalues and their multiplicities, needed to apply Theorem \ref{thm: general} for this case. Section \ref{sec: identity} shows that for the induced Markov chain on $S_n$, starting from $id \in S_n$, order $n$ steps are necessary and sufficient for convergence. Section \ref{sec: reversal} studies the chain starting from the reversal permutation, for which only order $\log(n)/2 \log(q)$ steps are required. Finally, Section \ref{sec: typical} considers starting from a `typical' element, according to the stationary distribution, for which $\log(n)/\log(q)$ steps are necessary and sufficient.

These results can be compared to Hildebrand's Theorem 1.1 \cite{hildebrand} which shows that the walk driven by $Q$ on $GL_n(q)$ converges in $n + c$ steps (uniformly in $q$). Our results thus contribute to the program of understanding how functions of a Markov chain behave and how the mixing time depends on the starting state. In this case, changing the starting state gives an exponential speed up.

\subsection{Hecke Algebras and the Metropolis Algorithm} \label{sec: heckAlgebraIntro}

The set of $\cB-\cB$ double cosets of $GL_n(q)$ has remarkable structure. For $\omega \in S_n$, let $T_\omega = \cB \omega \cB$.  Linear combinations of double cosets form an
algebra (over $\CC$, for example). 

\begin{definition} \label{def: IwahoriHecke}
The \emph{Iwahori-Hecke algebra} $H_n(q)$ is spanned by the symbols $\{T_\omega \}_{\omega \in S_n}$ and generated by $T_i = T_{s_i}$ for $s_i = (i, i+1), 1 \le i \le n -1$, with the relations
\begin{equation}
\begin{cases}
T_{s_i}T_\omega = T_{s_i\omega} & \hbox{if $I(s_i\omega) = I(\omega)+1$,} \\
T_{s_i}T_\omega = qT_{s_i \omega} +(q-1)T_{\omega} & \hbox{if $I(s_i\omega) = I(\omega)-1$} 
\end{cases},
\label{Hnqdefn}
\end{equation}
where $I(\omega)$ is the usual length function on $S_n$ ($I(s_i \omega) = I(\omega)\pm 1$).
\end{definition}

Consider the flag space $\cF = G/\cB$. The group $GL_n(q)$ acts on the left of $\cF$. One can see $H_n(q)$ acting on the right of $\cF$ and in fact
\[
H_n(q) = \End_{GL_n(q)}(GL_n(q)/\cB).
\]
The Hecke algebra is the full commuting algebra of $GL_n(q)$ acting on $GL_n(q)/\cB$. Because transvections form a conjugacy class, the sum of transvections is in the center of the group algebra $\C[GL_n(q)]$, and so it may be regarded as an element of $H_n(q)$. This will be explicitly delineated and the character theory of $H_n(q)$ used to do computations. 

\subsection{The Metropolis Connection}  \label{sec: metropolis}

The relations \eqref{Hnqdefn} can be interpreted probabilistically. Consider what \eqref{Hnqdefn} says as linear algebra: Left multiplication by $T_{s_i}$ can take $\omega$ to $\omega$ or $s_i \omega$. The matrix of this map (in the basis $\{T_\omega\}_{\omega \in S_n}$) has $\omega, \omega'$ entry
\begin{align*}
    \begin{cases} 1 & \text{if} \quad I(s_i \omega) = I(\omega) + 1, \omega' = s_i \omega \\
    q & \text{if} \quad I(s_i \omega) = I(\omega) - 1, \quad \omega' = \omega \\
    q - 1 & \text{if} \quad I(s_i \omega) = I(\omega) -1, \quad \omega' = s_i \omega \\
    0 & \quad \text{otherwise}
    \end{cases}.
\end{align*}

For example, on $GL_3(q)$ using the ordered basis $T_{id}, T_{s_1}, T_{s_2}, T_{s_1s_2}, T_{s_2s_1}, T_{s_1s_2s_1}$, the matrix of left multiplication by $s_1$ is
\[
T_{s_1} = \begin{pmatrix}
0 & q & 0 & 0 & 0 & 0 \cr
1 & q-1 & 0 & 0 & 0 & 0 \cr
0 & 0 & 0 & q & 0 & 0 \cr
0 & 0 & 1 & q-1 & 0 & 0 \cr
0 & 0 & 0 & 0 & 0 & q \cr
0 & 0 & 0 & 0 & 1 & q-1 \cr
\end{pmatrix}.
\]
The first column has a $1$ in row $s_1$ because $I(s_1) > I(id)$. The second column has entries $q$ and $q-1$ in the first two rows because $I(s_1^2) = I(id) < I(s_1)$.

We can also write the matrices for multiplication defined by $T_{s_2}$ and $T_{s_1s_2s_1}$ as
\begin{align*}
&T_{s_2} = \begin{pmatrix}
0 & 0 & q & 0 & 0 & 0 \cr
0 & 0 & 0 & 0 & q & 0 \cr
1 & 0 & q-1 & 0 & 0 & 0 \cr
0 & 0 & 0 & 0 & 0 & q\cr
0 & 1 & 0 & 0 & q-1 & 0 \cr
0 & 0 & 0 & 1 & 0 & q-1 \cr
\end{pmatrix} \\
& T_{s_1s_2s_1} = \begin{pmatrix}
0 & 0 & 0 & 0 & 0 & q^3 \cr
0 & 0 & 0 & q^2 & 0 & q^2(q - 1) \cr
0 & 0 & 0 & 0 & q^2 & q^2(q - 1) \cr
0 & q & 0 & q(q - 1) & q(q-1) & q(q-1)^2 \cr
0 & 0 & q & q(q-1) & q(q-1) & q(q-1)^2 \cr
1 & q-1 & q-1 & (q-1)^2 & (q-1)^2 & (q - 1)^3 + q(q-1)\cr
\end{pmatrix}.
\end{align*}

Observe that all three matrices above have constant row sums ($q, q^2$, and $q^3$ respectively). Dividing by these row sums gives three Markov transition matrices: $M_1, M_2,$ and $M_{121} = M_1M_2M_1$.

These matrices have a simple probabilistic interpretation: Consider, for $\overline{q} = 1/q$, the matrix defined
\[
M_1 = \begin{pmatrix}
0 & 1 & 0 & 0 & 0 & 0 \cr
\overline{q} & 1 - \overline{q} & 0 & 0 & 0 & 0 \cr
0 & 0 & 0 & 1 & 0 & 0 \cr
0 & 0 & \overline{q} & 1 - \overline{q} & 0 & 0 \cr
0 & 0 & 0 & 0 & 0 & 1 \cr
0 & 0 & 0 & 0 & \overline{q} & 1 - \overline{q} \cr
\end{pmatrix}.
\]
The description of this Markov matrix is: From $\omega$, propose $s_1 \omega$:
\begin{itemize}
    \item If $I(s_1 \omega) > I(\omega)$, go to $s_1 \omega$.
    \item If $I(s_1 \omega) < I(\omega)$, go to $s_1 \omega$ with probability $1/q$, else stay at $\omega$.
\end{itemize}
This is exactly the Metropolis algorithm on $S_n$ for sampling from $\pi_q(\omega)$ with the proposal given by the deterministic chain `multiply by $s_1$'. The matrices $M_i, 1 \le i \le n -1$ have a similar interpretation, and satisfy
\[
\pi_q(\omega) M_i(\omega, \omega') = \pi_q(\omega') M_i(\omega', \omega).
\] 
The Metropolis algorithm always results in a reversible Markov chain. See \cite{diaconisSCMetropolis} for background. It follows that any product of $\{M_i \}$ and any convex combination of such products yields a $\pi_q$ reversible chain. Note also that the Markov chain on $S_n$ is automatically reversible since it is induced by a reversible chain on $GL_n$.

\begin{corollary}
The random transvections chain on $GL_n(q)$ lumped to $\cB - \cB$ cosets gives a $\pi_q$ reversible
Markov chain on $S_n$.
\end{corollary}

\begin{proof}
Up to normalization, the matrix $D$ in Theorem \ref{thm: D} is a positive linear combination of
Markov chains corresponding to multiplication by 
\[
T_{(i,j)} = T_{s_i}T_{s_{i+1}} \cdots T_{s_{j-2}}T_{s_{j-1}}T_{s_{j-2}}\cdots T_{s_{i+1}}T_{s_i}.
\]
This yields a combination of the reversible chains $M_{i}M_{i+1} \hdots M_{j-1} \hdots M_{i}$.
\end{proof}

\begin{example}
The transition matrix of the transvections chain on $GL_3(q)$ lumped to $S_3$ is $\frac{1}{|\mathcal{T}_{n, q}|}D$, with $D$ defined
\begin{align*}
&D = (2q^2 - q - 1) \mathbf{I} +  
\\
&\,\,\,\,\,\, (q-1)\begin{pmatrix}
0 &q^2 &q^2 &0 &0 &q^3 \\
q &q(q-1) &0 &q^2 &q^2 &q^2(q-1) \\
q &0 &q(q-1) &q^2 &q^2 &q^2(q-1) \\
0 &q &q &2q(q-1) &q(q-1) &q^2+q(q-1)^2 \\
0 &q &q &q(q-1) &2q(q-1) &q^2+q(q-1)^2 \\
1 &q-1 &q-1 &q+(q-1)^2 &q+(q-1)^2 &(q-1)^3+3q(q-1)
\end{pmatrix}.
\end{align*}
When $q=2$, the lumped chain has transition matrix
$$\frac{1}{21} \begin{pmatrix}
5 &4 &4 &0 &0 &8 \\
2 &7 &0 &4 &4 &4 \\
2 &0 &7 &4 &4 &4 \\
0 &2 &2 &9 &2 &6 \\
0 &2 &2 &2 &9 &6 \\
1 &1 &1 &3 &3 &12
\end{pmatrix}.
$$
We report that this example has been verified by several different routes including simply
running the transvections chain, computing the double coset representative at each step and
estimating the transition rates from a long run of the chain.
\end{example}

\begin{remark}
The random transvections Markov chain on $S_n$ is the `q-deformation' of random transpositions on $S_n$. That is, as $q$ tends to $1$, the transition matrix tends to the transition matrix of random transpositions. To see this, recall $|\mathcal{T}_{n, q}| = (q^n - 1)(q^{n-1} - 1)/(q-1)$ and use L'Hopitals rule to note, for any integer $k$,
\begin{align*}
    &\lim_{q \to 1^+} \frac{(q - 1)^2 q^k(q-1)^\ell}{(q^n - 1)(q^{n-1} - 1)} =  \begin{cases} 0 & \text{if} \,\, \ell > 0 \\
    \frac{1}{n(n-1)} & \text{if} \,\,\, \ell = 0
    \end{cases}.
\end{align*}

\end{remark}

The interpretation of multiplication on the Hecke algebra as various `systematic scan' 
Markov chains is developed in \cite{diaconisRam}, \cite{bufetov}.  It works for other types in several variations.
We are surprised to see it come up naturally in the present work.

 The following corollary provides the connection to \cite[(3.16),(3.18),(3.20)]{Ra97}
 and  \cite[Proposition 4.9]{diaconisRam}.
 
 \begin{corollary}
Let $J_1 = 1$ and let
$J_k = T_{s_{k-1}}\cdots T_{s_2}T_{s_1}T_{s_1}T_{s_2}\cdots T_{s_{k-1}}$,
for $k\in \{2, \ldots, n\}$.
Then
$$D = \sum_{k=1}^n q^{n-k}(J_k-1).$$
 \end{corollary}
 
\begin{proof}
Using that $J_k = T_{s_{k-1}}J_{k-1}T_{s_{k-1}}$, check, by induction, that
$$J_k = q^{k-1} + (q-1)\sum_{j=1}^{k-1} q^{j-1}T_{(j,k)}.$$
Thus
\begin{align*}
\sum_{k=1}^n q^{n-k}(J_k-1) 
&= 0+ \sum_{k=2}^n \big(q^{n-k+(k-1)} -q^{n-k} + (q-1)\sum_{j=1}^{k-1} q^{n-k+j-1}T_{(j,k)}\big) \\
&=(n-1)q^{n-1}  - [n-1]_q + (q-1)D_{(21^{n-2})} = D,
\end{align*}
where $D_{(21^{n-2})} := \sum_{i < j} q^{(n - 1) - (j-i)} T_{(i, j)}$.
\end{proof}

\subsection{Eigenvalues and Multiplicities} \label{sec: evs}

Hildebrand \cite{hildebrand} determined the eigenvalues of the random walk driven by $Q$ on $GL_n(q)$. His arguments use MacDonald's version of J.A.\ Green's formulas for the characters of $GL_n(q)$ along with sophisticated use of properties of Hall-Littlewood polynomials. Using the realization of the walk on the Hecke algebra, developed below in Section \ref{sec: heckeComputations}, and previous work of Ram and Halverson \cite{halvseronRam}, we can find cleaner formulas and proofs. Throughout, we have tried to keep track of how things depend on both $q$ and $n$ (the formulas are more simple when $q = 2$).

\begin{theorem} \label{thm: evs_transvections}
\begin{enumerate}[(a)]
    \item The eigenvalues $\beta_\lambda$ of the Markov chain $P(x, y)$ driven by the random transvections measure $Q$ on $\cB \backslash GL_n(q) / \cB$ are indexed by partitions $\lambda \vdash n$. These are
    \begin{equation} \label{eqn: betalam}
        \beta_\lambda = \frac{1}{|\mathcal{T}_{n, q}|} \left( q^{n-1} \sum_{b \in \lambda} q^{ct(b)} - \frac{q^{n} - 1}{q - 1} \right),
    \end{equation}
    with $b$ ranging over the boxes of the partition $\lambda$. For the box in row $i$ and column $j$, the content is defined as $ct(b) = j - i$, as in \cite{macdonald}.
    
    \item The multiplicity of $\beta_\lambda$ for the induced Markov chain on $\cB \backslash GL_n(q) / \cB$ is $f_\lambda^2$, where
    \begin{equation}
        f_\lambda = \frac{n!}{\prod_{b \in \lambda} h(b)}.
    \end{equation}
    Here $h(b)$ is the hook length of box $b$ \cite{macdonald}. 
    
    \item The multiplicity of $\beta_\lambda$ for the Markov chain induced by $Q$ on $GL_n(q)/\cB$ is
    \begin{equation}
        d_{\lambda} = f_\lambda \cdot \frac{[n]_q!}{\prod_{b \in \lambda} [h(b)]_q}.
    \end{equation}
\end{enumerate}
\end{theorem}

The argument uses the representation of the Markov chain as multiplication on the Hecke algebra. This is developed further in Section \ref{sec: heckeComputations}. 
\begin{proof}
\textbf{(a):} Let $D$ be the sum of transvections as in Theorem \ref{thm: D}. By \cite{halvseronRam} (3.20), the irreducible representation $H_n^\lambda$ of $H_n$ indexed by $\lambda$ has a `seminormal basis' $\{ v_T \mid \widehat{S}_n^{\lambda} \}$ such that $J_k v_T = q^{k-1} q^{c(T(k))} v_T$, where $T(k)$ is the box containing $k$ in $T$. Thus,
\[
Dv_T = \sum_{k = 1}^n q^{n - k} (q^{k - 1}q^{c(T(k))} - 1)v_T = \left( - \left( \frac{q^n - 1}{q - 1} \right) + \sum_{k = 1}^n q^{n-1} q^{c(T(k))} \right) v_T.
\]

\textbf{(b):} The dimension of the irreducible representation of $H_n$ indexed by $\lambda$ is the same as the dimension of the irreducible representation of $S_n$ indexed by $\lambda$, which is well-known as the hook-length formula:
\[
\text{dim}(H_n^\lambda) = \text{Card}(S_n^\lambda) = \frac{n!}{\prod_{b \in \lambda} h(b)}.
\]

\textbf{(c):} With $G = GL_n(q), H = H_n$, the result follows since
\[
\text{dim}(G^\lambda) = \frac{[n]_q!}{\prod_{b \in \lambda} [h(b)]_q} \quad \text{and} \quad \mathbf{1}_\cB^G = \sum_{\lambda \in \widehat{S_n}} G^\lambda \otimes H^\lambda.
\]

\end{proof}

\begin{example}
Equation \ref{eqn: betalam} for some specific partitions gives
\begin{align*}
    &\beta_{(n)} = \frac{1}{|\mathcal{T}_{n, q}|} \left( q^{n-1} \cdot \left( \frac{q^n - 1}{q - 1} \right) - \frac{q^n - 1}{q - 1} \right) = 1 \\
    &\beta_{(n - 1, 1)} = \frac{1}{|\mathcal{T}_{n, q}|} \left( q^{n-1}  \left( q^{-1} + \frac{q^{n-1} - 1}{q - 1} \right) - \frac{q^n - 1}{q - 1} \right) = \frac{q^{n - 2} - 1}{q^{n-1} - 1} \\
    &\beta_{(n - 2, 1^2)} = \frac{1}{|\mathcal{T}_{n, q}|} \left( q^{n-1} \left( q^{-1} +q^{-2} + \frac{q^{n-2} - 1}{q - 1} \right)  - \frac{q^n - 1}{q - 1} \right) = \frac{q^{n - 3} - 1}{q^{n-1} - 1} \\
    &\beta_{(n - 2, 2)} = \frac{1}{|\mathcal{T}_{n, q}|} \left( q^{n-1} \left(1 + q^{-1}  + \frac{q^{n-2} - 1}{q - 1} \right) - \frac{q^n - 1}{q - 1} \right) = \frac{q^{n - 2} - 1}{q^{n} - 1} \\
    &\beta_{(2, 1^{n-2})} = \frac{1}{|\mathcal{T}_{n, q}|} \left( q^{n-1} \left( 1 + q + q^{-1} + q^{-2} +\hdots + q^{-(n-2)}\right)  - \frac{q^n - 1}{q - 1} \right) = \frac{q - 1}{q^{n-1} - 1} \\
    &\beta_{(1^{n})} = \frac{1}{|\mathcal{T}_{n, q}|} \left( q^{n-1}(q^0 + q^{-1} + \hdots + q^{-(n-1)}) - \frac{q^n - 1}{q - 1} \right) = 0.
\end{align*}
\end{example}

The following simple lemma will be used several times in the following sections. It uses the usual majorization partial order (moving up boxes) on partitions of $n$, \cite{macdonald}. For example, when $n = 4$ the ordering is $1111 \prec 211 \prec 22 \prec 31 \prec 4$.

\begin{lemma} \label{lem: majorization}
The eigenvalues $\beta_\lambda$ of Theorem \ref{thm: evs_transvections} are monotone increasing in the majorization order.
\end{lemma}

\begin{proof}
Moving a single box (at the corner of the diagram of $\lambda$) in position $(i, j)$ to position $(i', j')$ necessitates $i' < i$ and $j' > j$, and so $q^{j-i} < q^{j' - i'}$. Since any $\lambda \prec \lambda'$ can be obtained by successively moving up boxes, the proof is complete. For example, 

\begin{center}
\begin{tikzpicture}
\draw (-.5, -.5) node {$\lambda = $ };
\foreach \x in {0,...,3}
	\filldraw (\x*.25, 0) circle (.5mm);
\foreach \x in {0,...,3}
	\filldraw (\x*.25, -.25) circle (.5mm);
\foreach \x in {0,...,2}
	\filldraw (\x*.25, -.5) circle (.5mm);
\foreach \x in {0,...,0}
	\filldraw (\x*.25, -.75) circle (.5mm);
	\filldraw [red] (1*.25, -.75) circle (.5mm);

\draw (2, -.5) node {$\lambda^\prime = $ };
\foreach \x in {0,...,3}
	\filldraw (\x*.25 + 2.75, 0) circle (.5mm);
	\filldraw [red] (4*.25 + 2.75, 0) circle (.5mm);
\foreach \x in {0,...,3}
	\filldraw (\x*.25 + 2.75, -.25) circle (.5mm);
\foreach \x in {0,...,2}
	\filldraw (\x*.25 + 2.75, -.5) circle (.5mm);
\foreach \x in {0,...,0}
	\filldraw (\x*.25 + 2.75, -.75) circle (.5mm);

\end{tikzpicture}
\end{center}
\end{proof}

The partition $1^n = (1, 1, \dots, 1)$ is the minimal element in the partial order, and since $\beta_{(1^n)} = 0$, we have the following:

\begin{corollary}
For any $\lambda$,
\[
\beta_\lambda \ge 0.
\]
\end{corollary}

\begin{corollary} \label{cor: betaBound}
If $\lambda = \lambda_1 \ge \lambda_2 \ge \hdots \ge \lambda_k > 0$ is a partition of $n$ with $\lambda_1 = n - j, j \le n/2$, then 
\begin{align*}
    \beta_\lambda \le \beta_{(n-j, j)} &= \frac{q^j(q^{n - j - 1} - 1)(q^{n - j} - 1) + (q^{n - 2} - 1)(q^j - 1)}{(q^{n-1} - 1)(q^n -1)} \\
    &\le q^{-j}(1 + q^{-(n - 2j + 1)})(1 + q^{-(n - 1)})(1 + q^{-(n - 2)}).
\end{align*}
\end{corollary}

\begin{proof}
The first inequality follows from Lemma \ref{lem: majorization}. The formula for $\beta_{(n - j, j)}$ is a simple calculation from Equation \eqref{eqn: betalam}. Recall the elementary inequalities:
\[
\frac{1}{q^r} < \frac{1}{q^r - 1} < \frac{1}{q^{r-1}}, \quad \frac{1}{1 - q^{-r}} = 1 + \frac{1}{q^r - 1} < 1 + q^{-(r -1)}.
\]
These give the inequality for $\beta_{(n - j, j)}$:
\begin{align*}
    \beta_{(n - j, j)} &= \frac{q^j(q^{n - j - 1} - 1)(q^{n - j} - 1) + (q^{n - 2} - 1)(q^j - 1)}{(q^{n-1} - 1)(q^n -1)} \\
    &\le \left( \frac{q^j q^{n - j - 1}q^{n - j} + q^{n-2}q^j}{q^{n-1} q^n} \right) \frac{1}{(1 - q^{-(n-1)})(1 - q^{-n})} \\
    &\le (q^{-j} + q^{-(n - j + 1))})(1 + q^{-(n-2)})(1 + q^{-(n-1)}).
\end{align*}
\end{proof}

In the following sections we will use a further bound from Corollary \ref{cor: betaBound}.

\begin{corollary} \label{cor: alpha}
Define $\kappa_{n, q, j} := (1 + q^{-(n - 2j + 1)})(1 + q^{-(n - 1)})(1 + q^{-(n - 2)})$ and
\[
\alpha_{n, q} := \max_{1 \le j \le n/2} \frac{\log(\kappa_{n, q, j})}{j}.
\]
Then for all $n, q$,
\[
\alpha_{n, q} \le \frac{6}{n}.
\]
\end{corollary}

\begin{proof}
Using $1 + x \le e^x$, there is the initial bound
\begin{align*}
    \kappa_{n, q, j} &\le \exp \left( q^{-(n - 2j + 1)} + q^{-(n-1)} + q^{-(n-2)} \right) \\
    &\le \exp \left( 3 q^{-(n - 2j)} \right),
\end{align*}
which uses that $j \ge 1$ and so $(n - 2) \ge (n -2j)$. Then,
\begin{align*}
 \frac{\log(\kappa_{n, q, j})}{j} \le \frac{3 q^{-(n - 2j)}}{j}.
\end{align*}
With $f(x) = 3q^{-(n - 2x)}/x$, we have
\[
f'(x) = \frac{3q^{-(n - 2x)}(2x \log(q) - 1)}{x^2}.
\]
Since $2 \log(2) > 1$, we see that $f(x)$ is increasing for $x \ge  1$ and any $q \ge 2$. Thus, for $1 \le j \le n/2$, $f(j)$ is maximized for $j = n/2$, which gives
\[
\alpha_{n, q} \le \max_{1 \le j \le n/2} \frac{3 q^{-(n - 2j)}}{j} \le \frac{6}{n}.
\]
\end{proof}

\section{Mixing Time Analysis}

In this section, the eigenvalues from Section \ref{sec: evs} are used to give bounds on the distance to stationarity for the random transvections Markov chain on $S_n$. Section \ref{sec: repTheoryBounds} reviews the tools which are needed for the bounds from specific starting states. Section \ref{sec: identity} proves results for the chain started from the identity element, Section \ref{sec: reversal} proves results for the chain started from the reversal permutation, and Section \ref{sec: typical} contains bounds for the average over all starting states.

\subsection{Eigenvalue Bounds} \label{sec: repTheoryBounds}

The following result from \cite{diaconisRam} will be the main tool for achieving bounds on the chi-square distance of the chain from different starting states. 

\begin{proposition}[Proposition 4.8 in \cite{diaconisRam}] \label{prop: DRresult}
Let $H$ be the Iwahori-Hecke algebra corresponding to a finite real reflection group $W$. Let $K$ be a reversible Markov chain on $W$ with stationary distribution $\pi$ determined by left multiplication by an element of $H$ (also denoted by $K$).  The following identities are true:
\begin{enumerate}[(a)]
    \item $\chi^2_x(\ell) = q^{-2 I(x)}\sum_{\lambda \neq \mathbf{1}} t_\lambda \chi_H^\lambda(T_{x^{-1}} K^{2 \ell} T_x)$, $x \in W$,
    \item $\sum_{x \in W} \pi(x) \chi^2_x(\ell) = \sum_{\lambda \neq \mathbf{1}} d_\lambda \chi_H^\lambda(K^{2 \ell})$,
\end{enumerate}
where $\chi_H^\lambda$ are the irreducible characters, $t_\lambda$ the generic degrees, and $d_\lambda$ the dimensions of the irreducible representations of $H$.
\end{proposition}

\begin{figure}[ht]
    \centering
    \includegraphics[scale=0.5]{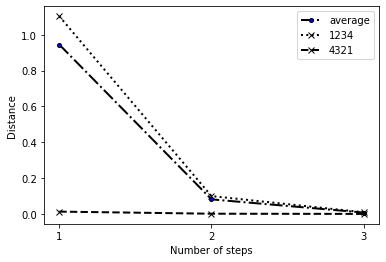}
    \caption{The quantities $\chi^2_x(\ell)$ for $x = id, \omega_0$, and $\sum_x \pi_q(x) \chi_x^2(\ell)$ when $n = 4, q = 3$, for timesteps $\ell = 1, 2, 3$.}
    
    \label{fig: compareTimes}
\end{figure}

In general the right-hand-side of (a) could be difficult to calculate, but it simplifies for the special cases $x = id, x = \omega_0$. These calculations, and the analysis of the sum, are contained in the following sections. Figure \ref{fig: compareTimes} shows the exact values for the chi-squared distance after several steps, from different starting states.

The right hand side of the equations in Proposition \ref{prop: DRresult} will involve the following quantities, defined for $\lambda \vdash n$:
\begin{itemize}
    \item $n_\lambda = \sum_{i = 1}^{|\lambda|} (i - 1) \lambda_i$,
    
    \item $c_\lambda = \sum_{b \in \lambda} c(b)$, where $c(b) = j - i$ if box $b$ is in column $j$ and row $i$
    
    \item $t_\lambda = q^{n_\lambda} \cdot r_\lambda$, where $r_\lambda = \frac{[n]_q!}{\prod_{b \in \lambda} [h_b]_q}$, $[k]_q = (q^k - 1)/(q - 1)$ and $[k]_q! = [k]_q \cdot [k-1]_q \hdots [2]_q$.
    
    \item $f_\lambda$, which is the number of standard Young tableau of shape $\lambda$. The formula is
    \[
    f_\lambda = \frac{n!}{\prod_{b \in \lambda} h(b)}.
    \]
\end{itemize}

Table \ref{tab: lamExamples} shows these values for $n = 4$ and general $q$. From this example we can observe that $c_\lambda$ is increasing with respect to the partial order on partitions, while $n_\lambda$ is decreasing. 


\begin{table}[ht]
    \centering
    \begin{tabular}{c|c|c|c | c||c}
         $\lambda$ & $f_\lambda$ & $n_\lambda$ & $t_\lambda$ & $c_\lambda$ & $\beta_\lambda$ \\ \hline \hline
         $(4)$ & 1 & 0 & 1 & 6 & 1 \\
        $(3, 1)$ & 3 & 1 & $\frac{q(q^3 - 1)}{q - 1}$ & 2 & $\frac{q^2 - 1}{q^3 - 1}$ \\
        $(2, 2)$ & 2 & 2 & $\frac{q^2(q^4 - 1)}{q^2 - 1}$ & 0 & $\frac{q^{2} - 1}{q^4 - 1}$ \\
        $(2, 1, 1)$ & 3  & 3 & $\frac{q^3(q^3 - 1)}{q - 1}$ & -2 & $\frac{q - 1}{q^3 - 1}$\\
        $(1, 1, 1, 1)$ & 1 & 6 & $q^6$ & -6 & 0\\
    \end{tabular}
    \caption{The quantities involved in the eigenvalue and multiplicity calculations for $n = 4$.} 
    \label{tab: lamExamples}
\end{table}

Let us record that
\begin{equation}
\beta_{(n-1,1)} = \frac{q^{n-2}-1}{q^{n-1}-1},
\qquad
t_{(n-1,1)} = q\cdot \frac{q^{n-1}-1}{q-1}
\qquad\hbox{and}\qquad
f_{(n-1,1)} = n-1.
\label{topvals}
\end{equation}
Since $n!\le n^n = e^{n\log n}$ then
\begin{equation}
\frac{n!}{e^{n\log(n)}} \le 1.
\label{Stirling}
\end{equation}

We will use the following bounds from \cite[Lemma 7.2]{diaconisRam}:
\begin{equation} \label{eqn: tlamBound}
t_\lambda \le q^{-\binom{\lambda_1}{2}+\binom{n}{2}}f_\lambda, \quad \sum_{\lambda \vdash n} f_\lambda^2 = n!.
\end{equation}
In addition, we need the following bounds for sums of $f_\lambda$. Part (b) of the below proposition is Lemma 7.2(b) in \cite{diaconisRam}; the proof there is incomplete so we give the simple proof below.  
\begin{proposition} \label{prop: bounds}

\begin{enumerate}[(a)]
    \item There is a constant $K > 0$ such that for all $1 \le j \le n$,
    \[
    \sum_{\lambda: \lambda_1=n-j} f_\lambda \le \frac{n^{j}}{\sqrt{j!}} \cdot \frac{K}{j} e^{2 \sqrt{2j}}.
    \]
    \item For $1\le j\le n$, 
    \[
    \sum_{\lambda: \lambda_1=n-j} f^2_\lambda \le \frac{n^{2j}}{j!}.
    \]
\end{enumerate}

\end{proposition}

\begin{proof}
Recall that $f_\lambda$ is the number of standard tableau of shape $\lambda$, i.e.\ the elements $1, 2, \dots, n$ are arranged in the shape $\lambda$ so that rows are increasing left to right and columns are increasing top to bottom. If $\lambda_1 = n -j$, then there are $\binom{n}{j}$ ways to choose the elements not in the first row of the tableau. For a fixed partition $\lambda$, the number of ways of arranging the remaining $j$ elements is at most the number of Young tableau corresponding to the partition of $j$ created from the remaining rows of $\lambda$. This number is at most $\sqrt{j!}$ (Lemma 3 in \cite{diaconisShah1981}). Thus,
\begin{align*}
    \sum_{\lambda: \lambda_1=n-j} f_\lambda \le \binom{n}{j} \cdot \sqrt{j!} \cdot p(j),
\end{align*}
where $p(j)$ is equal to the number of partitions of $j$. It is well-known that $\log(p(n)) \sim B \cdot \sqrt{n}$ for a constant $B$. More precisely, from (2.11) in \cite{hardy}, there is a constant $K >0$ such that for all $n \ge 1$,
\[
p(n) < \frac{K}{n} e^{2 \sqrt{2n}}.
\]
This gives (a).

For part (b), we again use the inequality $f_\lambda \le \binom{n}{j} f_{\lambda^*}$, where $\lambda^* = (\lambda_2, \dots, \lambda_k)$ is the partition of $j$ determined by the rest of $\lambda$ after the first row. Then,
\begin{align*}
    \sum_{\lambda: \lambda_1 = n - j} f_\lambda^2 \le \binom{n}{j}^2 \sum_{\lambda^* \vdash j} f_{\lambda^*}^2 = \binom{n}{j}^2 \cdot j! = \left( \frac{n!}{(n - j)!j!} \right)^2 \cdot j! \le \frac{n^{2j}}{j!}.
\end{align*}
\end{proof}

\begin{proposition} \label{prop: smonotone}
The function $s(\lambda) := q^{c_\lambda} t_\lambda$ is monotone with respect to the partial order on partitions. For any $\lambda \vdash n$,
\[
s(\lambda) \le q^{\binom{n}{2}}.
\]
\end{proposition}

\begin{proof}
 
Suppose that $\lambda \prec \widetilde{\lambda}$ and $\widetilde{\lambda}$ is obtained from $\lambda$ by `moving up' one box. Suppose the box at position $(i, j)$ is moved to $(i', j')$, with $i' < i, j' > j$. 

Let $g(\lambda) = c_\lambda + n_\lambda$. Then,
    \begin{align*}
        &c_{\widetilde{\lambda}} = c_\lambda + (j' - j) + (i - i') \\
        &n_{\widetilde{\lambda}} = n_\lambda - (i - 1) + (i' - 1) = n_\lambda + (i' - i).
    \end{align*}
This implies that 
    \[
    g(\widetilde{\lambda}) = g(\lambda) + (j' - j) + (i - i') - (i - 1) + (i' - 1) = g(\lambda) + (j' - j).
    \]
Now, consider the change in $r_\lambda$. The hook-lengths of $\widetilde{\lambda}$ are:
    \begin{align*}
        &\widetilde{h}(i, j) = h(i, j) - 1 = 0, \quad \widetilde{h}(i', j') = h(i, j) + 1 = 1 \\
        &\widetilde{h}(k, j) = h(k, j) - 1, \quad k < i, k \neq i' \\
        &\widetilde{h}(i, l) = h(i, l) - 1, \quad l < j \\
        &\widetilde{h}(k, j') = h(k, j') + 1, \quad k < i' \\
        &\widetilde{h}(i', l) = h(i', l) +1, \quad l < j', l \neq j
    \end{align*}
    and $\widetilde{h}(k, l) = h(k, l)$ for all other boxes $b = (k, l)$. Thus,
    \begin{align*}
        \sum_b \widetilde{h}(b) - \sum_b h(b) &= (j' - 1) + (i' - 1) - (i - 1) - (j - 1) \\
        &= (j' - j) + (i' - i).
    \end{align*}
Using the inequalities $q^{r-1} < (q^r - 1) < q^r$, we have
\begin{align*}
    r(\widetilde{\lambda}) &= \frac{[n]_q!(q - 1)^n}{\prod_b (q^{h(b)} - 1)} \\
    &\le \frac{[n]_q!(q - 1)^n}{q^{\sum_b \widetilde{h}(b) - n}} = \frac{1}{q^{(j' - j) + (i' - i) -n} } \frac{[n]_q!(q - 1)^n}{q^{\sum_b h(b)}} \\
    &\le \frac{1}{q^{(j' - j) + (i' - i) -n} } \frac{[n]_q!(q - 1)^n}{\prod_b (q^{h(b)} - 1)} = \frac{1}{q^{(j' - j) + (i' - i) -n} } r(\lambda).
\end{align*}
Combining this with the result for $g(\lambda)$:
\begin{align*}
    s(\widetilde{\lambda}) = q^{g(\widetilde{\lambda})} r_{\lambda} = q^{g(\lambda) + (j' - j)} \cdot \frac{1}{q^{(j' - j) + (i' - i) - n}} r(\lambda) = s(\lambda) q^{i - i' + n} > s(\lambda),
\end{align*}
since $i > i'$.

Assuming the monotonicity, then if $\lambda_1 \ge n/2$ we have $s(\lambda) \le s((\lambda_1, n - \lambda_1))$. To calculate this quantity: 
\begin{align*}
    g(\lambda) &\le g((\lambda_1, n - \lambda_1)) = c_{(\lambda_1, n - \lambda_1)} + n_{(\lambda_1, n - \lambda_1)} \\
    &= \left( \frac{\lambda(\lambda_1 - 1)}{2} + \frac{(n - \lambda_1 - 1)(n - \lambda_1 - 2)}{2} - 1 \right) + (n - \lambda_1) \\
    &= \binom{\lambda_1}{2} + \binom{n - \lambda_1}{2} = \frac{n(n-1)}{2} - \lambda_1(n - \lambda_1).
    \end{align*}
For $r_\lambda$, note that the hook lengths of $(\lambda_1, n - \lambda_1)$ are 
\begin{align*}
    &\lambda_1 + 1, \lambda_1, \hdots, 2 \lambda_1 - n + 2 \\
    &2\lambda_1 -n, 2 \lambda_1 - n - 1, \dots, 3, 2, 1 \\
    &n - \lambda_1, n- \lambda_1 - 1, \dots, 3, 2, 1
\end{align*}
If $\lambda_1 = n - j$, we see the terms that cancel:
\begin{align*}
    r_{(\lambda_1, n - \lambda_1)} = \frac{[n]_q!}{\prod_b [h(b)]_q} &= \frac{(q^n - 1)(q^{n-1} - 1) \hdots (q - 1)}{\prod_b (q^{h(b)} - 1)} \\
    &=\frac{(q^n - 1)(q^{n-1} - 1) \hdots (q^{\lambda_1 + 2} - 1) \cdot (q^{2 \lambda_1 - n + 1} - 1)}{(q^{n - \lambda_1} - 1) (q^{n - \lambda -1} - 1) \hdots (q - 1)} \\
    &< \frac{q^{\sum_{k = 0}^{j -2}(n - k) + (n - 2j + 1)}}{q^{\sum_{k = 1}^j k - j}} \\
    &= q^{\frac{1}{2}(j - 1)(2(n + 1) - j) + (n - 2j + 1) - \frac{1}{2}j(j+1) + j} \\
    &= q^{j(n - j - 1) + j} = q^{j(n - j)} = q^{\lambda_1(n - \lambda_1)}.
\end{align*}
This uses the inequality $q^{r-1} < q^r - 1 < q^r$. Thus, if $\lambda_1 \ge n/2$, we have shown
\[
s(\lambda) = q^{g(\lambda)} r_\lambda \le q^{\binom{n}{2} - \lambda_1(n - \lambda_1)} \cdot q^{\lambda_1(n- \lambda_1)} = q^{\binom{n}{2}}.
\]
Now suppose $\lambda_1 \le n/2$, so $s(\lambda) \le s((n/2, n/2))$ (assume $n$ is even). To calculate this,
    \begin{align*}
        g(\lambda) &\le g((n/2, n/2)) = c_{(n/2, n/2)} + n_{(n/2, n/2)} \\
        &= \frac{(n/2 - 1)n/2}{2} + \frac{(n/2 - 2)(n/2 - 1)}{2} - 1 + \frac{n}{2} = \frac{n^2}{4} - \frac{n}{2}.
    \end{align*}
To bound $r_\lambda$, use the same calculation as before to get
\begin{align*}
r_{(n/2, n/2)} \le q^{n^2/4},
\end{align*}
and so in total $s((n/2, n/2)) \le q^{n^2/2 - n/2} = q^{\binom{n}{2}}$.
\end{proof}

\subsection{Starting from $\id$} \label{sec: identity}

\begin{theorem} \label{thm: id}
Let $P$ be the Markov chain on $S_n$ induced by random transvections on $GL_n(q)$. 
\begin{enumerate}[(a)]
    \item For $t_\lambda, \beta_\lambda$ defined in Theorem \ref{thm: evs_transvections}, we have
    \[
    4 \|P_{id}^\ell - \pi_q \|_{TV}^2 \le \chi_{id}^2(\ell) = \sum_{\lambda \vdash n, \lambda \neq (n)} t_\lambda f_\lambda \beta_\lambda^{2 \ell}.
    \]
    
    \item Let $\alpha_{n, q}$ be as in Corollary \ref{cor: alpha} and $n, q$ such that $\log(q) > 6/n$. Then if $\ell \ge \frac{n \log(q)/2  + \log(n) + c}{(\log(q) - \alpha_{n, q})}, c > 0$, we have
    \begin{align*}
        \chi_{id}^2(\ell) &\le (e^{e^{-2c}} - 1) + e^{-cn}
    \end{align*}

    \item For any $\ell \ge 1$,
    \[
    \chi_{id}^2(\ell) \ge (q^{n-1} - 1)(n - 1)q^{-4 \ell}.
    \]
    
    \item If $\ell \le n/8$, then for fixed $q$ and $n$ large,
\begin{equation} \label{eqn: lower bound}
    \|P_{id}^\ell - \pi_q \|_{TV} \ge 1 - o(1).
\end{equation}

\end{enumerate}
\end{theorem}

Theorem \ref{thm: id} shows that restricting the random transvections walk from $GL_n(q)$ to the double coset space only speeds things up by a factor of $2$, when started from the identity. Hildebrand \cite{hildebrand} shows that the total variation distance on all of $GL_n(q)$ is only small after $n + c$ steps. Note this is independent of $q$. 



\begin{proof}

\textbf{(a):} The inequality follows from Proposition \ref{prop: DRresult} (b):
\begin{align*}
     \chi^2_{id}(\ell) &= q^{-2 I(id)} \sum_{\lambda \neq (n)} t_\lambda \chi_H^\lambda(T_{id} K^{2\ell}T_{id}) \\
     &=  \sum_{\lambda \neq (n)} t_\lambda \chi_H^\lambda(K^{2\ell}) = \sum_{\lambda \neq (n)} t_\lambda f_\lambda \beta_\lambda^{2 \ell}. 
\end{align*}

\textbf{(b):} From Corollary \ref{cor: betaBound} if $\lambda_1 = n - j$, then $\beta_\lambda \le \kappa_{n, q, j} q^{-j}$, where
\begin{align*}
    \kappa_{n,q, j} &= (1 + q^{-(n - 2j + 1)})(1 + q^{-(n-1)})(1 + q^{-(n-2)}). 
\end{align*}
Using the bound on $t_\lambda$ from \eqref{eqn: tlamBound}, for $1 \le j \le \lfloor n/2 \rfloor$, we have
\begin{align}
\sum_{\lambda: \lambda_1 = n - j} t_\lambda f_\lambda \beta_\lambda^{2 \ell} & \le  (\kappa_{n, q, j} q^{-j} )^{2\ell} \sum_{\lambda: \lambda_1 = n -j} q^{-\binom{n-j}{2} + \binom{n}{2}} f_\lambda^2 \notag \\
&\le (\kappa_{n, q, j} q^{-j} )^{2\ell} q^{nj - j(j+1)/2}   \cdot \frac{n^{2j}}{j!} \notag \\
&\le \frac{1}{j!} \exp \left( -2 \ell j \left(\log(q) - \frac{\log(\kappa_{n, q, j})}{j} \right) + \log(q) (nj - j(j + 1)/2) + 2j \log(n)    \right) \notag \\
&\le \frac{1}{j!} \exp \left( -2 \ell j \left(\log(q) - \alpha_{n, q} \right) + \log(q) (nj - j(j + 1)/2) + 2j \log(n)    \right) \label{eqn: idExponent}.
\end{align}

Recall the final inequality follows since $\alpha_{n, q} := \max_{1 \le j \le n/2} \log(\kappa_{n, q, j})/j$. If $\ell = \frac{n \log(q)/2  + \log(n) + c}{(\log(q) - \alpha_{n, q})}$, then the exponent in Equation \eqref{eqn: idExponent} is
\begin{align*}
- 2j \left( n\log(q)/2 + \log(n) + c \right) &+ \log(q) (nj - j(j + 1)/2) + 2j \log(n) \\
&= -j(2c + \log(q)(j+1)/2).
\end{align*}
This gives
\[
\sum_{j = 1}^{\lfloor n/2 \rfloor}\sum_{\lambda: \lambda_1 = n - j} t_\lambda f_\lambda \beta_\lambda^{2 \ell} \le  \sum_{j = 1}^{\lfloor n/2 \rfloor} \frac{e^{-2jc }}{j!} \le e^{e^{-2c}} - 1.
\]
Next we need to consider the partitions $\lambda$ with $\lambda_1 \le n/2$. For these partitions,
\[
\beta_{\lambda} \le \beta_{(n/2, n/2)} \le \kappa_{n, q, n/2} q^{-n/2},
\]
Then we have
\begin{align*}
 \sum_{j = n/2}^{n-1}  \sum_{\lambda: \lambda_1 = n -j} t_\lambda f_\lambda \beta_\lambda^{2 \ell} &\le  \beta_{(n/2, n/2)}^{2\ell} \sum_{j = n/2}^{n-1}  \sum_{\lambda: \lambda_1 = n -j} t_\lambda f_\lambda \\
   &\le  \beta_{(n/2, n/2)}^{2\ell} \sum_{j = n/2}^{n-1}  \sum_{\lambda: \lambda_1 = n -j} q^{nj - j(j+1)/2} f_\lambda^2 \\
    &\le n! \cdot \exp \left( -2 \ell (n/2) \left( \log(q) + \frac{\kappa_{n, q, n/2}}{n/2} \right) +  n^2/2 \right) \\
    &\le \exp \left( -n (n \log(q)/2  + \log(n) + c) +   n^2/2 \right) \sum_{\lambda \vdash n} f_\lambda^2 \\
    &\le n! \exp \left( -n (n \log(q)/2  + \log(n) + c) +   n^2/2 \right),
\end{align*}
using $\sum_{\lambda} f_\lambda^2 = n!$ and that if $n/2 \le j \le n$, then
\[
nj - j(j + 1)/2 \le n^2 - n^2/4 < n^2 - n(n + 1)/2 \le n^2/2,
\]
since the function is increasing in $j$. Note also that if $q \ge 3$ then $\log(q) > 1$ and $n! \le n^n = e^{\log(n)}$. To finish the bound,
\begin{align*}
      \sum_{j = n/2}^{n-1}  \sum_{\lambda: \lambda_1 = n -j} t_\lambda f_\lambda \beta_\lambda^{2 \ell}   &\le \exp\left( n \log(n) -n (n \log(q)/2  + \log(n) + c) +   n^2/2 \right) \\
    &= \exp\left(  -(n^2/2) (\log(q) - 1) - cn \right) \le e^{-cn}.
\end{align*}

\textbf{(c):} The lower bound comes from considering the $\lambda = (n - 1, 1)$ term from the sum in (a). Using the quantities \eqref{topvals}, this gives
\begin{align*}
    \chi^2_{id}(\ell) &= \sum_{\lambda \neq (n)} t_\lambda f_\lambda \beta_\lambda^{2\ell} \ge t_{(n-1, 1)} f_{(n-1, 1)} \beta_{(n-1, 1)}^{2 \ell} \\
    &= q \cdot \frac{q^{n-1} - 1}{q -1} (n-1) \left( \frac{q^{n-2} - 1}{q^{n-1} - 1} \right)^{2 \ell} \ge (q^{n-1} - 1)(n-1) q^{-4\ell}.
\end{align*}
This uses that $(q^{n-2} - 1)/(q^{n-1} - 1) \ge q^{-2}$.

\textbf{(d):} From the alternative version of the walk on the Hecke algebra, involving $D/|\mathcal{T}_{n, q}|$ with $D$ from \ref{thm: D}, the walk proceeds by picking a transposition $(i,j), i < j$ with probability proportional to
\[
q^{-(j-i)}
\]
and multiplying by $T_{ij}$. As described in Section \ref{sec: metropolis}, multiplication by $T_{ij}$ corresponds to proposing the transposition $(i, j)$ and proceeding via the Metropolis algorithm. Thus, multiplication by $T_{ij}$ induces at most $2(j - i)$ inversions, always less than $2n$. From \cite{diaconisShah1981} Theorem 5.1, under $\pi_q$, a typical permutation has $\binom{n}{2} - \frac{n - q}{q + 1} + O(\sqrt{n})$ inversions (and the fluctuations are Gaussian about this mean). If $\ell = n/8$, the measure $P_{id}^\ell(\cdot)$ is concentrated on permutations with at most $n^2/4$ inversions and $\pi_q$ is concentrated on permutations with order $n^2/2 - (n - 1)/(q + 1) + O(\sqrt{n})$ inversions. 
\end{proof}

\subsection{Starting from $\omega_0$} \label{sec: reversal}

\begin{theorem} \label{thm: reversal}
Let $P$ be the Markov chain on $S_n$ induced by random transvections on $GL_n(q)$, and let $\omega_0 \in S_n$ be the reversal permutation in $S_n$.
\begin{enumerate}[(a)]
    \item With $t_\lambda, c_\lambda, \beta_\lambda$ defined in Section \ref{sec: evs},
    \[
    4 \|P_{\omega_0}^{\ell} - \pi_q \|_{TV}^2 \le \chi_{\omega_0}^2(\ell) = q^{-\binom{n}{2}} \sum_{\lambda \neq (n)} q^{c_{\lambda}} t_\lambda f_\lambda \beta_\lambda^{2 \ell}.
    \]
    
    \item Let $\alpha_{n, q}$ be as in Corollary \ref{cor: alpha} and $n, q$ such that $\log(q) > 6/n$. If $\ell \ge (\log(n)/2 + c)/(\log(q) - \alpha_{n, q})$ for $c > 0$ with $c \ge 2 \sqrt{2}$ then 
    \[
    \chi_{\omega_0}^2(\ell) \le - 2K \log(1 - e^{-c}) + \sqrt{K} e^{-nc},
    \]
    for a universal constant $K > 0$ (independent of $q, n$).
    
    \item For any $\ell \ge 1$,
    \[
    \chi^2_{\omega_0}(\ell) \ge q^{-(n-2)}(n-1) (q^{n-1} -1) q^{-4 \ell}.
    \]
    
\end{enumerate}
\end{theorem}

\begin{remark}
Theorem \ref{thm: reversal} shows that the Markov chain has a cutoff in its approach to stationarity in the chi-square metric. It shows the same exponential speed up as the walk started at a typical position (Theorem \ref{thm: typical} below), and indeed is faster by a factor of $2$. This is presumably because it starts at the permutation $\omega_0$, at which the stationary distribution $\pi_q$ is concentrated, instead of `close to $\omega_0$'.
\end{remark}

\begin{proof}[Proof of Theorem \ref{thm: reversal}]
\textbf{(a):} By Proposition 4.9 in \cite{diaconisRam}, if $\omega_0$ is the longest element of $W$ then 
\[
\rho^\lambda(T_{\omega_0}^2) = q^{I(\omega_0) + c_\lambda}\mathrm{Id},
\]
where $\rho^\lambda$ is the irreducible representation indexed by $\lambda$. Using this and \ref{prop: DRresult} (a),
\begin{align*}
    \chi^2_{\omega_0}(\ell) &=  q^{-2 I(\omega_0)}\sum_{\lambda \neq \mathbf{1}} t_\lambda \chi_H^\lambda(T_{\omega_0^{-1}} K^{2 \ell} T_{\omega_0}) \\
    &=  q^{-2 I(\omega_0)}\sum_{\lambda \neq \mathbf{1}} t_\lambda \chi_H^\lambda( K^{2 \ell} T_{\omega_0} T_{\omega_0^{-1}}) \\
    &= q^{-2 I(\omega_0)}\sum_{\lambda \neq \mathbf{1}} t_\lambda \chi_H^\lambda( K^{2 \ell} T_{\omega_0}^2) \\
     &= q^{-2 I(\omega_0)}\sum_{\lambda \neq \mathbf{1}} t_\lambda \chi_H^\lambda( K^{2 \ell} q^{c_\lambda + I(\omega_0)}) \\
     &= q^{-I(\omega_0)} \sum_{\lambda \neq \mathbf{1}} q^{c_\lambda} t_\lambda \chi_H^\lambda( K^{2 \ell} ) = q^{-I(\omega_0)} \sum_{\lambda \neq \mathbf{1}} q^{c_\lambda} t_\lambda f_\lambda \beta_\lambda^{2 \ell},
\end{align*}
since $K \in Z(H)$, i.e.\ $K$ commutes with all elements of the Hecke algebra.

\textbf{(b):} Suppose $\lambda_1 = n -j$ for $1 \le j \le n/2$.  Recall the definition $s(\lambda) = q^{c_\lambda}t_\lambda$. From Proposition \ref{prop: smonotone}, $s(\lambda) \le q^{\binom{n}{2}}$. Then,
\begin{align*}
    q^{-\binom{n}{2}} \sum_{\lambda: \lambda_1 = n - j} q^{c_{\lambda}} t_\lambda f_\lambda \beta_\lambda^{2 \ell} & = q^{-\binom{n}{2}} \sum_{\lambda: \lambda_1 = n - j} s(\lambda) f_\lambda \beta_\lambda^{2 \ell} \\
    &    \le  (\kappa_{n, q, j} q^{-j} )^{2\ell} q^{-\binom{n}{2}} \sum_{\lambda: \lambda_1 = n -j} q^{ \binom{n }{2} }  f_\lambda \\
& \le  (\kappa_{n, q, j} q^{-j} )^{2\ell} \sum_{\lambda: \lambda_1 = n -j}  f_\lambda \\
&\le (\kappa_{n, q, j} q^{-j} )^{2\ell}   \cdot \frac{n^{j}}{\sqrt{j!}} \cdot \frac{K}{j}e^{2 \sqrt{2j}} \\
&\le \frac{K}{\sqrt{j!}} \exp \left( -2 \ell j \left(\log(q) - \frac{\log(\kappa_{n, q, j})}{j} \right)  + (j-1)\log(n)  + 2 \sqrt{2j}  \right).
\end{align*}
The third inequality uses Proposition \ref{prop: bounds} for $\sum_{\lambda: \lambda_1 = n - j} f_\lambda$. Recall $\alpha_{n, q} := \max_{1 \le j \le n/2} \log(\kappa_{n, q, j})/j$. If $\ell = (\log(n)/2 + c)/(\log(q) - \alpha_{n, q})$, then the bound becomes
\begin{align*}
   \sum_{j = 1}^{n/2} q^{-\binom{n}{2}} \sum_{\lambda: \lambda_1 = n - j} q^{c_{\lambda}} t_\lambda f_\lambda \beta_\lambda^{2 \ell} &\le \sum_{j = 1}^{n/2}\frac{K}{\sqrt{j!}} \exp \left( - j \log(n) - 2jc   + (j-1)\log(n)  + 2 \sqrt{2j}  \right) \\
   &\le  2 K \sum_{j = 1}^{n/2} \frac{e^{-2jc + 2 \sqrt{2j}}}{j}, 
\end{align*}
using the loose bound $\sqrt{j!} > j/2$ for all $j \ge 1$. With the assumption that $c \ge 2 \sqrt{2}$, we have $-2jc + 2 \sqrt{2j} \le -jc$ for all $j \ge 1$. Finally,
\begin{align*}
     2 K \sum_{j = 1}^{n/2} \frac{e^{-jc}}{j} & \le  2 K \sum_{j = 1}^{\infty} \frac{e^{-jc}}{j} \\
     &= - 2K \log(1 - e^{-c}).
\end{align*}

Now for the $\lambda$ with $\lambda_1 \ge n/2$, we have
\begin{align*}
      \sum_{j = n/2}^{n-1}  \sum_{\lambda: \lambda_1 = n -j}  q^{-\binom{n}{2}}  q^{c_{\lambda}} t_\lambda f_\lambda \beta_\lambda^{2 \ell} &\le  \beta_{(n/2, n/2)}^{2\ell} \sum_{j = n/2}^{n-1}  \sum_{\lambda: \lambda_1 = n -j}  f_\lambda \\
      &\le  \beta_{(n/2, n/2)}^{2\ell} \sum_{\lambda \vdash n}  f_\lambda \\
   &\le  \beta_{(n/2, n/2)}^{2\ell} \left(\sum_{\lambda \vdash n}  f_\lambda^2\right)^{1/2} \cdot p(n)^{1/2},
\end{align*}
where $p(n)$ is equal to the number of partitions of $n$ (the inequality is Cauchy-Schwartz). Since $p(n) \le \frac{K}{n} e^{2 \sqrt{2n}}$ for a constant $K > 0$ (\cite{hardy}) and $\sum_{\lambda \vdash n} f_\lambda^2 = n!$, this gives 
\begin{align*}
        &\le \sqrt{n!} \cdot \sqrt{\frac{K}{n}} \exp \left( \sqrt{2 n} -2 \ell (n/2) \left( \log(q) - \frac{\kappa_{n, q, n/2}}{n/2} \right) \right) \\
    &\le \sqrt{n!} \cdot \sqrt{\frac{K}{n}} \exp \left(\sqrt{2 n} -n \log(n)/2 -2n c \right) \\
    &\le \sqrt{K} \exp \left( -nc - n(c  - \sqrt{2}n^{-1/2}) \right),
\end{align*}
since $\sqrt{n!} \le e^{n \log(n)/2}$. Since $c > 2 \sqrt{2}$, then the bound is $\le \sqrt{K} e^{-nc}$ for any $n \ge 1$.


\textbf{(c):} A lower bound comes from using \eqref{topvals} for the lead term on the right hand side of (a):
\begin{align*}
    \chi_{\omega_0}^2(\ell) &\ge  q^{-\binom{n}{2}}  q^{c_{(n-1, 1)}} t_{(n-1, 1)} f_{(n-1, 1)} \beta_{(n-1, 1)}^{2 \ell} \\
    &= q^{-\binom{n}{2}}  q^{\binom{n-1}{2} + 1} \cdot q \frac{q^{n-1} - 1}{q - 1} \cdot (n-1) \left( \frac{q^{n-2} - 1}{q^{n-1} - 1} \right)^{2 \ell} \\
    &\ge q^{-(n-2)}(q^{n-1} -1)(n-1) q^{-4 \ell}.
\end{align*}
\end{proof}

\subsection{Starting from a Typical Site} \label{sec: typical}

In analyzing algorithms used repeatedly for simulations, as the algorithm is used, it approaches stationarity. This means the quantity
\[
\sum_{x \in S_n} \pi_q(x) \|P_x^\ell - \pi_q \|_{TV}
\]
is of interest. For the problem under study, $\pi_q$ is concentrated near $\omega_0$ so we expect rates similar to those in Theorem \ref{thm: reversal}.

\begin{theorem} \label{thm: typical}
Let $P$ be the Markov chain on $S_n$ induced by random transvections on $GL_n(q)$. 
\begin{enumerate}[(a)]
    \item With $f_\lambda, \beta_\lambda$ defined in Section \ref{sec: evs},
    \[
    \left( \sum_{x \in S_n} \pi_q(x) \|P_x^\ell - \pi_q \|_{TV} \right)^2 \le \frac{1}{4} \sum_{x \in S_n} \pi_q(x) \chi_{x}^2(\ell) = \frac{1}{4} \sum_{\lambda \neq (n)} f_\lambda^2 \beta_\lambda^{2 \ell}.
    \]
    
    \item Let $\alpha_{n, q}$ be as in Corollary \ref{cor: alpha} and $n, q$ such that $\log(q) > 6/n$. If $\ell \ge (\log(n) + c)/(\log(q) - \alpha_{n, q}), c > 0$, then
    \[
    \sum_{x \in S_n} \pi_q(x) \chi_{x}^2(\ell) \le (e^{e^{-c}} - 1) + e^{-cn}.
    \]
    
    \item For any $\ell \ge 1$, 
    \[
    \sum_{x \in S_n} \pi_q(x) \chi_{x}^2(\ell) \ge (n-1)^2 q^{-4 \ell}.
    \]
\end{enumerate}
\end{theorem}

\begin{proof}
\textbf{(a):} This is simply a restatement of \ref{prop: DRresult} part (b).


\textbf{(b):} 
We will divide the sum depending on the first entry of the partition. By Proposition \ref{prop: bounds},  we have the bound (true for any $1 \le j \le n$)
\[
\sum_{\lambda: \lambda_1 = n - j} f_\lambda^2 \le \frac{n^{2j}}{j!}.
\]
Combining this with Corollary \ref{cor: betaBound},   for $j \le \lfloor n/2 \rfloor$, 
\begin{align*}
    \sum_{\lambda: \lambda_1 = n - j} \beta_\lambda^{2 \ell} f_\lambda^2 &\le \sum_{\lambda: \lambda_1 = n - j} \beta_{(n - j, j)}^{2 \ell} f_\lambda^2 \le \beta_{(n - j, j)}^{2 \ell} \cdot \frac{n^{2j}}{j!} \\
    &\le \kappa_{n, q, j}^{2\ell} q^{-2 \ell j} \cdot \frac{n^{2j}}{j!} = \frac{1}{j!} \exp \left( 2 \ell( \log(\kappa_{n, q, j}) - j\log(q)) + 2 j \log(n) \right) \\
    &= \frac{1}{j!} \exp \left( 2 \ell j \left( \frac{ \log(\kappa_{n, q, j})}{j} - \log(q)\right) + 2 j \log(n) \right)
\end{align*}
 Define
\[
\alpha_{n, q} := \max_{1 \le j \le n/2} \frac{\log(\kappa_{n, q, j})}{j},
\]
so then if $\ell = (\log(n) + c)/(\log(q) - \alpha_{n, q})$, the bound is
\begin{align*}
    &\le \frac{1}{j!} \exp \left( 2 \ell j( \alpha_{n, q} - \log(q)) + 2 j \log(n) \right) \\
    & = \frac{1}{j!} \exp \left( -2 j \log(n) - 2j c + 2j \log(n) \right) = \frac{e^{-2jc}}{j!}.
\end{align*}

Then summing over all possible $j$ gives
\[
\sum_{j = 1}^{\lfloor n/2 \rfloor} \sum_{\lambda: \lambda_1 = n - j} f_\lambda^2 \beta_\lambda^{2 \ell} \le \sum_{j = 1}^{\lfloor n/2 \rfloor} \frac{e^{-2jc}}{j!} \le \left( e^{e^{-c}} - 1 \right).
\]

Now we have to bound the contribution from partitions $\lambda$ with $\lambda_1 \le n/2$. Because $\beta_\lambda$ is monotone with respect to the order on partitions, we have for all $\lambda$ such that $\lambda_1 \le n/2$, 
\[
\beta_{\lambda} \le \beta_{(n/2, n/2)} \le \kappa_{n, q, n/2} q^{-n/2},
\]
since $\lambda \preceq (n/2, n/2)$ (assuming without essential loss that $n$ is even). Then,
\begin{align*}
    \sum_{j = n/2}^{n-1}  \sum_{\lambda: \lambda_1 = n -j}f_\lambda^2 \beta_\lambda^{2 \ell} &\le \beta_{(n/2, n/2)}^{2\ell} \sum_{\lambda} f_\lambda^2 \\
    &\le n! \cdot \exp \left( 2 \ell (n/2)( \log(q) + \frac{\kappa_{n, q, n/2}}{n/2} ) \right) \\
    &\le n! \cdot \exp \left( -n (\log(n) + c)  \right) \\
    &\le e^{-cn},
\end{align*}
using $\sum_{\lambda} f_\lambda^2 = n! \le n^n$.


\textbf{(c):} 
The sum is bounded below by the term for $\lambda = (n - 1, 1)$. This is
\[
f_{(n-1, 1)}^2 \cdot \beta_{(n-1, 1)}^{2 \ell} = (n-1)^2 \left( \frac{q^{n-2} - 1}{q^{n-1} - 1} \right)^{2 \ell} \ge (n-1)^2 q^{-4 \ell}.
\]

\end{proof}

\newpage

\section{Hecke Algebra Computations} \label{sec: heckeComputations}

This section proves Theorem \ref{thm: D} which describes the transvections Markov chain on $S_n$ as multiplication in the Hecke algebra from Definition \ref{def: IwahoriHecke}. This is accomplished through careful row reduction.

\subsection{Overview }

Let $\CC[G]$ denote the group algebra for $G = GL_n(q)$. This is the space of functions $f: G \to \CC$, with addition defined $(f + g)(s) = f(s) + g(s)$ and multiplication defined by
\[
    f_1 \ast f_2(s) = \sum_{t\in G} f_1(t) f_2(st^{-1}).
\]
Equivalently, $\CC[G] = \mathrm{span} \{g \mid g \in G \}$ and we can write an element $f = \sum_{g} c_g g$ for $c_g \in \CC$, so $f(g) = c_g$.

Define elements in $\CC[G]$:
    \[
    \lB = \frac{1}{|\cB|} \sum_{x \in \cB} x, \quad T_\omega = \frac{1}{|\cB|} \sum_{x \in \cB\omega \cB} x, \,\,\,\,\, \omega \in S_n.
    \]
Note that if $b \in \cB$, then
    \[
    b\lB = \lB b = \lB, \quad \lB^2 = \lB.
    \]
    
If $g \in \cB \omega \cB$, so $g = b_1 \omega b_2$, then
    \begin{align*}
        \lB g \lB  &= \lB b_1 \omega b_2 \lB = \lB \omega \lB \\
        &= \frac{1}{|\cB|^2} \sum_{b_1, b_2 \in \cB} b_1 \omega b_2 \\
        &= \frac{1}{|\cB|^2} \frac{|\cB|^2}{|\cB \omega \cB|} \sum_{x \in \cB \omega \cB} x \\
        &= \frac{|\cB|}{|\cB \omega \cB|}  \frac{1}{|\cB|}  \sum_{x \in \cB \omega \cB} x = \frac{|\cB|}{|\cB \omega \cB|}  \frac{1}{|\cB|}  T_\omega = q^{-I(\omega)}T_\omega.
    \end{align*}
    
The Hecke algebra is $H_n(q) = \lB \CC[G] \lB$ and has basis $\{T_\omega \mid \omega \in S_n \}$. Note that $H_n(q)$ are all functions in $\C[G]$ which are $B-B$ invariant, i.e.\ $f(b_1 g b_2) = f(g)$.

Now let $P$ be the transition matrix for $G$ defined by multiplying by a random transvection. We can also write this
    \[
    P = \frac{1}{|\T_{n, q}|} \sum_{T \in \T_{n, q}} M_T,
    \]
    where $M_T$ is the transition matrix `multiply by $T$'. In other words, $M_T(x, y) = \textbf{1}(y = Tx)$.
    
Then $P$ defines a linear transform on the space $\C[G]$, with respect to the basis $\{g \mid g \in G \}$. The matrix $M_T$ is just the function: Multiply by $T$ in the group algebra. This means $P$ is equivalent to multiplication by $\frac{1}{|\T_{n, q}|} D$, with $D = \sum_{T \in \T_{n, q}} T$ as an element in $\C[G]$. The Markov chain lumped to $S_n = \cB \backslash GL_n(q) / \cB$ is then equivalent to multiplication by $D$ on $H_n(q)$. Since $D$ is the sum of all elements in a conjugacy class, it is in the center $Z(\CC[G])$. This means if $g \in B \omega B$, then 
    \[
    \lB (D g)  \lB= D  \lB  \omega  \lB = (D  \lB) q^{-I(\omega)}T_\omega.
    \]
In conclusion, to determine how $D$ acts in $H_n(q)$, we can compute $D \lB$. The remainder of the section proves the following.

\begin{theorem} \label{thm: secondD}
Let $D = \sum_{T \in \T_{n, q}} T \in \CC[G]$. Then,
\[
D \lB = \left((n - 1)q^{n-1} - [n-1]_q\right) \lB + (q-1) \sum_{1 \le i < j \le n} q^{n - 1 - (j-i)} T_{ij}.
\]
\end{theorem}

\subsection{Row Reduction}
Let $G = GL_n(q)$ and $\cB$ be the upper triangular matrices. For $1 \le i \le n -1$ and $c \in \FF_q$, define $y_i(c) \in G$ by
\[
y_i(c) = \begin{blockarray}{ccccccccc}
 &  & & & i & i + 1 &  & &  \\
\begin{block}{c(cccccccc)}
 & 1 & & &  &  &  & & \\
  &  & \ddots & &  &  &  & & \\
   &  & & 1 &  &  &  & & \\
  i  &  & &  & c & 1 &  & & \\
 i+1  &  & & & 1 & 0  &  & & \\
     &  & & &  &  & 1 & & \\
         &  & & &  &  &  & \ddots & \\
          &  & & &  &  &  & & 1 \\
\end{block}
\end{blockarray}.
\]
That is, multiplication on the left by $y_i(c)$ acts by adding $c$ times the $i$th row to the $(i+1)$th row, then permuting the $i$ and $i + 1$ rows.

Let $\omega \in S_n$. We can write the reduced word $\omega = s_{i_1} \hdots s_{i_\ell}$, for $s_i = (i, i+1)$ the simple reflections. Then,
\[
\cB \omega \cB = \left\lbrace \cB y_{i_1}(c_1) \hdots y_{i_\ell}(c_\ell) \cB \mid c_1, \dots, c_\ell \in \FF_q \right\rbrace,
\]
and $|\cB \omega \cB|/|\cB| = q^{I(\omega)}$ and $G = \bigsqcup_{\omega \in S_n} \cB \omega \cB$.

This provides a very useful way for determining the double coset that a matrix $M$ belongs to, which just amounts to performing row reduction by mulitplying by different matrices $y_i(c)$.


\subsection{Transvections}

For every transvection $T \in \T_{n, q}$ we can perform row operations, multiplying by $y_i(c)$, to determine which double coset $T$ belongs to. 

A transvection is defined by two vectors $\ma, \mv \in \FF_q^n$ with $\mv^\top \ma = 0$, with the last non-zero entry of $\mv$ normalized to be $1$. Assume this entry is at position $j$. Assume the first non-zero entry of $\ma$ is at position $i$. That is, the vectors look like:
\begin{align*}
    &\ma = (0, \dots, 0, a_i, a_{i+1}, \dots, a_n)^\top, \quad a_i \neq 0, \\
    &\mv = (v_1, \dots, v_{j-1}, 1, 0, \dots, 0)^\top.
\end{align*}
Let $T_{\ma, \mv} = Id + \mv \ma^\top$. We consider the possible cases for $i$ and $j$ to prove the following result.

\begin{proposition} \label{prop: TCases}
Let $T_{\ma, \mv} = I + \mv \ma^\top$ be the transvection defined by non-zero vectors $\ma, \mv$ with $\mv^\top \ma = 0$ and the last nonzero entry of $\mv$ equal to $1$. If $j > i$, then
$$T_{a,v} \in \cB s_{j-1}\cdots s_{i+1}s_is_{i+1}\cdots s_{j-1} \cB$$
exactly when the last nonzero entry of $\mv$ is $v_j$ and the first nonzero entry of $a$ is $\ma_i$.
\end{proposition}

\begin{proof}
\textbf{Case 1 $i > j$:} If $i > j$, then $T_{\ma, \mv} \in \cB$. To see this, suppose $k > l$. Then
\[
T_{\ma, \mv}(k, l) = (\mv \ma^\top)_{k, l} = v_k a_l = 0,
\]
because $v_k a_l$ can only be non-zero if $k \le j < i \le l$. 
Now, we can count how many transvections satisfy this. There are $(q - 1)q^{n - i}$ choices for $\ma$ (because $a_i$ must be non-zero) and $q^{j-1}$ choices for $\mv$ (because $v_j$ is fixed at $1$). In total for this case, there are
\begin{align*}
\sum_{j = 1}^{n-1} \sum_{i = j + 1}^n (q - 1) q^{n - i} q^{j - 1} &= (q - 1) \sum_{j = 1}^{n-1} q^{j-1}(q^{n - j - 1} + q^{n - j - 2} + \hdots + q + 1 ) \\
&= (q - 1) \sum_{j = 1}^{n-1} q^{j-1} \frac{q^{n - j} - 1}{q - 1} = \sum_{j = 1}^{n- 1} (q^{n - 1} - q^{j - 1}) \\
&=(n - 1)q^{n-1} - (1 + q + \hdots + q^{n - 1}) = (n - 1)q^{n-1} - [n-1]_q.
\end{align*}

\textbf{Case 2 $i = j$:} In this case, $\mv^\top \ma = a_i \neq 0$. This does not satisfy the condition $\mv^\top \ma = 0$, so this case cannot occur.

\textbf{Case 3 $i < j$:} If $i < j$, then the transvection is of the form 
\[
T_{\ma, \mv} = \begin{blockarray}{cccccccc}
 & 1 & \hdots & j & \hdots & i & \hdots & n \\
\begin{block}{c(cc|ccc|cc)}
 1 &  &  &  &  &  &  & \\
 \vdots &  & Id &  & \ast &  & \ast & \\ \cmidrule(lr){2-8}
  i &  &  &  &  &  &   & \\ 
   \vdots &  &  &  & N_{\ma, \mv} &  & \ast & \\
 j &  &  &  &  &  &  & \\ \cmidrule(lr){2-8}
  \vdots &  &  &  &  &  & Id & \\
  n &  &  &  &  &  &  &  \\
\end{block}
\end{blockarray}.
\]
Where
\[
N_{\ma, \mv} = \begin{pmatrix} 1+a_iv_i &a_{i+1}v_i &a_{i+2} v_i &\cdots &a_jv_i \\
a_iv_{i+1} &1+a_{i+1}v_{i+1} & a_{i+2} v_{i+1} &\cdots &a_j v_{i+1}\\
\vdots &&&&\vdots \\
a_i &a_{i+1} &&\cdots &1+ a_j
\end{pmatrix}
\]
Let $b$ be the $j - i + 1\times j - i + 1$ matrix
\[
b = \begin{pmatrix}
a_i & a_{i+1} & \hdots & a_{j-1} &-1 - a_j \\
0 & 1 & \hdots  &0 &-v_{i+1}  \\
0 &0 &\ddots &  & \vdots \\
0 &0 &0 &1 &-v_{j-1} \\
0 &0 &0 &0 &-a_i^{-1}
\end{pmatrix}
\]
Then,
\[
T_{\ma, \mv} = y_{j-1}(v_{j-1}) \hdots y_{i+1}(v_{i+1}) \cdot y_i(1 + a_i^{-1}) \cdot y_{i+1}(-a_{i+1}a_i^{-1}) \hdots y_{j-1}(-a_{j-1}a_{i}^{-1}) \cdot A,
\]
where $A$ is the upper-triangular matrix
\[
A = \begin{blockarray}{cccc}
 &  &  & \\
\begin{block}{c(c|c|c)}
 & Id & \ast & \ast \\ \cmidrule(lr){2-4}
  & 0 & b & \ast \\ \cmidrule(lr){2-4}
& 0 & 0 & Id \\
\end{block}
\end{blockarray}.
\]
See Appendix \ref{app: rowReduction} for details of this row reduction calculation.

To count how many transvections fit this case:
\begin{itemize}
\item  $q^{n-i}$ choices for $a_{i+1}, \dots, a_n$.
\item $q - 1$ choices of $a_i$
\item $q^{i-1}$ choices of $v_1, \dots, v_{i-1}$.
\item $q^{j - 1 - i}$ choices of $v_i, \dots, v_{j-1}$ to satisfy $a_i v_i + \hdots + v_{j-1} a_{j-1} + a_j = 0$.
\end{itemize}
The total is
\begin{align*}
q^{n - i}(q - 1)q^{i - 1}q^{j - 1 - i} &= (q - 1)q^{n - i + j - 2} \\
&= (q - 1)q^{n - 1 - (j - i)}q^{2(j - 1 - i) + 1} \\
&= (q - 1)q^{n - 1 - (j - i)}q^{I(s_{ij})}
\end{align*}

\end{proof}

Proposition \ref{prop: TCases} now enables the proof of Theorem \ref{thm: secondD}.

\begin{proof}[Proof of Theorem \ref{thm: secondD}]
Let $\CC[G]$ be the group algebra of $G = GL_n(q)$, and $Z(\CC[G])$ the center. Then since transvections are a conjugacy class, the sum of transvections is in the center,
\[
D = \sum_{\ma, \mv} T_{\ma, \mv} \in Z(\CC[G]).
\]
Since $D$ commutes with every element of $\CC[G]$, we can compute,
\begin{align*}
D \lB &= D \lB^2 = \lB D \lB \\
&= \sum_{\ma, \mv} \lB T_{\ma, \mv} \lB \\
&= \sum_{ \substack{\ma, \mv: \\ i > j}} \lB T_{\ma, \mv} \lB + \sum_{\substack{\ma, \mv: \\ i < j}} T_{\ma, \mv} \lB \\
&= \sum_{\substack{\ma, \mv: \\ i > j}}  \lB + \sum_{\substack{\ma, \mv: \\ i < j}} \lB s_{ij} \lB \\
&= ((n-1) q^{n-1} - [n-1]_q ) \lB + \sum_{1 \le i < j \le n} (q - 1) q^{n - 1 - (j - i)} q^{I(s_{ij})} \lB s_{ij} \lB \\
&= ((n-1) q^{n-1} - [n-1]_q ) \lB + (q-1) \sum_{1 \le i < j \le n}  q^{n - 1 - (j - i)} T_{s_{ij}}.
\end{align*}

\end{proof}


\section{Symmetric Function Calculation} \label{sec: symmetricFunction}
Section \ref{sec: symFunc} introduces symmetric functions which will be used in Section \ref{sec: characterComputation} to compute $D$, the element in the center of the group algebra $\C[GL_n(q)]$ which corresponds to the conjugacy class of transvections. This gives a different proof of Theorem \ref{thm: secondD}. Sections \ref{sec: n=2example} and \ref{sec: n=3example} check and motivate the computation by computing $D$ explicitly for $n = 2$ and $n = 3$.

\subsection{Symmetric Functions} \label{sec: symFunc}

For $\lambda = (\lambda_1, \ldots, \lambda_n)\in \ZZ_{\ge 0}^n$ with
$\lambda_1\ge \cdots \ge \lambda_n\ge 0$  
define the Hall-Littlewood polynomial by
$$P_\lambda = P_\lambda(x;0,t) 
= \frac{1}{v_\lambda(t)} \sum_{w\in S_n} w
\Big(x_1^{\lambda_1}\cdots x_n^{\lambda_n} \prod_{1\le i<j\le n} \frac{x_i-tx_j}{x_i-x_j} \Big)
$$
where $v_\lambda(t)$ is the normalization that make the coefficient of $x^\lambda$ in 
$P_\lambda$ equal to $1$.  For $\lambda = (1^{m_1}2^{m_2}\ldots)$ define
$$Q_\lambda = b_\lambda(t) P_\lambda
\qquad\hbox{where}\qquad
b_\lambda(t) = \prod_{i=1}^\ell \varphi_{m_i}(t)
\quad\hbox{with}\quad
\varphi_m(t) = (1-t)(1-t^2)\cdots (1-t^m).
$$
The \emph{monomial symmetric functions} $m_\lambda$ and 
the \emph{Schur functions} $s_\lambda$ are given by 
$$m_\lambda = P_\lambda(x;0,1)
\qquad\hbox{and}\qquad
s_\lambda = P_\lambda(x;0,0),
$$
the evaluations of  $P_\lambda$ at $t=1$ and $t=0$, respectively.
For $r\in \ZZ_{\ge 0}$ define $q_r$ by the generating function
$$\sum_{r\in \ZZ_{\ge 0}} q_rz^r = \prod_{i=1}^n \frac{1-tx_iz}{1-x_iz}
\qquad\hbox{and define}\qquad
q_\nu = q_{\mu_1}\cdots q_{\mu_n},
\qquad
\hbox{for $\nu = (\nu_1, \ldots, \nu_n)\in \ZZ_{\ge 0}^n$.}
$$
The \emph{Big Schurs} are defined by the $n\times n$ determinant (see \cite[Ch.\ III (4.5)]{macdonald})
$$S_\lambda = S_\lambda(x;t) = \det( q_{\lambda_i - i+j}).$$
Define $a_{\mu\nu}(t)$, $K_{\lambda\mu}(0,t)$ and $L_{\nu\lambda}(t)$ by
$$Q_\mu(x;0,t) = \sum_\nu a_{\mu\nu}(t) m_\nu
\qquad Q_\mu(x;0,t) = \sum_\lambda K_{\lambda\mu}(0,t) S_\lambda
\qquad\hbox{and}\qquad
q_\nu(x;t) = \sum_\lambda L_{\nu\lambda}(t) s_\lambda.$$
By \cite[(4,8) and (4.10)]{macdonald},
$\langle q_\nu, m_\mu\rangle_{0,t} = \delta_{\nu\mu}$ and
$\langle s_\lambda, S_\mu\rangle = \delta_{\lambda\mu}$
in the inner product $\langle , \rangle_{0,t}$ of \cite[Ch.\ III]{macdonald}
so that 
$$L_{\nu\lambda}(t) = \langle q_\nu, S_\lambda\rangle
\qquad\hbox{which gives}\qquad
S_\lambda = \sum_{\nu} L_{\nu\lambda}(t) m_\nu.$$
Thus
$$
Q_\mu
= \sum_\nu
\Big(\sum_\lambda K_{\lambda\mu}(0,t) L_{\nu\lambda}(t) \Big) m_\nu.
$$
We have 
\begin{align*} s_{(1^n)} &= m_{(1^n)} = P_{(1^n)} \qquad\hbox{and} \\
s_{(21^{n-2})} &= P_{(21^{n-2})} + (t+t^2+\cdots+ t^{n-1})P_{(1^n)}
= m_{(21^{n-2})}+(n-1)m_{(1^n)}
\end{align*}
giving
$$P_{(1^n)}=m_{(1^n)} \qquad\hbox{and}\qquad
P_{(21^{n-2})} = m_{(21^{n-2})} + ((n-1) - (t+t^2+\cdots +t^{n-1}))m_{(1^n)}.$$
Thus
\begin{align*} Q_{(1^n)} &= b_{(1^n)}(t) m_{(1^n)}
\qquad\hbox{and}  \\
Q_{(21^{n-2})} &= b_{(21^{n-2})}(t) \big( m_{(21^{n-2}])} 
+ ((n-1) - (t+t^2+\cdots +t^{n-1}))m_{(1^n)}\big),
\end{align*}
and so $a_{(21^{n-2})\nu}(t) = b_{(21^{n-2})}(t) \tilde a_{(21^{n-2})\nu}$, where
\begin{equation}
\tilde a_{(21^{n-2})\nu}(t) 
= \sum_\lambda K_{\lambda\mu}(0,t) L_{\nu\lambda}(t) 
= \begin{cases}
1, &\hbox{if $\nu = (21^{n-2})$,} \\
(n-1) - (t+t^2+\cdots +t^{n-1}), &\hbox{if $\nu = (1^n)$,} \\
0, &\hbox{otherwise.}
\end{cases}
\label{afor21n}
\end{equation}

\subsection{A Character Computation} \label{sec: characterComputation}

Let $\FF_q$ be the finite field with $q$ elements,
$G = GL_n(q)$ the group of $n\times n$ invertible matrices with entries in $\FF$
and let $B$ be the subgroup of upper triangular matrices.
Then 
\begin{align*}
\Card(B) &= (q-1)^n q^{\frac12n(n-1)} \quad\hbox{and} \\
\Card(G) &= (q^n-1)(q^n-q)\cdots (q^n-q^{n-1}) = q^{\frac12 n (n-1)} (q-1)(q^2-1)\cdots (q^n-1)
\end{align*}
By \cite[Ch.\ II (1.6) and Ch.\ IV (2.7)]{macdonald},
\begin{align*}
\Card(Z_G(u_{(21^{n-2})})) 
&= q^{\vert (21^{n-2}) \vert + 2n(21^{n-2})} (1-q^{-1})\cdot (1-q^{-1})(1-q^{-2})\cdots (1-q^{-(n-2)}) \\
&= q^{n+(n-1)(n-2)}q^{-1-\frac12(n-1)(n-2)}(q-1)^2(q^2-1)(q^3-1)\cdots (q^{n-2}-1) \\
&= q^{n-1+\frac12(n-1)(n-2)}(q-1)^2(q^2-1)(q^3-1)\cdots (q^{n-2}-1) \\
&= q^{\frac12n(n-1)}(q-1)^2(q^2-1)(q^3-1)\cdots (q^{n-2}-1),
\end{align*}
 so that
$$\Card(\T) =\frac{\Card(G)}{\vert Z_G(u_{(21^{n-2})})\vert}
= \frac{(q^{n-1}-1)(q^n-1)}{(q-1)}
= \frac{(q^{n-1}-1)(q^n-1)}{(q-1)}.
$$
and
$$\frac{\vert B\vert}{\vert Z_G(u_{(21^{n-1})}) \vert}
= \frac{(q-1)}{[n-2]!},
\qquad\hbox{where}\quad
[n-2]! = \frac{(1-q)\cdots (1-q^{n-2})}{(1-q)^{n-2}}.
$$

Let $\mathbf{1}_B^G = \Ind_B^G(triv)$ be the representation of $G$ obtained by inducing
the trivial representation from $B$ to $G$.  As above, the Hecke algebra is $H = \End_G(\mathbf{1}_B^G)$.
The algebra $H$ has basis $\{ T_w\ |\ w\in S_n\}$ with multiplication determined by 
$$T_{s_k} T_w = \begin{cases}
(q-1)T_w + qT_{s_kw}, &\hbox{if $I(s_kw)<I(w)$,} \\
T_{s_kw}, &\hbox{if $I(s_kw)>I(w)$.}
\end{cases}
$$
As a $(G,H)$-bimodule,
$$\mathbf{1}_B^G \cong \bigoplus_\lambda G^\lambda \otimes H^\lambda,
$$
where the sum is over partitions of $n$, $G^\lambda$ is the irreducible $G$-module indexed
by $\lambda$ and $H^\lambda$ is the irreducible $H$-module indexed by $\lambda$.

Let $\T$ be the conjugacy class of transvections in $G=GL_n(q)$.  
A favorite representative of $\T$ is $u_{(21^{n-2})} = I + E_{12}$, where $I$ denotes the
identity matrix and $E_{12}$ is the matrix with $1$ in the $(1,2)$ entry and $0$ elsewhere.  
Let
$$D = \sum_{g\in \T} g,\qquad\quad
\hbox{an element of $Z(\CC[G])$,}$$
where $Z(\CC[G])$ denotes the center of the group algebra of $G$.
The element $D$ acts on $\mathbf{1}_\cB^G$ by
$$\sum_\lambda 
\frac{\chi^\lambda_G(\T)}{\chi^\lambda_G(1)}\cdot \id_{G^\lambda\otimes H^\lambda}\ ,$$
where $\id_{G^\lambda\otimes H^\lambda}$ is the operator that acts by the identity on
the component $G^\lambda \otimes H^\lambda$ and by $0$ on 
$G^\mu \otimes H^\mu$ for $\mu \ne \lambda$.

Let $z_\lambda$ be the element of $H$ that acts on $H^\lambda$ by the identity
and by $0$ on $H^\mu$ for $\mu \ne \lambda$.  With an abuse of notation, let $D\in H$ be given by
$$D = \sum_\lambda  \frac{\chi^\lambda_G(\T)}{\chi^\lambda_G(1)}z_\lambda,
\qquad\hbox{so that $D$ acts on $\mathbf{1}_B^G$ the same way that $C$ does.}
$$
(Note previously $D$ referred to element in the group algebra.)

Use cycle notation for permutations in $S_n$ so that 
$(i,j)$ denotes the transposition in $S_n$ that switches $i$ and $j$.
$$D_{(21^{n-2})} = \sum_{i<j} q^{(n-1)-(j-i)} T_{(i,j)}.$$
Note that at $q=1$ this specializes to the sum of the transpositions in the group algebra of the 
symmetric group.
Use cycle notation for permutations.   Then 
$$T_{s_k} (T_{(k,j)} + qT_{(k+1,j)}) 
= qT_{(k,k+1,j)} + q(q-1)T_{(k,j)}+qT_{(k+1,k,j)} = (T_{(k,j)} + qT_{(k+1,j)}) T_{s_k},
\qquad\hbox{if $k+1<j$,}
$$ 
$$T_{s_k} (qT_{(i,k)} + T_{(i,k+1)}) 
= qT_{(i, k,k+1)} + q(q-1)T_{(i,k+1)}+qT_{(i,k+1,k)} = (qT_{(i,k)} + T_{(i,k+1)}) T_{s_k},
\qquad\hbox{if $ i < k$,}
$$ 
and $T_{s_k}T_{(k,k+1)} = T_{(k,k+1)}T_{s_k}$ so that
$$\hbox{
if $k\in \{1, \ldots, n-1\}$ then $T_{s_k}D_{(21^{n-2})} = D_{(21^{n-2})}T_{s_k}$.}$$
So $D_{(21^{n-2})}\in Z(H)$.

These arguments now give a different proof of Theorem \ref{thm: secondD}, restated below. 
\begin{theorem} \label{Dvalue} Assume $n\ge 2$. Then
$$D = (n-1)q^{n-1} - [n-1] + (q-1) D_{(21^{n-2})}.
$$
\end{theorem}
\begin{proof}
In terms of the 
favorite basis $\{T_w\ |\ w\in S_n\}$ of $H$, the element $z_\lambda$ is given by
(see \cite[(68.29)]{curtisReinerBook} or \cite[(1.6)]{HLR})
$$z_\lambda = 
\frac{1}{\vert G/\cB\vert}
\sum_{w\in S_n} \chi^\lambda_G(1) q^{-I(w)} \chi^\lambda_H(T_{w^{-1}})T_w.$$
Thus, the expansion of $D$ in terms of the basis $\{T_w\ |\ w\in S_n\}$ of $H$ is 
\begin{align}
D &= \frac{1}{\vert G/\cB \vert} \sum_\lambda \sum_{w\in S_n} 
\frac{\chi^\lambda_G(\T)}{\chi^\lambda_G(1)}
\chi^\lambda_G(1) q^{-I(w)} \chi^\lambda_H(T_{w^{-1}})T_w  
\nonumber \\
&= \frac{\vert \cB\vert}{\vert G\vert} \sum_{w\in S_n} 
\Big( \sum_\lambda \chi^\lambda_G(\T) \chi^\lambda_H(T_{w^{-1}})\Big) 
q^{-I(w)}T_w 
\nonumber \\
&= \frac{\vert \cB\vert}{\vert G\vert} \Card(\T) \sum_{w\in S_n} 
\Big( \sum_\lambda \chi^\lambda_G(u_{21^{n-2}}) \chi^\lambda_H(T_{w^{-1}})\Big) 
 q^{-I(w)}T_w 
\nonumber \\
&= \frac{\vert \cB\vert}{\vert Z_G(u_{(21^{n-1})}) \vert} \sum_{w\in S_n} 
\Big( \sum_\lambda \chi^\lambda_G(u_{21^{n-2}}) \chi^\lambda_H(T_{w^{-1}})\Big) 
 q^{-I(w)}T_w.
 \label{Dexp}
\end{align}
For a partition $\nu$ let $\gamma_\nu$ be the favorite permutation of cycle type $\nu$
(a minimal length permutation which has cycle type $\nu$). 
By \cite[Th.\ 4.14]{Ra91},
$$\frac{q^{n} }{(q-1)^{I(\nu)} } q_\nu(0;q^{-1}) 
= \sum_\lambda \chi^\lambda_H(T_{\gamma_\nu^{-1}}) s_\lambda
\qquad\hbox{so that}\qquad
\chi^\lambda_H(T_{\gamma_\nu^{-1}}) = 
\frac{q^{n} }{(q-1)^{I(\nu)} }  L_{\nu\lambda}(q^{-1}),
$$
since $\nu$ is a partition of $n$.
By \cite[Theorem 4.9(c)]{HR99},
$$
\chi^\lambda_G(u_\mu) = q^{n(\mu)} K_{\lambda\mu}(0,q^{-1})
$$
Thus
\begin{align*}
\sum_\lambda &\chi^\lambda_G(u_{21^{n-2}}) \chi^\lambda_H(T_{\gamma_\nu^{-1}}) 
= \sum_\lambda q^{n(21^{n-2})} K_{\lambda(21^{n-2})}(0,q^{-1})
\frac{q^n}{(q-1)^{I(\nu)}}   L_{\nu\lambda}(q^{-1}) \\
&= \sum_\lambda 
\frac{q^{\frac12(n-1)(n-2) + n} }{ (q-1)^{I(\nu)} } 
K_{\lambda(21^{n-2})}(0,q^{-1})
L_{\nu\lambda}(q^{-1}) 
= \frac{q^{\frac12(n-1)(n-2) + n} }{ (q-1)^{I(\nu)} } a_{(21^{n-2})\nu}(q^{-1}).
\end{align*}
Plugging this into \eqref{Dexp},
the coefficient of $T_{\gamma_\nu}$ in $D$ is
\begin{align*}
\frac{\vert \cB\vert}{\vert Z_G(u_{(21^{n-1})})\vert}  
\frac{q^{\frac12(n-1)(n-2) + n} }{ (q-1)^{I(\nu)} } 
&a_{(21^{n-2})\nu}(q^{-1}) q^{-I(\gamma_\nu)}
=\frac{(q-1)}{[n-2]!} \frac{q^{\frac12(n-1)(n-2)+n} }{ (q-1)^{I(\nu)} } 
a_{(21^{n-2})\nu}(q^{-1}) q^{-(\vert \nu\vert - I(\nu))}
\\
&=\frac{q^{\frac12(n-1)(n-2)} (q-1)}{[n-2]!} \frac{q^{I(\nu)} }{ (q-1)^{I(\nu)} } 
a_{(21^{n-2})\nu}(q^{-1})
\end{align*}
Now use \eqref{afor21n},
$$b_{(21^{n-2})}(q^{-1}) = (1-q^{-1})(1-q^{-1})\cdots (1-q^{-(n-2)})
= q^{-1-\frac12(n-1)(n-2)}(q-1)^{n-1} [n-2]!,
$$
and
$$
\tilde a_{(21^{n-2})\nu}(q^{-1}) 
= \begin{cases}
1, &\hbox{if $\nu = (21^{n-2})$,} \\
\displaystyle{ (n-1) - q^{1-n}\frac{(q^{n-1}-1)}{q-1},} &\hbox{if $\nu = (1^n)$}, \\
0, &\hbox{otherwise},
\end{cases}
$$
to get that the coefficient of $T_{\gamma_\nu}$ in $D$ is
\begin{align*}
\frac{(q-1)}{[n-2]!} 
&q^{\frac12(n-1)(n-2)} \frac{q^{I(\nu)} }{ (q-1)^{I(\nu)} } 
q^{-1-\frac12(n-1)(n-2)}(q-1)^{n-1} [n-2]!
\cdot \tilde a_{(21^{n-2})\nu}(q^{-1})
\\
&= (q-1)^{n-I(\nu)} q^{I(\nu)-1} 
\tilde a_{(21^{n-2})\nu}(q^{-1})
= \begin{cases}
(q-1)q^{n-2}, &\hbox{if $\nu = (21^{n-2})$,} \\
\displaystyle{ 
(n-1)q^{n-1} - [n-1]_q, } &\hbox{if $\nu = (1^n)$}, \\
0, &\hbox{otherwise},
\end{cases}
\end{align*}
Since $\tilde a_{(21^{n-2})\nu}(q^{-1})$ is 0 unless $\nu=(1^n)$ or $\nu = (21^n)$
and since the coefficient of $T_{\gamma_{(21^{n-2})}} = T_{s_1}$ in $D_{(21^{n-2})}$ is
$q^{n-2}$ then it follows from \cite[(2.2) Main Theorem]{Fr99} that 
$$D = (n-1)q^{n-1} - [n-1]_q + (q-1) D_{(21^{n-2})}.
$$
 \end{proof}

\subsection{The Example for $n=2$} \label{sec: n=2example}

If $n=2$ then $H = \hbox{span}\{1, T_{s_1}\}$ with 
$(1+T_{s_1})(q-T_{s_1}) = 0$ and 
$$z_{(2)} = \frac{1}{1+q}(1+T_{s_1})
\qquad\hbox{and}\qquad
z_{(1^2)} = \frac{q}{1+q}(1-q^{-1}T_{s_1}).
$$
In this case $Z(H) = H$ and
$D_{(2)} = T_{(1,2)} = T_{s_1} = qz_{(2)} - z_{(1^2)}.$

So $D_{(2)} = T_{s_1}$ and $H$ has basis $\{1, T_{s_1}\}$ and
\begin{align*}
\begin{array}{rl}
1\cdot D_{(2)} &= 1\cdot T_{s_1} = T_{s_1}, \\
T_{s_1}\cdot D_{(2)} &= T_{s_1}\cdot T_{s_1} = (q-1)T_{s_1}+q,
\end{array}
\qquad\hbox{so that}\qquad
\begin{pmatrix} 0 &q \\ 1 &q-1 \end{pmatrix}
\end{align*}
is the matrix for multiplying by $D_{(2)}$.
The character values are
\begin{align*}
\chi^{(1^2)}_G(u_{(1^2)}) = q^{n(1^2)}K_{(1^2)(1^2)}(0,q^{-1}) = q^1\cdot 1 = q, \\
\chi^{(1^2)}_G(u_{(2)}) = q^{n(2)}K_{(1^2)(2)}(0,q^{-1}) = q^0\cdot 0 = 0, \\
\chi^{(2)}_G(u_{(1^2)}) = q^{n(1^2)}K_{(2)(1^2)}(0,q^{-1}) = q^1\cdot q^{-1} = 1, \\
\chi^{(2)}_G(u_{(2)}) = q^{n(2)}K_{(2)(2)}(0,q^{-1}) = q^0\cdot 1 = 1,
\end{align*}
giving that $C$ acts on $\mathbf{1}_B^G$ the same way as
\begin{align*}
D&=\Card(\T)\Big(
\frac{\chi^{(1^2)}(u_{(2)})}{\chi^{(1^2)}(u_{(1^2)})} z_{(1^2)}
+\frac{\chi^{(2)}(u_{(2)})}{\chi^{(2)}(u_{(1^2)})} z_{(2)}\Big)
= \frac{\Card(G)}{\Card(Z_G(u_{(2)}))} \Big( \frac{0}{q} z_{(1^2)} + \frac{1}{1}z_{(2)}\Big) \\
&=\frac{(q^2-1)(q^2-q)}{ q^{2+0}(1-q^{-1})} z_{(2)} 
=\frac{(q^2-1)(q-1)}{ q-1} z_{(2)} 
= (q^2-1)\frac{1}{1+q}(1+T_{s_1}) \\
&= (q-1)(1+D_{(2)}).
\end{align*}
and the matrix for multiplying by $D$ in the basis $\{1, T_{s_1}\}$ is
$$(q-1) \left( 
\begin{pmatrix} 1 &0 \\ 0 &1 \end{pmatrix}
+\begin{pmatrix} 0 &q \\ 1 &q-1 \end{pmatrix}
\right)
=(q-1) \begin{pmatrix} 1 &q \\ 1 &q \end{pmatrix}.$$
In this case $\vert \cB \vert = (q-1)^2q$ and $\vert Z_G(u_2)\vert = q(q-1)$
so that 
$$\frac{ \vert B\vert }{\vert Z_G(u_2)\vert } = (q-1).$$
Then
$$\begin{array}{ll}
\chi^{(1^2)}(T_{\gamma_{(1^2)}}) =\chi^{(1^2)}(1) = 1,\qquad
&\chi^{(1^2)}(T_{\gamma_{(2)}}) = \chi^{(1^2)}(T_{s_1}) = -1, \\
\chi^{(2)}(T_{\gamma_{(1^2)}}) = \chi^{(2)}(1) = 1,\qquad
&\chi^{(2)}(T_{\gamma_{(2)}}) = \chi^{(2)}(T_{s_1}) = q,
\end{array}
$$
so that
\begin{align*}
\sum_{\lambda} \chi^\lambda_G(u_2)\chi^\lambda_H(1)
&= \chi^{(2)}_G(u_2) \chi^{(2)}_H(T_{\gamma_{(1^2)}})
+ \chi^{(1^2)}_G(u_2) \chi^{(1^2)}_H(T_{\gamma_{(1^2)}}) 
= 1\cdot 1
+ 0\cdot 1 = 1, \\
\sum_{\lambda} \chi^\lambda_G(u_2)\chi^\lambda_H(T_{s_1})
&= \chi^{(2)}_G(u_2) \chi^{(2)}_H(T_{\gamma_{(2)}})
+ \chi^{(1^2)}_G(u_2) \chi^{(1^2)}_H(T_{\gamma_{(2)}}) 
= 1\cdot q + 0\cdot (-1) = q, 
\end{align*}
so that the coefficient of $1$ in $D$ is $(q-1)\cdot 1\cdot q^{I(1)} = (q-1)\cdot 1\cdot 1 = (q-1)$ 
and the coefficient of $T_{s_1}$ is  $(q-1)\cdot q\cdot q^{-I(s_1)} = (q-1)$.  Thus
$$ D = (q-1)(1+T_{s_1}).$$
For $n=2$, Theorem \ref{Dvalue} says that 
\begin{align*}
D &= (q-1) \big( (q-1)^{-2}q^2a_{(2)(1^2)}(q^{-1})+(q-1)^{-1}q D_{(2)}\big)  
= q^{-2} \big( (q-1)^2 (1-q^{-1}) + (q-1)D_{(2)} \big) 
\end{align*}

\subsection{The Example of $n=3$} \label{sec: n=3example}

If $n=3$ then
$H = \mathrm{span}\{1, T_{s_1}, T_{s_2}, T_{s_1s_2}, T_{s_2s_1}, T_{s_1s_2s_1}\}$ with 
$$(1+T_{s_1})(q-T_{s_1}) = 0,
\quad
(1+T_{s_2})(q-T_{s_2}) = 0
\quad\hbox{and}\quad
T_{s_1}T_{s_2}T_{s_1} = T_{s_2}T_{s_1}T_{s_2}.$$
Then 
\begin{align*}
z_{(3)} &= \frac{1}{[3]!}(1+T_{s_1}+T_{s_2}+T_{s_2s_1}+T_{s_1s_2}+T_{s_1s_2s_1}), \\
z_{(1^3)} &= \frac{q^3}{[3]!}(1-q^{-1}T_{s_1}-q^{-1}T_{s_2}+q^{-2}T_{s_2s_1}+q^{-2}T_{s_1s_2}-q^{-3}T_{s_1s_2s_1}) \\
z_{(21)} &= 1 - z_{(3)} - z_{(1^3)}.
\end{align*}
Then
\begin{align*}
(q+1)D_{(21)} &= (q+q^2)T_{s_1}+(q+q^2)T_{s_2}+(q+1)T_{s_1s_2s_1}
= [3]!(q z_{(3)} - z_{(1^2)}) - (q-q^3)\cdot 1 \\
&= \frac{[3]!}{ [3]!} ((q-q^3) + (q+q^2)T_{s_1}+(q+q^2)T_{s_2}+0+0+(q+1)T_{s_1s_2s_1}) - (q-q^3)
\\
&= (q+1)(qT_{s_1} + qT_{s_2} + T_{s_1s_2s_1})
\end{align*}
and the brute force computation showing that $D_{(21)} = qT_{s_1}+qT_{s_2}+T_{s_1s_2s_1}$
commutes with $T_{s_1}$ is
\begin{align*}
T_{s_1}D_{(21)} &= q(q-1)T_{s_1}+q^2 + qT_{s_1s_2} + (q-1)T_{s_1s_2s_1} + q T_{s_2s_1} \\
D_{(21)}T_{s_1} &= q(q-1)T_{s_1}+q^2 + qT_{s_2s_1} + (q-1)T_{s_1s_2s_1} + q T_{s_1s_2} .
\end{align*}
In this computation already, some flags in $\cB s_1\cB$ are ending up in $\cB s_1s_2s_1\cB$
after the application of $D_{(21)}$.

In the basis $\{1, T_{s_1}, T_{s_2}, T_{s_1s_2}, T_{s_2s_1}, T_{s_1s_2s_1}\}$, the matrices of
multiplication by $T_{s_1}$ and $T_{s_2}$ are
$$\begin{pmatrix}
0 &q &0 &0 &0 &0 \\
1 &q-1 &0 &0 &0 &0 \\
0 &0 &0 &q &0 &0 \\
0 &0 &1 &q-1 &0 &0 \\
0 &0 &0 &0 &0 &q \\
0 &0 &0 &0 &1 &q-1
\end{pmatrix}
\qquad\hbox{and}\qquad
\begin{pmatrix}
0 &0 &q &0 &0 &0 \\
0 &0 &0 &0 &q &0 \\
1 &0 &q-1 &0 &0 &0 \\
0 &0 &0 &0 &0 &q \\
0 &1 &0 &0 &q-1 &0 \\
0 &0 &0 &1 &0 &q-1
\end{pmatrix}
$$
The matrix of multiplication by $T_{s_1s_2s_1}$ is
$$\begin{pmatrix}
0 &0 &0 &0 &0 &q^3 \\
0 &0 &0 &0 &q^2 &q^2(q-1) \\
0 &0 &0 &q^2 &0 &q^2(q-1) \\
0 &0 &q &q(q-1) &q(q-1) &q(q-1)^2 \\
0 &q &0 &q(q-1) &q(q-1) &q(q-1)^2 \\
1 &q-1 &q-1 &(q-1)^2 &(q-1)^2 &(q-1)^3+q(q-1)
\end{pmatrix}
$$
where the last column comes from the computation
\begin{align*}
T_{s_1s_2s_1}T_{s_1}T_{s_2}T_{s_1}
&= \big((q-1)T_{s_1s_2s_1} + q T_{s_1s_2}\big) T_{s_2s_1} \\
&= \big( ((q-1)^2T_{s_1s_2s_1}+q(q-1)T_{s_2s_1})+(q(q-1)T_{s_1s_2}+q^2T_{s_1})\big)T_{s_1} \\
&= ((q-1)^3 T_{s_1s_2s_1} + q(q-1)^2 T_{s_1s_2}) 
+(q(q-1)^2T_{s_2s_1}+q^2(q-1)T_{s_2}) \\
&\qquad +q(q-1)T_{s_1s_2s_1}
+(q^2(q-1)T_{s_1}+q^3)
\end{align*}
Thus the matrix of multiplication by $D_{21} = qT_{s_1}+qT_{s_2}+T_{s_1s_2s_1}$ is
$$\begin{pmatrix}
0 &q^2 &q^2 &0 &0 &q^3 \\
q &q(q-1) &0 &0 &2q^2 &q^2(q-1) \\
q &0 &q(q-1) &2q^2 &0 &q^2(q-1) \\
0 &0 &2q &2q(q-1) &q(q-1) &q^2+q(q-1)^2 \\
0 &2q &0 &q(q-1) &2q(q-1) &q^2+q(q-1)^2 \\
1 &q-1 &q-1 &q+(q-1)^2 &q+(q-1)^2 &(q-1)^3+3q(q-1)
\end{pmatrix}
$$
with $1+2(q-1)+2q+2(q-1)^2=2q^2-4q+2+4q-2+1=2q^2+1$
and $(q-1)^3+3q(q-1)=q^3-3q^2+3q-1+3q^2-3q=q^3-1$ so that the bottom row sums to
$q^3-1+2q^2+1=q^3+2q^2$
Thus the matrix of $D = (n-1)q^{n-1}-[n-1]+(q-1)D_{(21)}
=2q^2-(q+1)+(q-1)D = (q-1)(2q+1)+(q-1)D$ which has row sums
$$
(q-1)\begin{pmatrix}
(2q+1)+(2q^2+q^3) \\
(2q+1)+(2q^2+q^3) \\
(2q+1)+(2q^2+q^3) \\
(2q+1)+(2q+3q(q-1)+q^2+q^3-2q^2+q) \\
(2q+1)+(2q+3q(q-1)+q^2+q^3-2q^2+q) \\
(2q+1)+(q^3+2q^2)) 
\end{pmatrix}
$$
and
$$\Card(\T) = \frac{(q^{n-1}-1)(q^n-1)}{q-1} = \frac{(q^2-1)(q^3-1)}{q-1}
=(q-1)(q+1)(q^2+q+1) = (q-1)(q^3+2q^2+2q+1).$$

Let
$$M_2 = T_{s_1}T_{s_1}
\qquad\hbox{and}\qquad
M_3 = T_{s_2}T_{s_1}T_{s_1}T_{s_2}
$$
so that
$$M_2= (q-1)T_{s_1}+q
\qquad\hbox{and}\qquad
M_3 = T_{s_2}((q-1)T_{s_1}+q)T_{s_2}
=(q-1)T_{s_2s_1s_2}+q(q-1)T_{s_2}+q^2,$$
giving
\begin{align*}
q(M_2-1)+(M_3-1) 
&= (q-1)(T_{s_1s_2s_1}+qT_{s_2}+qT_{s_1})+2q^2-q-1 \\
&=(3-1)q^{3-1}-[3-1]+(q-1)D_{(21)}=D.
\end{align*}

Let's work out the $q_\mu$ and the matrix $L$.  By definition
\begin{align*}
\sum_r q_r(x;t) z^r 
&= \prod_{i=1}^n \frac{1-tx_iz}{1-x_iz}
\Big( \sum_k (-1)^ke_kt^kz_k\Big) \big(\sum_{\ell} h_\ell z^\ell\big) \\
&=(1-zte_1+z^2t^2e_2-t^3z^3e_3+\cdots)(1+h_1z+h_2z^2+\cdots)
\end{align*}
giving $q_0 = 1$, $q_1 = h_1-te_1 = (1-t)s_{(1)}$, 
$$q_2 = t^2e_2-tzh_1e_1+h_2 = t^2(s_{(1^2)}-t(s_{(1^2)}+s_{(2)})+s_{(2)}
=(1-t)\big ( (-t)s_{(1^2)} + s_{(2)} \big)
$$
and
\begin{align*}
q_3 &= (-t)^3 e_3+(-t)^2 h_1e_2 + (-t)h_2e_1+h_3 \\
&=(-t)^3 s_{(1^3)} + (-t)^2(s_{(21)}+s^{(1^3)}) +(-t)(s_{(3)}+s_{(21)}+s_{(3)}) \\
&=(1-t)\big((-t)^2s_{(1^3)} +(-t) s_{(21)} + s_{(3)}\big).
\end{align*}
Thus
\begin{align*}
q_{(1^2)} &= q_1^2 = (1-t)^2(s_{(1^2)}+s_{(2)}), \qquad\qquad
q_{(2)}=(1-t)\big ( (-t)s_{(1^2)} + s_{(2)} \big) \\
q_{(1^3)} &= q_1^3 = (1-t)^3(s_{(1^3)}+2s_{(21)} + s_{(3)}\big), \\
q_{(21)} &= q_2q_1 = (1-t)^2((-t)(s_{(1^3)}+s_{(21)})+(s_{(3)}+s_{(21)}))
=(1-t)^2( (-t)s_{(1^3)} + (1-t)s_{(21)} + s_{(3)}), \\
q_3 &=(1-t)\big((-t)^2s_{(1^3)} +(-t) s_{(21)} + s_{(3)}\big).
\end{align*}
Then
\begin{align*}
\chi^{(1^2)}_H(T_{\gamma_{(1^2)}}) 
&= 1 = \frac{1}{(1-q^{-1})^2} (1-q^{-1})^2 = \frac{q^2}{(q-1)^2} L_{(1^2)(1^2)}(q^{-1}), \\
\chi^{(2)}_H(T_{\gamma_{(1^2)}}) 
&= 1 = \frac{1}{(1-q^{-1})^2} (1-q^{-1})^2 = \frac{q^2}{(q-1)^2} L_{(1^2)(2)}(q^{-1}), \\
\chi^{(1^2)}_H(T_{\gamma_{(2)}}) 
&= -1 = \frac{q}{(1-q^{-1})} (1-q^{-1})(-q^{-1}) = \frac{q^2}{(q-1)} L_{(2)(1^2)}(q^{-1}), \\
\chi^{(2)}_H(T_{\gamma_{(2)}}) 
&= q = \frac{q}{(1-q^{-1})} (1-q^{-1}) = \frac{q^2}{(q-1)}L_{(2)(2)}(q^{-1}), \\
\end{align*}
\begin{align*}
\chi^{(1^3)}_H(T_{\gamma_{(1^3)}}) 
&= 1 = \frac{1}{(1-q^{-1})^3} (1-q^{-1})^3 = \frac{1}{(1-q^{-1})^3} L_{(1^3)(3)}(q^{-1}), \\
\chi^{(21)}_H(T_{\gamma_{(1^3)}}) 
&= 2 = \frac{1}{(1-q^{-1})^3} (1-q^{-1})^3\cdot 2 = \frac{1}{(1-q^{-1})^3} L_{(1^3)(21)}(q^{-1}), \\
\chi^{(3)}_H(T_{\gamma_{(1^3)}}) 
&= 1 = \frac{1}{(1-q^{-1})^3} (1-q^{-1})^3 = \frac{1}{(1-q^{-1})^3} L_{(1^3)(3)}(q^{-1}), \\
\chi^{(1^3)}_H(T_{\gamma_{(21)}}) 
&= -1 = \frac{q}{(1-q^{-1})^2} (1-q^{-1})^2(-q^{-1}) = \frac{q^3}{(q-1)^2} L_{(21)(1^3)}(q^{-1}), \\
\chi^{(21)}_H(T_{\gamma_{(21)}}) 
&= 1-q = \frac{q}{(1-q^{-1})^2} (1-q^{-1})^2(1-q^{-1}) = \frac{q^3}{(q-1)^2}L_{(21)(21)}(q^{-1}), \\
\chi^{(3)}_H(T_{\gamma_{(21)}}) 
&= q = \frac{q}{(1-q^{-1})^2} (1-q^{-1})^2\cdot 1 = \frac{q^3}{(q-1)^2}L_{(21)(3)}(q^{-1}), \\
\chi^{(1^3)}_H(T_{\gamma_{(3)}}) 
&= 1 = \frac{q^2}{(1-q^{-1})} (1-q^{-1})(-q^{-1})^2 = \frac{q^3}{(q-1)} L_{(3)(1^3)}(q^{-1}), \\
\chi^{(21)}_H(T_{\gamma_{(3)}}) 
&= -q = \frac{q^2}{(1-q^{-1})} (1-q^{-1})(-q^{-1}) = \frac{q^3}{(q-1)}L_{(3)(21)}(q^{-1}), \\
\chi^{(3)}_H(T_{\gamma_{(3)}}) 
&= q^2 = \frac{q^2}{(1-q^{-1})} (1-q^{-1})\cdot 1 = \frac{q^3}{(q-1)}L_{(3)(21)}(q^{-1}), \\
\end{align*}

Since
\begin{align*}
\chi^{(1^3)}_G(u_{(21)}) &= q^{n(21)} K_{(1^3)(21)}(0,q^{-1}) = q \cdot 0=0, 
\\
\chi^{(21)}_G(u_{(21)}) &= q^{n(21)} K_{(21)(21)}(0,q^{-1}) = q \cdot 1 = q, 
\\
\chi^{(3)}_G(u_{(21)}) &= q^{n(21)} K_{(3)(21)}(0,q^{-1}) = q \cdot q^{-1} = 1, 
\end{align*}
and
\begin{align*}
\chi^{(1^3)}_G(u_{(1^3)}) &= q^{n(1^3)} K_{(1^3)(1^3)}(0,q^{-1}) = q^3 \cdot 1 = q^3, 
\\
\chi^{(21)}_G(u_{(1^3)}) &= q^{n(1^3)} K_{(21)(1^3)}(0,q^{-1}) = q^3 \cdot (q^{-1}+q^{-2}) = q^2+q, 
\\
\chi^{(3)}_G(u_{(1^3)}) &= q^{n(1^3)} K_{(3)(1^3)}(0,q^{-1}) = q^3 \cdot q^{-3} = 1, 
\end{align*}
then $C$ acts the same way as 
\begin{align*}
D 
&= \frac{(q^{n-1}-1)(q^n-1)}{(q-1)}
\big(\frac{0}{q^3}z_{(1^3)}+\frac{q}{q^2+q}z_{(21)} + z_{(3)}\big) 
= \frac{(q^2-1)(q^3-1)}{(q-1)}
\big(\frac{1}{q+1}z_{(21)}+ z_{(3)}\big) \\
&= \frac{(q^2-1)(q^3-1)(q-1)}{(q-1)(q^2-1)}
\big(1 - z_{(3)} - z_{(1^3)}+(q+1)z_{(3)}\big) 
= (q^3-1) \big(1 - z_{(1^3)}+qz_{(3)}\big) \\
&= \frac{(q^3-1)}{[3]!} \left(
\begin{array}{l}
[3]! - (q^3-q^2T_{s_1}-q^2T_{s_2}+qT_{s_1s_2}+qT_{s_2s_1}-T_{s_1s_2s_1}) \\
+q (1+T_{s_1}+T_{s_2}+T_{s_1s_2}+T_{s_2s_1}+T_{s_1s_2s_1})  
\end{array}
\right) \\
&= \frac{(q^3-1)}{[3]!} \big( ([3]!-q^3+q)+(q^2+q)T_{s_1}+ (q^2+q)T_{s_2}+(1+q)T_{s_1s_2s_1} \big) \\
&= \frac{(q^3-1)}{[3]!}([3]!-q^3+q)
+\frac{(q^3-1)}{[3]!}(1+q)D_{(21)}
= (q^3-1) + \frac{(q-1)}{(1+q)}(q-q^3)
+\frac{(q^3-1)}{[3]!}\frac{(q^2-1)}{(q-1)}D_{(21)} \\
&= (q^3-1)-q(q-1)^2 +(q-1) D_{(21)} 
= 2q^2-q-1 + (q-1) D_{(21)}  \\
&= (3-1)q^{3-1}-(q+1) + (q-1) D_{(21)}.
\end{align*}

Then
\begin{align*}
\sum_{\lambda} &\chi^\lambda(u_{(21)}) \chi^\lambda_H(T_{\gamma_{(1^3)}})
= \chi^{(1^3)}(u_{(21)}) \chi^{(1^3)}_H(T_{\gamma_{(1^3)}})
+\chi^{(21)}(u_{(21)}) \chi^{(21)}_H(T_{\gamma_{(1^3)}})
+\chi^{(3)}(u_{(21)}) \chi^{(3)}_H(T_{\gamma_{(1^3)}}) \\
&=0\cdot 1+ q\cdot 2+ 1\cdot 1 = 2q+1, \\
\sum_{\lambda} &\chi^\lambda(u_{(21)}) \chi^\lambda_H(T_{\gamma_{(21)}})
= \chi^{(1^3)}(u_{(21)}) \chi^{(1^3)}_H(T_{\gamma_{(21)}})
+\chi^{(21)}(u_{(21)}) \chi^{(21)}_H(T_{\gamma_{(21)}})
+\chi^{(3)}(u_{(21)}) \chi^{(3)}_H(T_{\gamma_{(21)}}) \\
&=0\cdot (-1) + q\cdot (1-q) + 1\cdot q = 2q-q^2, \\
\sum_{\lambda} &\chi^\lambda(u_{(21)}) \chi^\lambda_H(T_{\gamma_{(3)}})
= \chi^{(1^3)}(u_{(21)}) \chi^{(1^3)}_H(T_{\gamma_{(3)}})
+\chi^{(21)}(u_{(21)}) \chi^{(21)}_H(T_{\gamma_{(3)}})
+\chi^{(3)}(u_{(21)}) \chi^{(3)}_H(T_{\gamma_{(3)}}) \\
&=0\cdot 1 + q\cdot (-q) + 1\cdot q^2 = 0,
\end{align*}

\section{A Single Coset Lumping} \label{sec: topCard}

This section develops the correspondence between $\cB \backslash GL_n(q) / \cB$ double-cosets and flags, and describes the random transpositions Markov chain in this setting. This description is useful for analyzing specific features of the Markov chain. If $\{\omega_t \}_{t \ge 0}$ is the induced chain on $S_n$, then $\{\omega_t(1) \}_{t \ge 0}$ is a process on $\{1, \dots, n \}$. Thinking of $\omega \in S_n$ as a deck of cards, this is the evolution of the label of the `top card'. It is the lumping of the chain on $S_n$ onto cosets  $S_n/S_{n-1}$. The main result, Lemma \ref{lem: SSTR} below, shows that the top card takes only a bounded number of steps to reach stationarity.

\subsection{Flag Representation} \label{sec: flag}

The subgroup $\cB$ gives rise to the quotient $GL_n(q)/\cB$. This may be pictured as the space of `flags'. 

\begin{definition}
Here a flag $F$ consists of an increasing sequence of subspaces $F = F_1 \subset F_2 \subset \hdots \subset F_n$ with $\text{dim}(F_i) = i$. The \emph{standard flag} is 
\[
E = \left\langle \me_1 \right\rangle \subset \left\langle \me_1, \me_2 \right\rangle \hdots \subset \left\langle \me_1, \dots, \me_n \right\rangle.
\]
\end{definition}

Indeed, $GL_n(q)$ operates transitively on flags and the subgroup fixing the \emph{standard flag} is exactly $\cB$. There is a useful notion of `distance' between two flags $F, F'$ which defines a permutation.

\begin{definition}
Let $F, F'$ be two flags. The \emph{Jordan-Holder permutation} $\omega(F, F')$ is a permutation $\omega = \omega(F, F') \in S_n$ defined by $\omega(i) = j$ if $j$ is the smallest index such that 
\[
F_{i - 1} + F_j' \supseteq F_i.
\]
\end{definition}
The Jordan-Holder permutation distance satisfies
\begin{equation}
    \omega(F, F') = \omega(MF, MF'), \, M \in GL_n(q) \quad \omega(F, F') = \omega(F', F)^{-1}.
\end{equation}
A thorough development with full proofs is in \cite{abels}. 

\begin{lemma}
For $M \in GL_n(q)$ and $E$ the standard flag,
\[
M \in \cB \omega \cB \iff \omega = \omega(E, ME).
\]
\end{lemma}


\subsection{Top Label Chain}

This representation is useful for analyzing a further projection of the chain on $S_n$: Let $\{\omega_t \}_{t \ge 0}$ denote the Markov chain on $S_n$. Let $P_1(\cdot, \cdot)$ denote the marginal transition probabilities of the first position. That is,
\[
P_1(j, k) = \P(\omega_{t+1}(1) = k \mid \omega_t(1) = j ), \quad j, k \in \{1, \dots, n \}.
\]
If $\omega \sim \pi_q$, then the marginal distribution on $\{1, \dots, n \}$ of the first card $\omega(1)$:
\begin{align*}
\pi_{q, 1}(j) := P(\omega(1) = j) = \sum_{\omega: \omega(1) = j} \frac{q^{I(\omega)}}{[n]_q!} = \frac{q^{j-1}}{1 + q + \hdots + q^{n-1}}.
\end{align*}

\begin{lemma} \label{lem: P1}
For $1 \le i, j \le n$,
\begin{align*}
P_1(i, j) &= \begin{cases}
\frac{(q-1)^2q^{n + j - 3}}{(q^n - 1)(q^{n-1} - 1)} & i \neq j \\
\frac{(q^{n-1} - 1)^2 + (q - 1)^2(q^{i - 1} - 1)q^{n - 2}}{(q^n - 1)(q^{n-1} - 1)} & i = j
\end{cases}.
\end{align*}
\end{lemma}

\begin{remark} \label{remark: pq1}
Observe that for $i \neq j$ the transition probabilities can be written 
\begin{align*}
&P_1(i, j) = \frac{q^{j-1}}{(q^n - 1)/(q-1)} \left( \frac{(q - 1) q^{n-2}}{q^{n-1} - 1} \right) = \pi_{q, 1}(j) \cdot \left( \frac{(q - 1) q^{n-2}}{q^{n-1} - 1} \right) \\
&P_1(j, j) = \pi_{q, 1}(j) \cdot \left( \frac{(q - 1) q^{n-2}}{q^{n-1} - 1} \right) + \left(1 - \frac{(q - 1) q^{n-2}}{q^{n-1} - 1} \right).
\end{align*}
Write $p:= ((q - 1) q^{n-2})/(q^{n-1} - 1)$. This provides another description of the Markov chain: At each step, flip a coin which gives heads with probability $p$, tails with probability $1-p$. If heads, move to a random sample from $\pi_{q, 1}$. Otherwise, don't move. 
\end{remark}

\begin{remark}
Though the Markov chain $P_1$ on $\{1, \dots, n \}$ was defined via lumping from a chain on the group $GL_n(q)$ with $q > 1$ a prime power, note that the transitions are well-defined even for $q < 1$. If $q < 1$ then the Mallows measure $\pi_q$ is concentrated at the identity permutation and $P_1(i, j)$ is largest for $j = 1$. Note also that the description in Remark \ref{remark: pq1} is also valid, since $p = ((q - 1) q^{n-2})/(q^{n-1} - 1) < 1$ for an $q > 0$.
\end{remark}

\begin{proof}
If $\omega(1) = i$, then a flag representing the double coset has $F_1 = \mathbf{e}_i$. Recall that a transvection $T$ defined by vectors $\mv, \ma$ has $T(\mathbf{e}_i) = \mathbf{e}_i + a_i \mv$. The first coordinate in the new permutation is the smallest $j$ such that
\[
\me_i + a_i \mv \subset \left\langle e_1, \dots, e_j \right\rangle.
\]
There are two cases when this smallest $j$ is equal to $i$:
\begin{enumerate}
    \item $a_i = 0$, and $\mv$ can be anything: There are $(q^{n - 1} - 1)$ possibilities for $\ma$ such that $a_i = 0$. Then there are $(q^{n-1} - 1)/(q - 1)$ possibilities for $\mv$. 
    
    \item $a_i \neq 0$ and $\mv$ is such that $a_i v_i \neq - 1$ and $v_k = 0$ for all $k > i$. Note that if $i = 1$, then this is not possible. Since the nonzero entry at the largest index in $\mv$ is normalized to be $1$, there are two further possibilities:
    
    \begin{itemize}
        \item If $v_i = 1$, then there are $(q - 2)$ possibilities for $a_i$,  $(q^{i - 1} -1)$ possibilities for the rest of $\mv$ (note that it's not allowed for the rest of $\mv$ to be $0$, because then we could not get $\ma^\top \mv = 0$), and then $q^{n-2}$ possibilities for the rest of $\ma$. 
        
        \item If $v_i = 0$, then there are $(q^{i - 1} - 1)/(q - 1)$ possibilities for $\mv$ and then $(q - 1)q^{n-2}$ possibilities for $\ma$.
    \end{itemize}
    
\end{enumerate}
In summary, if $i \neq 1$, then
\begin{align*}
    P_1(i, i) &= \frac{1}{|\T_{n, q}|} \left( \frac{(q^{n-1} - 1)(q^{n - 1} -1)}{q - 1} + (q^{i - 1} -1)(q - 2)q^{n-2} + \frac{(q^{i - 1} -1)(q - 1)q^{n-2}}{q - 1} \right) \\
    &= \frac{(q^{n-1} - 1)^2 + (q - 1)^2(q^{i - 1} - 1)q^{n - 2}}{(q^n - 1)(q^{n-1} - 1)}.
\end{align*}


Now if $i \neq j$ we consider the two possibilities.
\begin{enumerate}
    \item $i > j$: This transition will occur if $a_i\cdot v_i = -1$ and $\mv$ is such that $v_j > 0$ and $v_k = 0$ for all $j < k < i$. Since the $\mv$ is normalized so that the entry at the largest index is. $1$, this requires $v_i = 1$ and $a_i = -1$. There are then $(q - 1)q^{j - 1}$ such $\mv$, and for each $\mv$ there are $q^{n-2}$ possibilities for $\ma$. This gives in total
    \[
    (q - 1) q^{n + j - 3}.
    \]
    
    \item $i < j$: This transition will occur if $a_i \neq 0$, $v_j = 1$ and $v_k = 0$ for $k > j$. There are $q^{j - 1}$ possibilities for $\mv$ and then $(q - 1)q^{n-2}$ possibilities for $\ma$, so again in total
    \[
    (q - 1)q^{n+ j - 3}.
    \]
\end{enumerate}
In summary, for any $i \neq j$,
\begin{align*}
    P_1(i, j) = \frac{(q-1)q^{n + j - 3}}{\vert \T_{n, q} \vert} = \frac{(q-1)^2q^{n + j - 3}}{(q^n - 1)(q^{n-1} - 1)}. 
\end{align*}

\end{proof}



\begin{lemma}
Let $P_1$ be the Markov chain from Lemma \ref{lem: P1} with stationary distribution $\pi_{q, 1}$. Then $P_1$ has eigenvalues $\{1, \beta \}$ with
\[
\beta = \frac{q^{n -2} -1}{q^{n-1} - 1}.
\]
The eigenvalue $\beta$ has multiplicity $n -1$. If $q > 1$ is fixed and $n$ large, then $\beta \sim 1/q$.
\end{lemma}

\begin{proof}
Call $P_1(i, j) = c_j$, and remember this is constant across all $i \neq j$. Then $P_1(i, i) = p_i = 1 - \sum_{j \neq i} c_j$. Let $M$ be the matrix with column $j$ constant $c_j$, so the rows are constant $[c_1, \dots, c_n]$ and let $r = \sum c_j$ be the constant row sum. Note we can write our transition matrix as
\[
T = (1 - r)I + M.
\] 
Note that the diagonal entries of this matrix match what we want since $T(i, i) = (1 - r) + c_i = 1 - \sum_{j \neq i} c_j = p_i$. 

Since $M$ has constant columns/rows it has null-space of dimension $n-1$. That is, we can find $n-1$ vectors $v_2, v_3, \dots v_n$ such that $vM = 0$. Then 
\[
vT = v((1 - r)I + M) = (1-r) v 
\]
(and remember of course the other eigenvector is the stationary distribution, with eigenvalue $1$). This says that there is only one eigenvalue $\beta = 1-r$, with multiplicity $n -1$ -- eigenfunctions are basis of null space of $M$. The eigenvalue $\beta$ is
\begin{align*}
\beta = 1 - \sum_j \frac{(q - 1)^2q^{n + j -3}}{(q^n - 1)(q^{n-1} - 1)} = 1 - \frac{q^{n-2}(q - 1)}{q^{n-1} - 1} = \frac{q^{n-2} - 1}{q^{n-1} - 1}.
\end{align*}

Note that this is the largest eigenvalue of the full chain on $S_n$.
\end{proof}

\begin{remark}
This lumping comes from the embedding $S_{n-1} \subset S_n$ as all permutations which fix the first coordinate. Then the coset space $S_n / S_{n-1}$ consists of equivalence classes of permutations which have the same label in the first element. Similarly, we could consider the embedding $S_{n-2} \subset S_n$ as all permutations which fix the first two coordinates. This would induce a Markov chain on the space of $\{(a, b): 1 \le a, b, \le n, a \neq b \}$; we have not attempted to find the more complicated transitions of this chain. 
\end{remark}

\subsection{Mixing Time}

The mixing time of the Markov chain $P_1$ is very fast; the chain reaches stationarity within a constant number of steps (the constant depends on $q$ but not on $n$). This can be proven using a strong stationary time.

\paragraph{Strong Stationary Time} 

An \emph{strong stationary time} for a Markov chain is a random stopping time $\tau$ at which the chain is distributed according to the stationary distribution. That is, if $(X_t)_{t \ge 0}$ is a Markov chain on $\Omega$ with stationary distribution $\pi$, then $\tau$ satisfies
\[
\P \left( \tau = t, X_\tau = y  \mid X_0 = x\right) = \P \left( \tau = t \mid X_0 = x \right) \pi(y), \quad x,y \in \Omega.
\]
In words, $X_\tau$ has distribution $\pi$ and is independent from $\tau$. See Section 6.4 of \cite{levinPeres}.

A strong stationary time is very powerful for bounding convergence time. Intuitively, once the strong stationary time is reached the chain has mixed, so the mixing time is bounded by the random time $\tau$. This idea is formalized in the following.

\begin{proposition}[Proposition 6.10 from \cite{levinPeres}]
If $\tau$ is a strong stationary time for a Markov chain $P$ on state space $\Omega$, then for any time $t$,
\[
\max_{x \in \Omega} \|P^t(x, \cdot) - \pi \|_{TV} \le \max_{x \in \Omega} \P \left( \tau > t \mid X_0 = x \right).
\]
\end{proposition}

For most Markov chains it is very difficult to find an obvious strong stationary time. A simple one exists for the chain $P_1$, using the alternate description from Remark \ref{remark: pq1}. 

To restate the alternate description, let $p := ((q - 1) q^{n-2})/(q^{n-1} - 1)$. Let $\{R_t \}_{t \ge 1}$ be a sequence Bernoulli$(p)$ random variables. The Markov chain $\{X_t \}_{t \ge 0}$ defined by $P_1$ can be coupled with the random variables $R_t$ by: Given $X_t$,
\begin{enumerate}
    \item If $R_{t+1} = 1$, sample $Z \sim \pi_{q, 1}$ and set $X_{t+1} = Z$.
    \item If $R_{t+1} = 0$, set $X_{t+1} = X_t$.
\end{enumerate}

\begin{lemma} \label{lem: SSTR}
With $\{R_t \}_{t \ge 1}$ defined above, the random time $
\tau = \inf \{ t > 0 : R_t = 1 \}. 
$
is a strong stationary time for $\{X_t \}_{t \ge 0}$. If \[
t =  c \cdot \frac{q^{n-1} - 1}{(q-1)q^{n-2}}, \quad c > 0, 
\]
then
\[
\P(\tau > t) < e^{-c}.
\]
\end{lemma}

\begin{proof}
By the alternate distribution of the Markov chain $P_1$, whenever $R_t = 1$, the next state $X_t$ is a sample from $\pi_{q, 1}$. That is,
\[
X_t \mid \{R_t = 1 \} \sim \pi_{q, 1}.
\]
Clearly then
\[
\P(\tau = t, X_\tau = y \mid X_0 = x) = \P(\tau =t \mid X_0 = x) \pi_{q, 1}(y),
\]
and $\tau$ is a strong stationary time. The time $\tau$ is a Geometric random variable with parameter $p$. Note also that $\tau$ is independent of the starting state $X_0$. Then,
\begin{align*}
    \P( \tau > t) &= (1-p)^{t} \le e^{-pt} \\
    &= \exp \left( -t \frac{(q - 1) q^{n-2}}{(q^{n-1} - 1)} \right) \le e^{-c},
\end{align*}
for $ t = c \cdot \frac{q^{n-1} - 1}{(q-1)q^{n-2}}$.
\end{proof}


\section{Some Extensions} \label{sec: extensions}

The main example treated above has $G = GL_n(q)$ and $H = K$ the Borel subgroup. As explained in Theorem \ref{thm: general} for \emph{any} finite group, \emph{any} subgroups $H, K$, and \emph{any} $H$-conjugacy invariant probability measure $Q$ on $G$ ($Q(hsh^{-1}) = Q(s)$ for all $s \in G, h \in H$), the walk on $G$ generated by $Q$, lumped to double cosets $H \backslash G / K =: \mathcal{X}$ gives a Markov chain on $\mathcal{X}$ with transition kernel
\[
K(x, A) = Q(HAKx^{-1}), \quad A \subset \mathcal{X},
\]
and stationary distribution the image of Haar measure on $G$.

There are many possible choices of $G, H, K,$ and $Q$. This gives rise to the problem of making choices and finding interpretations that will be of interest. This section briefly describes a few examples: Gelfand pairs, contingency tables, the extension from $GL_n(q)$ to finite Chevally groups, and a continuous example $\mathcal{O}_{n-1} \backslash \mathcal{O}_n / \mathcal{O}_{n-1}$. We have high hopes that further interesting examples will emerge.

\subsection{Parabolic Subgroups of $GL_n$}
In, \cite{karpThomas} the authors enumerate the double cosets of $GL_n(q)$ generated by parabolic subgroups.

\begin{definition}
Let $\alpha = (\alpha_1, \dots, \alpha_k)$ be a partition of $n$. The parabolic subgroup $P_\alpha \subset GL_n(q)$ consists of all invertible block upper-triangular matrices with diagonal block sizes $\alpha_1, \dots, \alpha_k$
\end{definition}

Section 4 of \cite{karpThomas} shows that if $\alpha, \beta$ are two partitions of $n$, the double cosets $P_\alpha \backslash GL_n(q) / P_\beta$ are indexed by contingency tables with row sums $\alpha$ and column sums $\beta$. Proposition 4.37 contains the size of a double-coset which corresponds to the table $T$, with $\theta = 1/q$
\begin{equation} \label{eqn: parabolicGLn}
\theta^{- n^2 + \sum_{1 \le i < i' \le I, 1 \le j < j' \le J} T_{ij} T_{i'j'}}(1 - \theta)^n \cdot \frac{\prod_{i = 1}^I[\alpha_i]_\theta! \prod_{j = 1}^J [\beta_j]_\theta!}{\prod_{i, j} [M_{ij}]_\theta!}, \quad I = |\alpha|, J = |\beta|.
\end{equation}
Dividing \eqref{eqn: parabolicGLn} by $(1-\theta)^n$ and sending $\theta \to 1$ recovers the usual Fisher-Yates distribution for partitions $\alpha, \beta$.

For $\alpha = (n-1, 1), \beta = (1^n)$, the contingency tables $P_\alpha \backslash GL_n(q) / P_\beta$ are uniquely determined by the position of the single $1$ in the second row. That is, the space is in bijection with the set $\{1, 2, \dots, n \}$. If $T^j$ denotes a table with the $1$ in the second row in column $j$, i.e.\ $T_{2j}^j = 1, T_{2k}^j = 0$ for all $k \neq j$, then Equation \eqref{eqn: parabolicGLn} becomes
\begin{align*}
    q^{-n + n^2 - (j-1)}(q -1)^n \cdot [n-1]_{1/q}! = q^{\binom{n}{2} - (j-1)} \frac{q^n - q^{n-1}}{q^n - 1} \cdot  \prod_{k = 1}^n (1 - q^k),
\end{align*}
which uses
\begin{align*}
    [n-1]_{1/q}! &= \frac{((1/q)^{n-1} - 1)((1/q)^{n-2} - 1) \hdots (1/q - 1)}{(1/q - 1)^{n-1}} = \prod_{i = 1}^{n-1} \frac{q^i - 1}{q^i - q^{i-1}} \\
    &= \frac{q^n - q^{n-1}}{q^n - 1} \cdot \prod_{i = 1}^{n} \frac{q^i - 1}{q^i - q^{i-1}} = \frac{q^n - q^{n-1}}{q^n - 1} \cdot \frac{1}{q^{\binom{n}{2}}(q - 1)^n} \cdot \prod_{k = 1}^n (1 - q^k).
\end{align*}
Dividing by $|GL_n(q)| = (q - 1)^n q^{\binom{n}{2}} \cdot [n]_q$ gives
\begin{align*}
    q^{-(j-1)} \cdot  \frac{(q^n - q^{n-1})(q - 1)}{(q^n - 1)^2(q - 1)^n} \prod_{k = 1}^{n} (q^k - 1).
\end{align*}
We initially thought that this distribution would match $\pi_{q, 1}$, and that the `follow the top card' chain in Section \ref{sec: topCard} was equivalent to the induced chain on double cosets $P_\alpha \backslash GL_n(q) / P_\beta$ from random transvections on $GL_n(q)$. However, this does not appear to be the case; $\pi_{q, 1}(j)$ is proportional to $q^{(j-1)}$, whereas the distribution induced on $\{1, 2, \dots, n \}$ via $P_\alpha \backslash GL_n(q) / P_\beta$ gives probability of $j$ proportional to $q^{(j-1)}$. It remains an open problem to investigate probability distributions and Markov chains on the double cosets of $GL_n(q)$ from parabolic subgroups. In a reasonable sense, for finite groups of Lie type, parabolic subgroups or close cousins are the only systematic families which can occur; see \cite{seitz}.

\subsection{Gelfand Pairs}

A group $G$ with subgroup $H$ is a \emph{Gelfand pair} if the convolution of $H$-biinvariant functions is commutative. Probability theory for Gelfand pairs was initiated by Letac \cite{letac} and Sawyer \cite{sawyer} who studied the induced chain on
\[
GL_n(\mathbb{Z}_p) \backslash GL_n(\mathbb{Q}_p) / GL_n(\mathbb{Z}_p)
\]
as simple random walk on a $p$-ary tree. Many further examples of finite Gelfand pairs appear in \cite{diaconisShah1987}, \cite{diaconisBook}, \cite{ceccherini}. These allow analysis of classical problems such as the Bernoulli-Laplace Model of diffusion and natural walks on phylogenetic trees. The literature cited above contains a large number of concrete examples waiting to be interpreted. For surprising examples relating Gelfand pairs, conjugacy class walks on a `dual group', and `folding', see \cite{lawther}.

We also note that Gelfand pairs occur more generally for compact and non-compact groups. For example, $\mathcal{O}_n/\mathcal{O}_{n-1}$ is Gelfand and the spherical functions become the spherical harmonics of classical physics (this example is further discussed in Section \ref{sec: continuous} below).  Gelfand pairs are even useful for large groups such as $S_\infty$ and $U_\infty$, which are not locally compact; see \cite{borodinOlshanki} and \cite{neretin} for resesarch in this direction.


\subsection{Contingency Tables}

Simper \cite{simperThesis} considers the symmetric group $S_n$ with parabolic subgroups $S_\lambda$ and $S_\mu$, for $\lambda = (\lambda_1, \dots, \lambda_I), \mu = (\mu_1, \dots, \mu_J)$ partitions of $n$. Then $S_\lambda \backslash S_n / S_\mu$ can be interpreted as $I \times J$ contingency tables $\{T_{ij} \}_{1 \le i \le I, 1 \le j \le J}$ with row sums $\lambda$ and column sums $\mu$. Such tables appear in every kind of applied statistical work. See \cite{diaconisSimper} Section 5 for references. The stationary distribution,
\[
\pi(T) = \frac{1}{n!} \prod_{i, j} \frac{\lambda_i! \mu_j!}{T_{ij}!},
\]
is the familiar Fisher-Yates distribution underlying `Fisher's exact test for independence'. Markov chains with the Fisher-Yates distribution as stationary were studied in \cite{diaconisSturmfels}. If $Q$ is the random transposition measure on $S_n$, \cite{simperThesis} uses knowledge of the $Q$ chain to give an eigen-analysis of the chain induced by $Q$ on contingency tables.

\subsection{Chevalley Groups}
Let $G$ be a finite Chevalley group (a finite simple group of Lie type). These come equipped with natural notions of Borel subgroups $B$ and Weyl groups $W$. The Bruhat decomposition
\[
G = \bigsqcup_{\omega \in W} B \omega B
\]
is in force, and conjugacy invariant probabilities $Q$ on $G$ induce Markov chains on $W$. Let $U$ be a minimal unipotent conjugacy class in $G$ (\cite{carterBook}, Chapter 5) and $Q$ the uniform probability on $U$. Conjugacy invariance implies that convolving by $Q$ induces an element of $\End_G(G/B)$. This may be indentified with the Hecke algebra of $B$-biinvariant functions. James Parkinson has shown us that the argument of Section \ref{sec: transvections} (for the computation
of $D$ in the Hecke algebra) goes through for a general Chevalley groups $G$ over a finite field $\FF_q$. Although a similar formula holds in full generality, for simplicity
we will state it in the equal parameter case (non twisted) and 
when $G$ is not of type $C_n$ with character lattice equal to the root lattice.
Let $\theta$ be the highest root and let $\cX_\theta = \{ x_\theta(c)\ |\ c\in \FF_q\}$
be the corresponding root subgroup.  The conjugacy class of $x_\theta(1)$ is 
an analogue, for general Chevalley groups, of the conjugacy class of transvections
in $GL_n(q)$.
The sum of the elements in the conjugacy class of $x_\theta(1)$ provides
a Markov chain on $B\backslash G/B$ which acts the same way as
\begin{align*}
D 
&= (q-1)\sum_{\alpha\in \Phi^+_l} q^{\frac12(I(s_\theta)-I(s_\alpha))}(1+T_{s_\alpha}),
\end{align*}
where $\theta$ is the highest root, $\Phi^+_l$ is the set of positive long roots,
and $s_\beta$ denotes the reflection in the root $\beta$ and $T_w$ is the element of the
Iwahori-Hecke algebra for $B\backslash G/B$ corresponding to the element $w$ in the
Weyl group.

Following the ideas in \cite{diaconisRam}, the $T_{s_\alpha}$ can be interpreted via the Metropolis algorithm applied to the problem of sampling from the stationary distribution $\pi(x) = Z^{-1}q^{I(x)}$ on $W$ by choosing random generators. We have not worked out any further examples but would be pleased if someone would.

\subsection{A Continuous Example} \label{sec: continuous}

Most of the generalities above extend to compact groups $G$ and closed subgroups $H, K \neq G$. Then, $X = H \backslash G /K$ is a compact space and an $H$-conjugacy invariant probability $Q$ on $G$ induces a Markov chain on $X$.

To consider a simple example, let $G = \O_n$, the usual orthogonal group over $\R$ and $H = K = \O_{n-1}$ embedded as all $m \in \O_n$ fixing the `North pole' $e_1 = (1, 0, 0 \dots, 0)^\top$. Then, $\O_n/\O_{n-1}$ can be thought of as the $(n-1)$-sphere $\mathcal{S}_{n-1}$. The double coset space $\O_{n-1} \backslash \O_n / \O_{n-1}$ codes up the `latitude' -- Consider the sphere $\O_n / \O_{n-1}$ defined by `circles' orthogonal to $e_1$.

\begin{figure}[ht]
    \centering
    \includegraphics{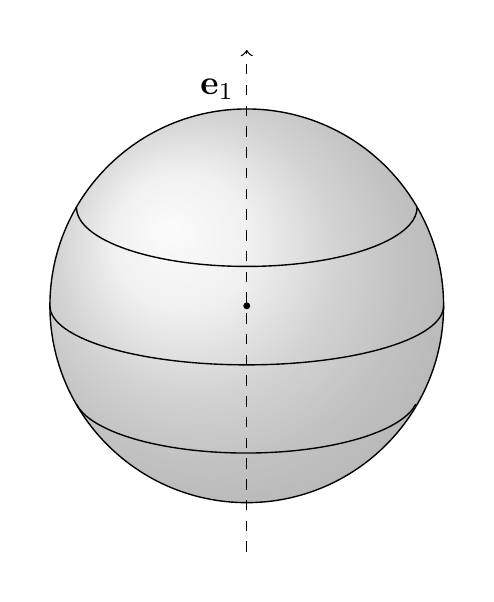}
    \caption{The space $\O_n/\O_{n-1}$ is defined by the circles on the sphere orthogonal to $e_1$.}
    \label{fig:sphere}
\end{figure}

Then $\O_{n-1} \backslash \O_n / \O_{n-1}$ simply codes which circle contains a given point on the sphere. Thus, $\O_{n-1} \backslash \O_n / \O_{n-1}$ may be identified with $[-1, 1]$.

Represent a uniformly chosen point on the sphere as $x = z/\|z \|$ with $z = (z_1, z_2, \dots, z_n)$ independent standard normals. The latitude is $z_1/\|z\|$ and so $\pi(x)$ is the distribution of the square root of a $\beta(1/2, (n-1)/2)$ distribution on $[-1, 1]$. When $n = 3$, $\pi(x)$ is uniform on $[-1, 1]$ (Theorem of Archimedes).

One simple choice for a driving measure $Q$ on $\O_n$ is `random reflections'. In probabilistic language, this is the distribution of $I - 2 U U^\top$, with $U$ uniform on $\mathcal{S}_{n-1}$. There is a nice probabilistic description of the induced walk on the sphere.

\begin{lemma} \label{lem: sphere}
The random reflections measure on $\mathcal{S}_{n-1}$ has the following equivalent description:
\begin{itemize}
    \item From $x \in \mathcal{S}_{n-1}$, pick a line $\ell$ through $x$ uniformly. With probability $1$, $\ell$ intersects the sphere in a unique point $y$. Move to $y$.
    
\end{itemize}
\end{lemma}

\begin{figure}[ht]
    \centering
    \includegraphics{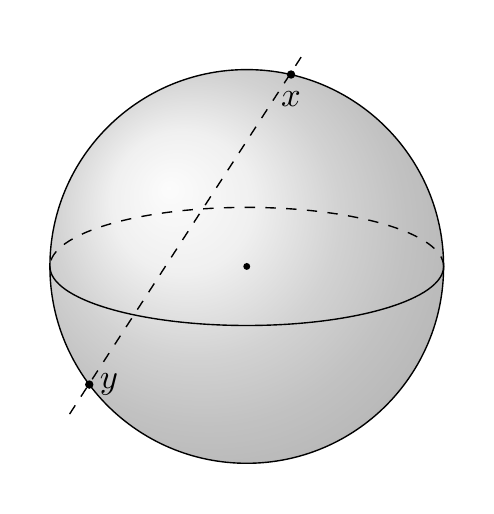}
    \caption{Illustration of procedure from Lemma \ref{lem: sphere}.}
\end{figure}

\begin{remark}
As the lemma shows, there is a close connection between the walk generated from $Q$ and the popular `princess and monster' algorithm. See \cite{comets}. These algorithms proceed in general convex domains. We know all the eigenvalues of the walk on the sphere and can give sharp rates of convergence.
\end{remark}

The induced chain on $[-1, 1] \cong \O_{n-1} \backslash \O_n / \O_{n-1}$ is obtained by simply reporting the latitude of $y$. Thanks to Sourav Chatterjee for the following probabilistic description. For simplicity, it is given here for $n = 3$ (so $\pi(x)$ is uniform on $[-1, 1]$).

\begin{lemma}
For $n = 3$, the Markov chain on $[-1, 1]$ described above is; From $x_0 \in [-1, 1]$, the chain jumps to
\[
X_1 := (1 - 2U_1^2)x_0 + 2(\cos(\pi U_2)) |U_1| \sqrt{ 1 - U_1^2} \sqrt{1 - x_0^2}
\]
where $U_1, U_2$ are i.i.d.\ uniform on $[-1, 1]$ random variables. (We have checked that the uniform distribution is stationary using Monte Carlo.)
\end{lemma}

\begin{remark}
For a detailed analysis of the random reflections walk on $\O_n$ (including all the eigenvalues), see \cite{porod}.
\end{remark}

\bibliographystyle{abbrv}
\bibliography{citations}

\appendix

\section{Row reduction for $GL_n(q)$} \label{app: rowReduction}

Explicit row reduction can be performed on a general transvection to determine which $\cB-\cB$ double coset it lies in. This allows computing transition probabilities using the combinatorics of possible transvections, which is elementary but tedious. 

\begin{proposition}
Let $T_{\ma, \mv} = I + \mv \ma^\top$ be the transvection defined by non-zero vectors $\ma, \mv$ with $\mv^\top \ma = 0$ and the last nonzero entry of $\mv$ equal to $1$. If $j > i$, then
$$T_{a,v} \in B s_{j-1}\cdots s_{i+1}s_is_{i+1}\cdots s_{j-1} B$$
exactly when the last nonzero entry of $\mv$ is $v_j$ and the first nonzero entry of $a$ is $\ma_i$.
\end{proposition}

The transvection corresponding to $\ma, \mv$ is
$$T_{\ma, \mv}(\mx) = x+\mv(\ma^\top \mx)
\qquad\hbox{so that}\qquad T_{\ma, \mv}(\mathbf{e}_i) = \mathbf{e}_i+a_i \mv,
$$
and the $i$th column of $T_{\ma, \mv}$ is $a_i \mv$ except with an extra $1$ added to the $i$th entry.
So
$$T_{\ma, \mv} = \begin{pmatrix} 1+a_1v_1 &a_2v_1 &a_3v_1 &\cdots &a_nv_1 \\
a_1v_2 &1+a_2v_2 &a_3v_2 &\cdots &a_nv_2 \\
\vdots &&&&\vdots \\
a_1v_n &a_2v_n &&\cdots &1+a_nv_n
\end{pmatrix}
=1 + \big( a_iv_j \big)_{1\le i,j\le n}.
$$
As an example of the row reduction, take $n = 5$ with $v_5 = 1, a_1 \neq 0$. Then,
\begin{align*}
T_{\ma, \mv} &= \begin{pmatrix}
1+a_1v_1 &a_2v_1 &a_3v_1 &a_4v_1 &a_5v_1 \\
a_1v_2 &1+a_2v_2 &a_3v_2 &a_4v_2 &a_5v_2 \\
a_1v_3 &a_2v_3 &1+a_3v_3 &a_4v_3 &a_5v_3 \\
a_1v_4 &a_2v_4 &a_3v_4 &1+a_4v_4 &a_5v_4 \\
a_1v_5 &a_2v_5 &a_3v_5 &a_4v_5 &1+a_5v_5 
\end{pmatrix}
\end{align*}
\begin{align*}
&= y_4(v_4v_5^{-1}) 
\begin{pmatrix}
1+a_1v_1 &a_2v_1 &a_3v_1 &a_4v_1 &a_5v_1 \\
a_1v_2 &1+a_2v_2 &a_3v_2 &a_4v_2 &a_5v_2 \\
a_1v_3 &a_2v_3 &1+a_3v_3 &a_4v_3 &a_5v_3 \\
a_1v_5 &a_2v_5 &a_3v_5 &a_4v_5 &1+a_5v_5 \\
0 &0 &0 &1 &-v_4v_5^{-1}
\end{pmatrix}
\\
&= y_4(v_4v_5^{-1}) y_3(v_3v_5^{-1}) 
\begin{pmatrix}
1+a_1v_1 &a_2v_1 &a_3v_1 &a_4v_1 &a_5v_1 \\
a_1v_2 &1+a_2v_2 &a_3v_2 &a_4v_2 &a_5v_2 \\
a_1v_5 &a_2v_5 &a_3v_5 &a_4v_5 &1+a_5v_5 \\
0 &0 &1 &0 &-v_3v_5^{-1} \\
0 &0 &0 &1 &-v_4v_5^{-1}
\end{pmatrix}
\\
&= y_4(v_4v_5^{-1}) y_3(v_3v_5^{-1}) y_2(v_2v_5^{-1}) 
\begin{pmatrix}
1+a_1v_1 &a_2v_1 &a_3v_1 &a_4v_1 &a_5v_1 \\
a_1v_5 &a_2v_5 &a_3v_5 &a_4v_5 &1+a_5v_5 \\
0 &1 &0 &0 &-v_2v_5^{-1} \\
0 &0 &1 &0 &-v_3v_5^{-1} \\
0 &0 &0 &1 &-v_4v_5^{-1}
\end{pmatrix}
\\
&= y_4(v_4v_5^{-1}) y_3(v_3v_5^{-1}) y_2(v_2v_5^{-1}) y_1(a_1^{-1}v_5^{-1}+v_1v_5^{-1})
\begin{pmatrix}
a_1v_5 &a_2v_5 &a_3v_5 &a_4v_5 &1+a_5v_5 \\
0 &-a_1^{-1}a_2 &-a_1^{-1}a_3 &-a_1^{-1}a_4 &z \\
0 &1 &0 &0 &-v_2v_5^{-1} \\
0 &0 &1 &0 &-v_3v_5^{-1} \\
0 &0 &0 &1 &-v_4v_5^{-1}
\end{pmatrix}
\end{align*}
with $(a_1^{-1}v_5^{-1}+v_1v_5^{-1})(1+a_5v_5)+z = a_5v_1$ so that
\begin{align*}
z
&= a_5v_1 - (a_1^{-1}v_5^{-1}+a_1^{-1}a_5+v_1v_5^{-1}+a_5v_1) \\
&=  - a_1^{-1}a_5 - a_1^{-1}v_5^{-1}  -v_1v_5^{-1}.
\end{align*}

Thus
\begin{align*}
T_{\ma, \mv} 
&= y_4(v_4v_5^{-1}) y_3(v_3v_5^{-1}) y_2(v_2v_5^{-1}) y_1(a_1^{-1}v_5^{-1}+v_1v_5^{-1}) \\
&\quad \cdot y_2(-a_2a_1^{-1})
\begin{pmatrix}
a_1v_5 &a_2v_5 &a_3v_5 &a_4v_5 &1+a_5v_5 \\
0 &1 &0 &0 &-v_2v_5^{-1}  \\
0 &0 &-a_3a_1^{-1} &-a_4a_1^{-1} &z-a_2v_2a_1^{-1}v_5^{-1} \\
0 &0 &1 &0 &-v_3v_5^{-1} \\
0 &0 &0 &1 &-v_4v_5^{-1}
\end{pmatrix}
\\
&= y_4(v_4v_5^{-1}) y_3(v_3v_5^{-1}) y_2(v_2v_5^{-1}) y_1(a_1^{-1}v_5^{-1}+v_1v_5^{-1}) \\
&\quad \cdot y_2(-a_2a_1^{-1})y_3(-a_3a_1^{-1})
\begin{pmatrix}
a_1v_5 &a_2v_5 &a_3v_5 &a_4v_5 &1+a_5v_5 \\
0 &1 &0 &0 &-v_2v_5^{-1}  \\
0 &0 &1 &0 &-v_3v_5^{-1} \\
0 &0 &0 &-a_4a_1^{-1} &z-(a_2v_2+a_3v_3)a_1^{-1}v_5^{-1} \\
0 &0 &0 &1 &-v_4v_5^{-1}
\end{pmatrix}
\\
&= y_4(v_4v_5^{-1}) y_3(v_3v_5^{-1}) y_2(v_2v_5^{-1}) y_1(a_1^{-1}v_5^{-1}+v_1v_5^{-1}) \\
&\quad \cdot y_2(-a_2a_1^{-1})y_3(-a_3a_1^{-1})y_4(-a_4a_1^{-1})
\begin{pmatrix}
a_1v_5 &a_2v_5 &a_3v_5 &a_4v_5 &1+a_5v_5 \\
0 &1 &0 &0 &-v_2v_5^{-1}  \\
0 &0 &1 &0 &-v_3v_5^{-1} \\
0 &0 &0 &1 &-v_4v_5^{-1} \\
0 &0 &0 &0 &z-(a_2v_2+a_3v_3+a_4v_4)a_1^{-1}v_5^{-1} 
\end{pmatrix}
\\
&= y_4(v_4v_5^{-1}) y_3(v_3v_5^{-1}) y_2(v_2v_5^{-1}) y_1(a_1^{-1}v_5^{-1}+v_1v_5^{-1}) \\
&\quad \cdot y_2(-a_2a_1^{-1})y_3(-a_3a_1^{-1})y_4(-a_4a_1^{-1})
\begin{pmatrix}
a_1v_5 &a_2v_5 &a_3v_5 &a_4v_5 &1+a_5v_5 \\
0 &1 &0 &0 &-v_2v_5^{-1}  \\
0 &0 &1 &0 &-v_3v_5^{-1} \\
0 &0 &0 &1 &-v_4v_5^{-1} \\
0 &0 &0 &0 &-a_1^{-1}v_5^{-1}
\end{pmatrix}
\end{align*}

since
\begin{align*}
z-(a_2&v_2+a_3v_3+a_4v_4)a_1^{-1}v_5^{-1} 
=(-a_5v_5-1-a_1v_1)a_1^{-1}v_5^{-1}-(a_2v_2+a_3v_3+a_4v_4)a_1^{-1}v_5^{-1} \\
&=-a_1^{-1}v_5^{-1}
\end{align*}
Thus, if $j>i$ then
$$T_{\ma, \mv} \in Bs_{j-1}\cdots s_{i+1}s_is_{i+1}\cdots s_{j-1} B$$
exactly when the last nonzero entry of $v$ is $v_j$ and the first nonzero entry of $a$ is $a_i$.

\end{document}